\documentclass[3p,times]{elsarticle}

\usepackage{graphicx}
\usepackage{psfrag}
\usepackage[dvipsnames]{xcolor}

\usepackage{amssymb}
\usepackage{amsmath}
\usepackage{dsfont}
\usepackage{rotating}
\usepackage{multirow}

\DeclareMathAlphabet\mathbfcal{OMS}{cmsy}{b}{n}  

\newcommand{\ta}{\tilde{a}}
\newcommand{\tw}{\tilde{w}}

\newcommand{\dt}{\Delta t}

\newcommand{\Q}{\mathbf{Q}} 
 
\newcommand{\f}{\mathbf{f}}

\newcommand{\B}{\mathbf{B}}

\newcommand{\x}{\mathbf{x}}

\DeclareMathOperator{\argmin}{argmin}

\def\R{{\cal R}}

\def\be{\begin{equation}}
\def\ee{\end{equation}}
\def\bea{\begin{eqnarray}}
\def\eea{\end{eqnarray}}
\def\ba{\begin{array}{l}\displaystyle}
\def\ea{\end{array}}
\def\E{{\cal E}}
\def\ca{\\[+0.3cm]\displaystyle}

\newcommand{\SSS}{\mathbb{S}}

\def\supp{{\rm Supp}}
\def\supp{\mbox{supp }}
\newcommand{\ens}[1]{\mathbb{#1}}
\def\Ball{ {\cal B}}

\def\f{\hat f}

\def\bb{\hat \beta}
\def\B{\hat B}

\begin{document}
	
	\begin{frontmatter}
		
		\journal{Computer methods in Applied Mechanics and Engineering} 
		
		
		
		\title{High order finite volume schemes with IMEX time stepping for the Boltzmann model on unstructured meshes} 
		
		\author[UNIFE,Centro]{Walter Boscheri}
		\ead{walter.boscheri@unife.it}
		\author[UNIFE,Centro]{Giacomo Dimarco}
		\ead{giacomo.dimarco@unife.it}
		\address[UNIFE]{Department of Mathematics and Computer Science, University of Ferrara, Via Machiavelli 30, 44121 Ferrara, Italy}
		\address[Centro]{Center for Modeling, Computing and Statistic CMCS, University of Ferrara, Via Muratori 9, 44121 Ferrara, Italy}

%
\begin{abstract}
In this work, we present a family of time and space high order finite volume schemes for the solution of the full Boltzmann equation. The velocity space is approximated by using a discrete ordinate approach while the collisional integral is solved by spectral methods. The space reconstruction is realized by integrating the distribution function, describing the state of the system, over arbitrary shaped and closed control volumes using a Central Weighted ENO (CWENO) technique. Compared to other reconstruction methods, this approach permits to keep compact stencil sizes which is a remarkable property in the context of kinetic equations due to the considerable demand of computational resources. The full discretization is then obtained by combining the previous phase-space approximation with high order Implicit-Explicit (IMEX) Runge Kutta schemes. These methods guarantee stability, accuracy and preservation of the asymptotic state. 
Comparisons of the Boltzmann model with simpler relaxation type kinetic models (like BGK) is proposed showing the capability of the Boltzmann equation to capture different physical solutions. The theoretical order of convergence is numerically measured in different regimes and the methods are tested on several standard two-dimensional benchmark problems in comparison with Direct Simulation Monte Carlo results. The article ends with a prototype engineering problem consisting of a subsonic and a supersonic flow around a NACA 0012 airfoil. All test cases are run with MPI parallelization on several threads, thus making the proposed methods suitable for parallel distributed memory supercomputers. 
\end{abstract}
%
\end{frontmatter}

%
%
\section{Introduction}
Kinetic models permit to describe a very large set of physical phenomena involving fluid and gases \cite{cercignani}. This description relies on a different point of view compared to standard fluid models such as compressible Euler or Navier-Stokes equations. In fact, the state of the system is determined through the time evolution of a distribution function giving the probability for a particle of the gas to be at a given position, with a given velocity, at a fixed instant of time. Consequently, the full problem depends, in general, on seven independent variables: three for the space, three for the velocity plus the time. This gives a very large spectrum of possibilities regarding the modeling: typically rarefied gases and high speed flows can be described only through these approaches while fluid models fail. The counterpart is however the requirement of very expensive simulations and very large variable storages compared to the compressible fluid dynamic case \cite{ACTA}.

In realistic simulations, the problem of the excessive computational cost is traditionally bypassed using Direct Simulation Monte Carlo methods \cite{Cf, Nanbu80} which are both versatile and very efficient at least for describing steady problems. However, the simulation of unsteady flows remains a real challenge for probabilistic methods due to the large statistical error introduced in the computation. Concerning deterministic methods, along with finite volume approaches \cite{Titarev1,Titarev2,Titarev3}, a recent strategy, based on a Lagrangian technique, has been shown to be able to deal efficiently with multidimensional kinetic equations \cite{FKS,FKS_HO,FKS_GPU}. This method has been also employed to compute the solution of a full six dimensional Boltzmann model \cite{FKS_Bolt}. Nevertheless, even if very efficient, this approach is only first order in time and space and it is based on a regular Cartesian mesh. 

In this work, we instead deal with high order accurate space-time finite volume methods and unstructured meshes. The methods here discussed are clearly more demanding in terms of computational resources compared to \cite{FKS} but, on the other hand, more accurate and versatile. In particular, the results presented in this work follows the approach recently forwarded in \cite{BosDim}, where a simple kinetic model of relaxation type is considered. Here, we deal with the much more challenging and physically realistic Boltzmann kinetic model \cite{cercignani}. We remark that on this subject, up to our knowledge, the number of articles dealing with multidimensional high order in space and time numerical methods for the full Boltzmann equation is extremely limited, in particular for high order methods on unstructured meshes we are aware of only few examples \cite{Hu3,Hu4,WU} that are all concerned with discontinuous Galerkin methods and explicit time discretizations or instead low order implicit ones. In this work we focus also on high order in time implicit methods which are of paramount importance in kinetic equations \cite{ACTA}.
 
To deal with the Boltzmann equation, we introduce four different level of discretization: two different type of discretization for the velocity space, discretization of the physical space and time discretization. For the velocity space, we first introduce a traditional discrete velocity model (DVM) \cite{bobylev} which replaces the continuous kinetic equation with a finite set of fixed velocities. This idea is then integrated by relying on a spectral approach for efficiently approximate the Boltzmann integral \cite{FiMoPa:2006,wu2013deterministic,FilbetRusso} as detailed later. The transport part of the DVM model is solved by a finite volume strategy, i.e. we evolve the cell averages of the distribution function. High order space accuracy requires a reconstruction that produces high order piecewise polynomials from the neighboring cell averages. Many reconstruction strategies do exist \cite{Shu1}, and we concentrate on a Central WENO reconstruction (CWENO) because of its compactness property. These have been originally introduced in the context of conservation laws in the one-dimensional case \cite{LPR:99} and also used in our recent work \cite{BosDim}. In particular, we develop these techniques on arbitrary polygonal meshes in the physical space. The CWENO reconstruction procedure has also been successfully used in a context of moving mesh schemes with grid regeneration \cite{Gaburro2017} and embedded in a quadrature-free finite volume scheme for solving all Mach compressible flows \cite{BosPar21}.

The fourth type of discretization, i.e. the time one, is handled by using Implicit-Explicit (IMEX) Runge-Kutta (RK) methods \cite{Dimarco_stiff2}. As detailed in \cite{ACTA}, kinetic equations are particularly difficult to solve due to their multiscale nature in which collision and transport scales coexist. Close to the fluid limit, the collision rate grows exponentially, while the fluid dynamic time scale conserves a much lower pace \cite{cercignani}. Several authors have tackled the above problem in the recent decades \cite{Jin2, Jin_review,Filbet_Jin, BLM,LemouAP, degondrev, DegondAP, Dimarco_stiff1, Dimarco_stiff3, Dimarco_stiff4, Qin,Hu1,Hu2, Puppo} by developing and studying the class of methods known as Asymptotic Preserving (AP). These techniques allow the full problem to be solved for time steps independent of the collisional fast scale identified by the Knudsen number. Moreover, they are stable and consistent with the limit model provided by the compressible Euler equation in the limit in which the collisional scale grows to infinity. 

A last difficulty is present when the full Boltzmann model is considered. One has to deal with the discretization of a five fold integral (for the three-dimensional case, or a three fold integral for the two-dimensional one) for each fixed value of the physical mesh. This is known to be a very challenging and computationally expensive problem when treated by deterministic techniques and it is one of the main reason why particles methods are often preferred in practice \cite{ACTA}. Here, this problem is overcome by Fourier techniques \cite{MoPa:2006,FiMoPa:2006} which can be shown to have a complexity of the order of $\mathcal O(N^{d_v} \log(N^{d_v}))$ where $N$ is the number of modes taken for the velocity space in one direction and $d_v$ the dimension of the velocity space. We recall that the research field in this direction is very active and many recent contributions are present in the literature (see \cite{wu2015fast,gamba:2010} and the references therein).

Summarizing, in this paper we extend the high order CWENO-IMEX Runge-Kutta Asymptotic Preserving methods on arbitrary shaped unstructured grids developed in \cite{BosDim} to the Boltzmann case. These methods are successively tested on benchmark multidimensional rarefied gas dynamic problems comparing relaxation type kinetic equations (like the simple BGK model) with the Boltzmann model. Comparisons of the proposed method with Direct Simulation Monte Carlo simulations applied to the solution of the Boltzmann model are also proposed showing an excellent agreement between the two different techniques. The theoretical accuracy is also numerically verified for different regimes. Finally, a more realistic test case involving a flow around a NACA 0012 airfoil for different values of the Knudsen number is performed in the subsonic and in the supersonic regime. A MPI parallelization is realized distributing the space variable on different threads and the computational cost related to the Boltzmann collision operator is discussed. 

The article is organized as follows. In Section \ref{sec_Boltzmann}, we introduce the Boltzmann and the BGK models, their properties and their fluid dynamic limit. In Section \ref{sec_Disc}, we present the four type of discretization, namely the discrete ordinate discretization, the Fourier approximation of the Boltzmann operator, the CWENO reconstruction and the IMEX Runge-Kutta schemes applied to the phase-space discretization of the kinetic equations. The Section \ref{sec.validation} is devoted to present several numerical examples of the schemes. In particular, up to third order accuracy in space and time is proved, the performances of the method are measured and its capability to deal with such equations is shown. Conclusions and future investigations are discussed in a final section.

%
%

\section{The Boltzmann and the BGK equations}
\label{sec_Boltzmann}
We consider kinetic equations belonging to the following class
\be
\frac{\partial f}{\partial t} + v\cdot \nabla_x f= Q(f).
\label{eq:B}
\ee
This equation supplemented by the initial condition $f(x,v,t=0)=f_{0}(x,v)$ furnishes the time evolution of a non-negative function $f=f(x,v,t)$ which gives the distribution of particles with velocity $v \in \R^{d_v} $ in the space $x \in \Omega \subset \R^{d_x}$ at time $ t > 0$. In the following, for simplicity, we will fix $d_x=d_v=d$, i.e. we consider the same dimension for both the physical and the velocity space. Moreover all numerical simulations are performed in the $d_x=d_v=2$ case which means in a four-dimensional space. The operator $Q(f)$ describes the effects of particle interactions and its form depends on the details of the microscopic dynamics. In particular, we focus on two models. The first one is a relaxation type model
\be
\partial_t f + v\cdot\nabla_{x}f = \frac{\nu}{\varepsilon} (M_{f}-f),
\label{eq:BGK}
\ee
which is known as the BGK equation \cite{Gross}. In this model, the complex interactions between particles are substituted by a relaxation towards the local thermodynamical equilibrium defined by the Maxwellian distribution function $M_{f}$ \be
M_f=M_{f}[\rho,u,T](v)=\frac{\rho}{(2\pi \theta)^{d/2}}\exp\left(\frac{-|u-v|^{2}}{2\theta}\right) ,\label{eq:M} \ee where $\rho \in \R, \ \rho>0$ and $u \in \R^d$ are respectively the density and mean velocity, while $\theta=RT$ with $T$ the temperature of the gas and $R$ the gas constant fixed to $R=1$ in the rest of the paper. These quantities are related to the distribution function through the following relations \be \rho=\int_{\R^d} fdv, \qquad u=\frac{1}{\rho}\int_{\R^d} vfdv, \qquad  \theta=\frac{1}{\rho d}
\int_{\R^d}|v-u|^{2}fdv. \label{eq:Mo} 
\ee Additionally, the parameter $\nu> 0$ in (\ref{eq:BGK}) is the relaxation frequency which is taken in this work
as $\nu=\rho$. Let observe that this is a simplified version of a more general frequency law of the type $\nu=\rho T^{1-\beta}$ (see \cite{Cow} for a discussion).

The second model we consider, which we focus our entire study on, is represented by the more challenging case of the full Boltzmann equation
\be
\partial_t f + v\cdot\nabla_{x}f =\frac{1}{\varepsilon} \int_{\R^{d_v}}
\int_{\SSS^{d_v-1}} B(|v-v_*|,\omega) \left( f(v')f(v'_*)-f(v)f(v_*) \right) dv_* d \omega.
\label{eq:Boltzmann}
\ee
In this model, $\omega$ is a vector of the unitary sphere $\SSS^{d_v-1} \subset \R^{d_v}$,  whereas $(v',v'_*)$ are given by the relations \be v'=\frac{1}{2}(v+v_*+|q|\omega),\quad
v'_*=\frac{1}{2}(v+v_*-|q|\omega), \ee
where $q=v-v_*$ is the
relative velocity and represents the post-collisional velocities. The kernel $B$ characterizes the
details of the binary interactions. Here we concentrate on the case of the so-called {Maxwell pseudo-molecules model} \cite{cercignani}, i.e. \be B(v,v_*,\omega)=
b_{0}(\cos\theta), \ee
which implies that the collision rate is independent with respect to the relative velocities between the particles in the gas. 

One of the key property of the collision operator which permits to establish the connection between kinetic and fluid models is that it guarantees the conservation of density, momentum and energy.
This fact can be resumed by writing
\be
\int_{\R^{d_v}} Q(f)\phi(v)\,dv=0,
\label{eq:i5}
\ee
where $\phi(v)=(1,v,|v|^2)$. Let observe that the energy is related to the other macroscopic quantities through
\be  \frac{1}{2}d_v\rho T=E-\frac{1}{2}\rho|u|^2.  \ee
The parameter $\varepsilon$ appearing in both models \eqref{eq:BGK} and \eqref{eq:Boltzmann} is the so-called Knudsen number which is used to rescale the equation in time and space to ease the transition from the different time scales which characterize a kinetic model. Since
the Boltzmann operator satisfies \be Q(f)=0\quad \iff \quad f=M_f, \label{eq:i6}\ee 
one can get in the limit $\varepsilon\rightarrow 0$ from (\ref{eq:Boltzmann}) or equivalently from \eqref{eq:BGK} the compressible Euler equations by replacing $f$ with $M_f$ and integrating in velocity space 
\be \ba \frac{\partial \rho}{\partial t} + \nabla_{x}
\cdot(\rho u) = 0, \ca \frac{\partial \rho u}{\partial t} +
\nabla_{x} \cdot (\rho u \otimes u+pI) = 0, \ca \frac{\partial
	E}{\partial t} +\nabla_{x} \cdot((E+p)u) = 0, \ca p=\rho \theta,
\quad E=\frac{d}{2}\rho \theta +\frac{1}{2} \rho |u|^{2}, \ea
\label{eq:sys1} \ee
where $p$ is the gas pressure. Using the expansion  $f=\sum_{n=0}^{\infty}\varepsilon^n f_n$ according to \cite{cercignani}, the Navier-Stokes equations can be derived as well. We stress that the schemes used in this work permits to recover the same property here observed for the continuous case: in the discrete case if $\varepsilon\to 0$ the numerical schemes derived in this work become automatically high order space-time discretizations of \eqref{eq:sys1}. In addition, these schemes enjoy the property of being stable independently of the size of the Knudsen number $\varepsilon$ (see \cite{ACTA, Dimarco_stiff1, Dimarco_stiff3} for details).

We need finally to define suitable boundary conditions in physical space in order to determine the solutions to kinetic equations of type \eqref{eq:B}. In this context, the only conditions that need to be discussed are the ones required when the fluid meets an object (the surface of a wing for instance) or equivalently hits a wall and particles interact with the atoms of the surface before being reflected backward. For $v \cdot n \ge 0$ and $x\in\partial \Omega$, where {$n$} denotes the unit normal pointing inside the domain {$\Omega$}, such boundary conditions are modelled by 
\begin{equation}
|v \cdot n| f(x,v,t)
= \int_{v_{\ast} \cdot n<0} |v_{\ast}\cdot n(x)| K(v_{\ast}
\rightarrow v, x,t) f(x,v_{\ast},t)\,dv_{\ast}. \label{eq:KE} 
\end{equation}
The above expression means that the ingoing flux is defined in terms of the outgoing flux modified by a given boundary kernel {$K$}. The only condition imposed is that positivity and mass conservation at the
boundaries are guaranteed. A general condition gives for the ingoing velocities the following relation
\be
f(x,v,t)=(1-\alpha)Rf(x,v,t)+\alpha Mf(x,v,t),
\label{eq:BOU}
\ee 
in which {$x \in \partial\Omega$, $v
	\cdot  n(x) \ge 0$}. Therefore a fraction $1-\alpha$ of the outgoing molecules are reflected while the ones belonging to the remaining part $\alpha$ are thermalized. The coefficient
{$\alpha$}, with {$0\leq \alpha\leq 1$}, is called the
{accommodation coefficient} and it holds \be Rf(x,v,t)=f(x,v-2 n( n
\cdot v),t),\quad  Mf(x,v,t)=\mu(x,t) M_w(v), \label{eq:MAXB} \ee
where $M_w(v)$ is a Maxwellian distribution with unit mass, fixed temperature and mean velocity corresponding to the speed of the object or of the wall. The value of {$\mu$} is instead determined by mass conservation \be \mu(x,t) \int_{v \cdot  n
	\geq 0}M_w(v)\vert v \cdot  n \vert dv = \int_{v \cdot  n <
	0}f(x,v,t)\vert v \cdot  n \vert dv. \label{eq:MU} \ee
In the numerical test involving boundaries, we will consider specifically the case {$\alpha=1$} corresponding to full accommodation in which the re-emitted molecules have completely lost memory of the incoming molecules.

%
%
\section{The numerical method}\label{sec_Disc}
This section is divided into four parts. The first one is about discrete velocity models, the second one is about the finite volume method, the third one discusses the spectral method for the the Boltzmann operator. The last one is about the time integration techniques.
\subsection{The Discrete Velocity Models (DVM)} \label{sec_DVM}
The unbounded velocity space is truncated and the tails of the distribution function are then not considered. The consequence is that the exact conservation of macroscopic quantities is impossible, because in general the support of the distribution function is non-compact as for instance for the Maxwellian equilibrium distribution $M_f$. The same holds true for the Boltzmann collision integral defined on the entire space $\mathcal{R}^d$. Different
strategies can be adopted to overcome the lack of conservation \cite{gamba, Mieussens}. Let observe that however if the truncation is performed with sufficiently large bounds, the loss of conservation is typically negligible, in fact distribution functions tend to zero exponentially fast. This is the direction pursued in this work. However, the rest of the scheme can be easily adapted to different choices.

This new bounded space is discretized with a finite number of points. These are called the discrete velocities, the only velocities the particles can assume. Taking inspiration from \cite{Mieussens}, we introduce a Cartesian grid $\mathcal{V}$  \be\label{disc_space}
\mathcal{V}=\left\{ v_{k}=k\Delta v+a, \ k=k^{(i)}, \ i=1,..,d, \ a=(a_1,..,a_d)\right\},\ee where $a$ is an arbitrary vector, $\Delta v$ is a constant mesh size in velocity and where the components of the index $k$ have some given bounds $K^{(i)}, \ i=1,..,d$. In this setting, the continuous distribution function $f$ is replaced by the vector $f_{\mathcal{K}}(x,t)$ of size $N$. Each component of this vector is assumed to be an approximation of the distribution function $f$ at location $v_{k}$: \be
f_{\mathcal{K}}(x,t)=(f_{k}(x,t))_{k},\qquad f_{k}(x,t) \approx
f(x,v_{k},t). \ee 
The discrete ordinate kinetic model consists then of the following system of equations to be solved 
\be
\partial_t f_{k} + v_{k} \cdot\nabla_{x}f_{k} = Q_k(f), \ k=1,..,N.
\label{eq:DM1_gen} \ee
For the two cases here considered we have either $ Q_k(f)=Q_{k,BGK}(f_k)\frac{\nu}{\varepsilon} (\E_{k}[U]-f_{k})$ for the BGK model, or $Q_k(f)=Q_{k,B}(f)$ for the Boltzmann model. In this case $Q_B(f_k)$ corresponds to the solution given by the spectral approximation of the collision integral \eqref{eq:Boltzmann} projected over the discrete space. Finally, the function $\E_{k}[U]$ represents a suitable approximation of $M_{f}$, e.g. $\E_{k}[U]=M_f(x,v_k,t)$, and $U=(\rho,\rho u ,E)^T$ is the vector of the macroscopic quantities. These quantities are recovered from the knowledge of the distribution function $f_{\mathcal{K}}(x,t)$ thanks to discrete summations on the discrete velocity space.

%
%
%

\subsection{Space discretization}
\label{sec_num_approx}
\begin{figure}[!htbp]
	\begin{center}
		\begin{tabular}{cc}
			\includegraphics[width=0.47\textwidth]{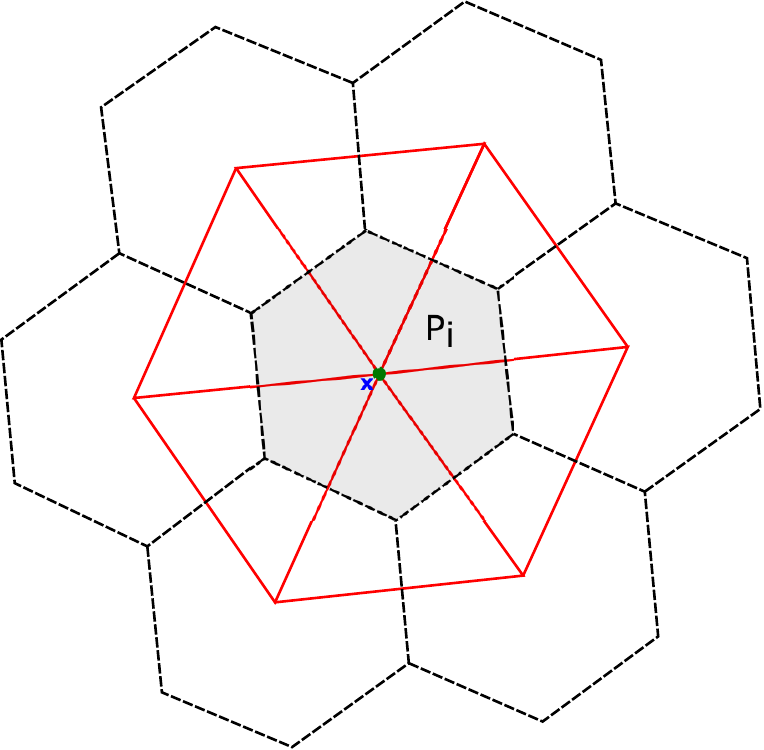}  &          
			\includegraphics[width=0.47\textwidth]{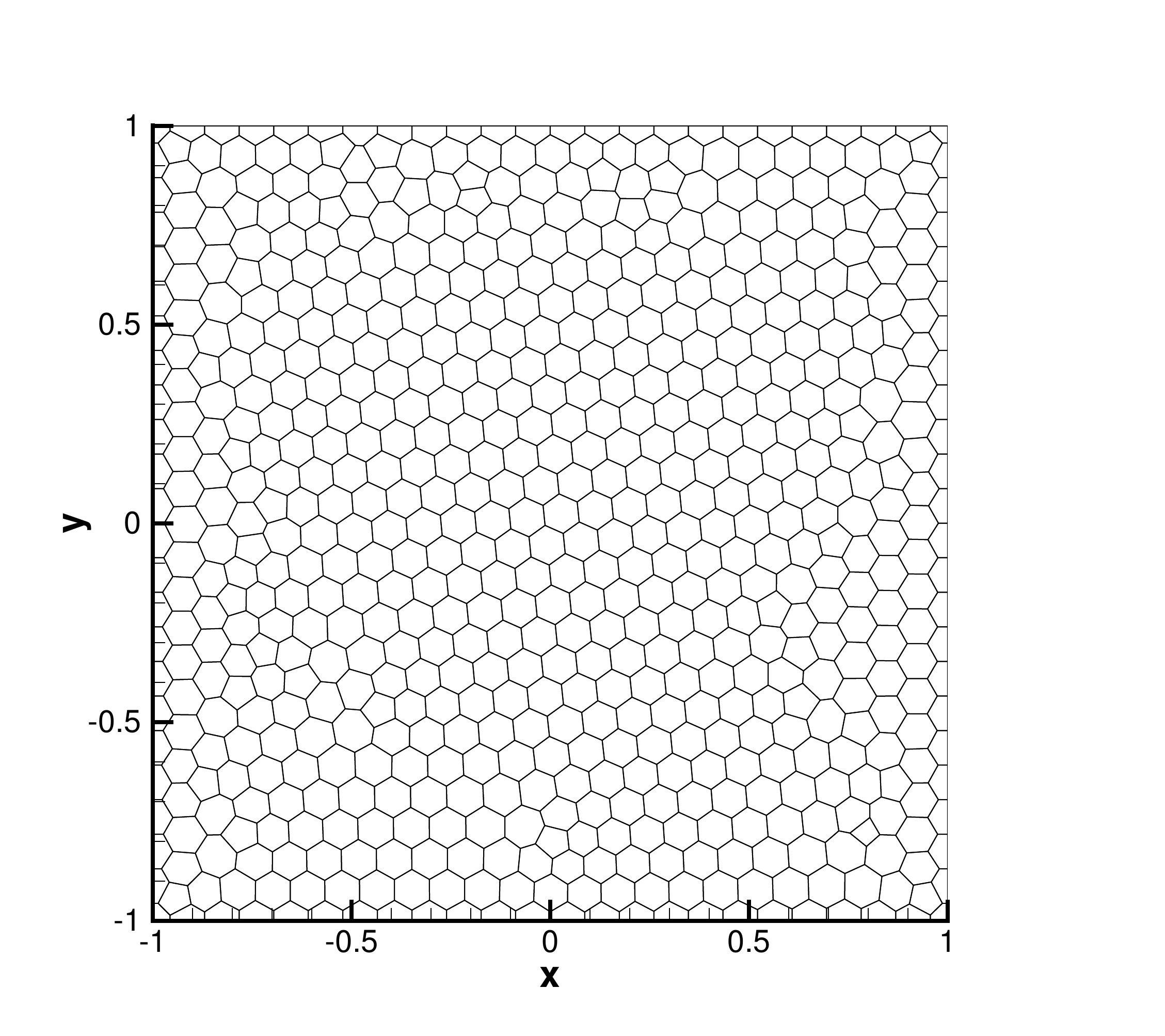} \\ 
		\end{tabular}
		\caption{Left: example of the construction of the Voronoi tessellation (black dotted line) starting from the Delaunay triangulation (red solid line). The generator point $x_{c,i}$ is marked with a blue cross, the barycenter $x_{b_i}$ of the Voronoi control volume is represented by a green circle. Right: example of an arbitrarily polygonal unstructured mesh.}
		\label{fig.grid}
	\end{center}
\end{figure}

We discuss now the space discretization of the discrete system of equations \eqref{eq:DM1_gen}. This is based on a finite volume framework and on a Central WENO (CWENO) reconstruction \cite{Boscheri_Russo, BosDim}. We start from the definition of the grid on domains $\Omega$ which are considered two-dimensional and arbitrarily shaped. In particular, we realize a centroid based Voronoi-type tessellation made of $N_P$ non overlapping polygons $P_i, i=1, \dots N_P$. The mesh is constructed by firstly fixing the position $x_{c_i}, i=1,\dots, N_P$ of the generator points and successively by introducing a Delaunay triangulation having these generators $x_{c_i}$ as vertexes of the triangles. Then each generator point $x_{c_i}$ is associated to a Voronoi element $P_i$ by connecting the centers of mass of all the Delaunay triangles having this generator point as a vertex. 
It is interesting to notice that, in this construction, the barycenter  $x_{b_i}$ of the obtained Voronoi polygon $P_i$
\be
x_{b_i} = \frac{1}{|P_i|} \int_{P_i} x \, dx,
\ee
usually differs with the position of the generator point, see Figure \ref{fig.grid} for a visual explanation. Let \be\mathcal{D}(P_i) = \{d_{i_1}, \dots, d_{i_j}, \dots, d_{i_{N_{V_i}}} \},\ee denote the set of $N_{V_i}$ vertexes of polygon $P_i$. Once the grid is set, we connect the barycenters $x_{b_i}$ with each vertex of $\mathcal{D}(P_i)$ and we subdivide the Voronoi polygon $P_i$ in $N_{V_i}$ triangles.  
This sub-triangulation is the one used in practice to numerically integrate the unknowns. Figure \ref{fig.grid} shows an example of such mesh. 

\begin{figure}[!htbp]
	\begin{center}
		\begin{tabular}{ccc} 
			\includegraphics[width=0.3\textwidth]{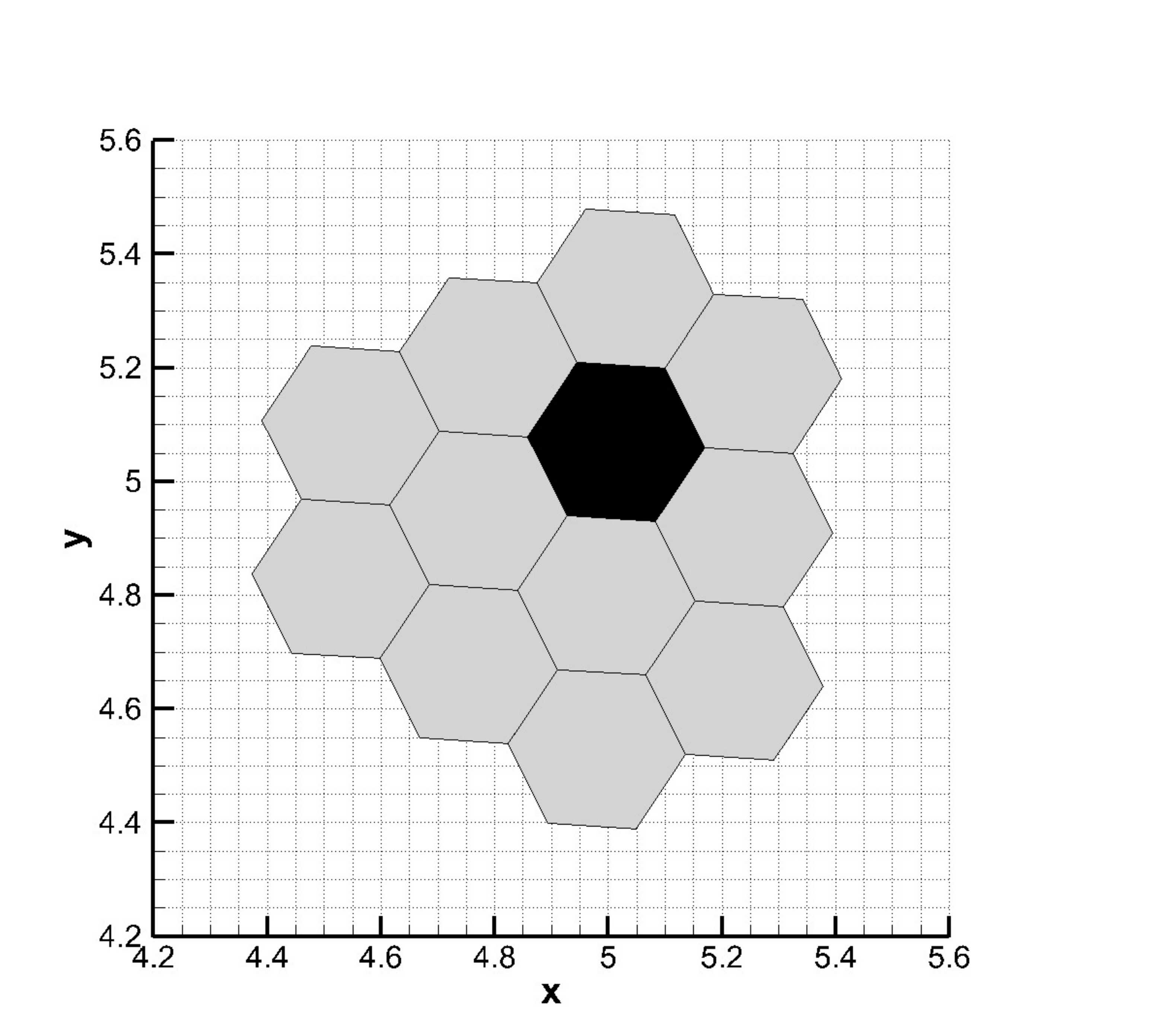}  &           
			\includegraphics[width=0.3\textwidth]{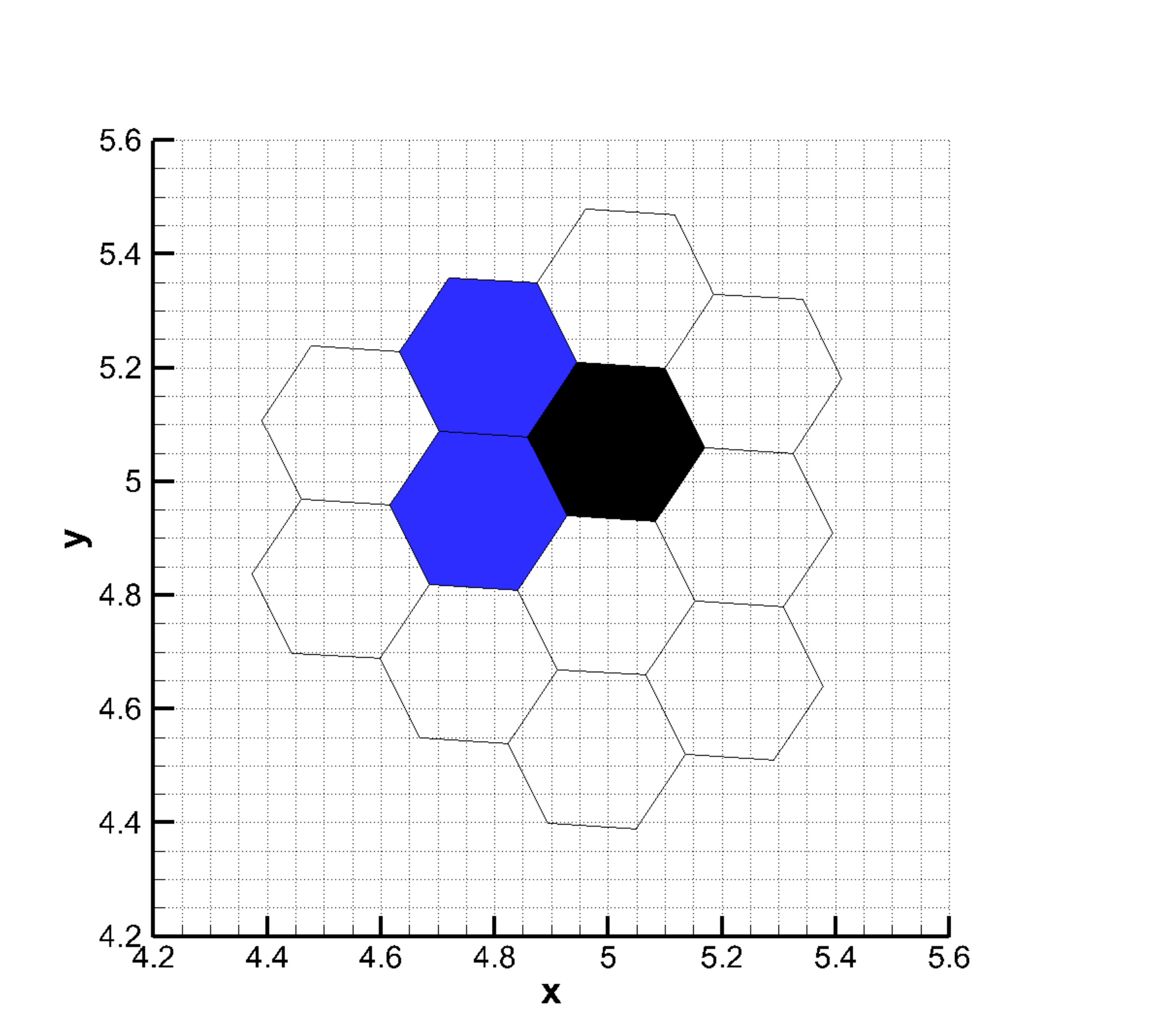} &
			\includegraphics[width=0.3\textwidth]{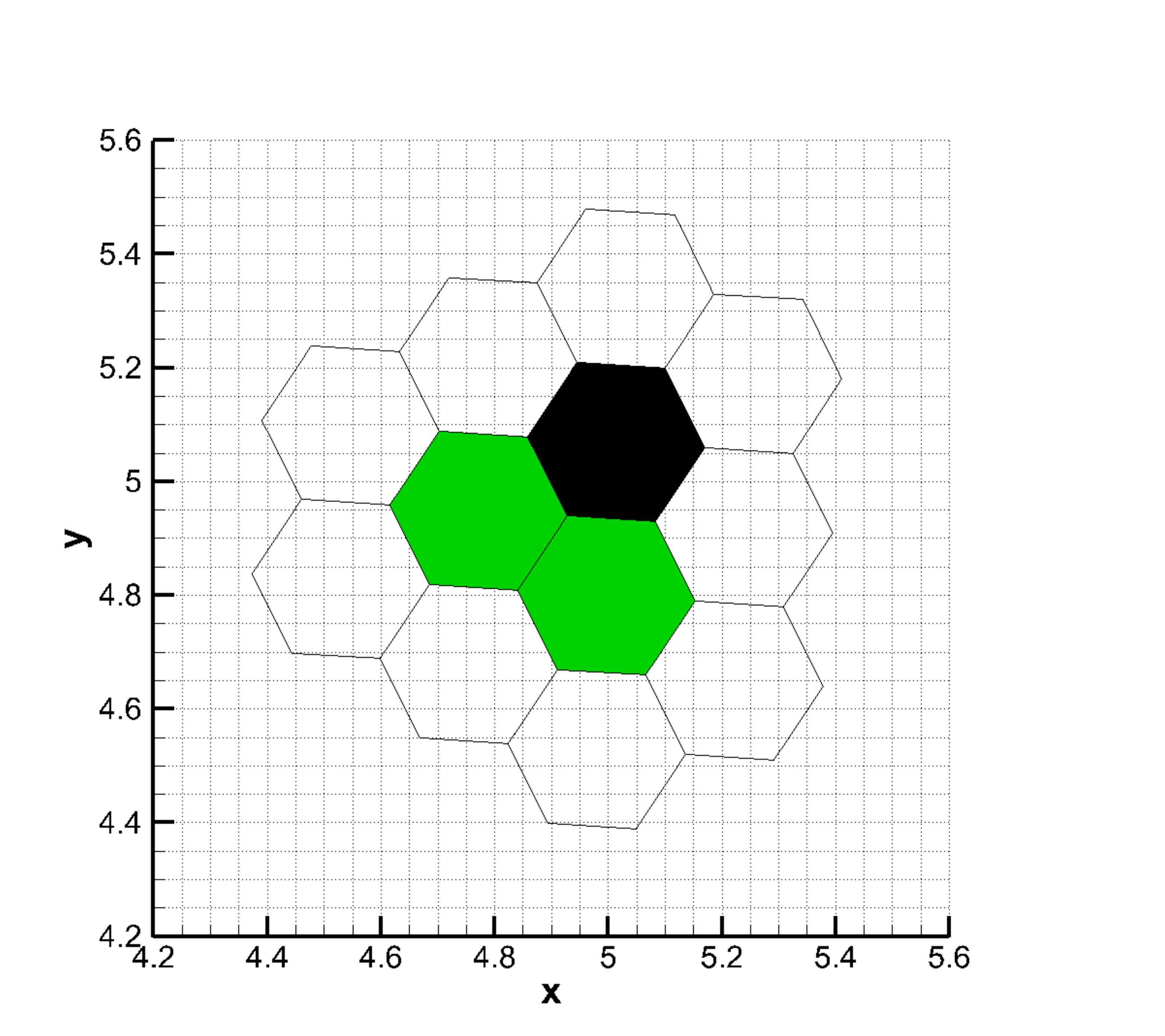}  \\           
			\includegraphics[width=0.3\textwidth]{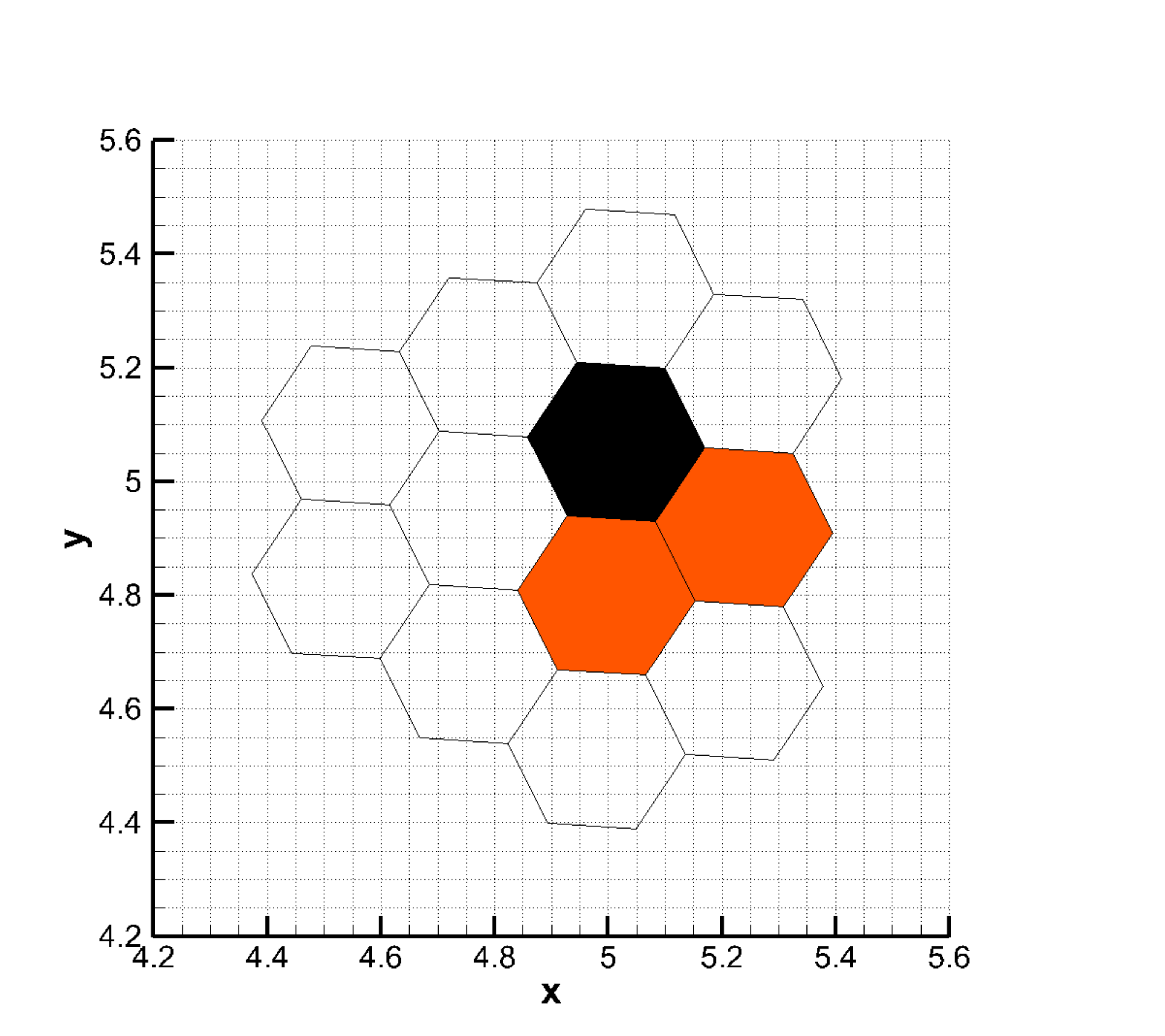}  &           
			\includegraphics[width=0.3\textwidth]{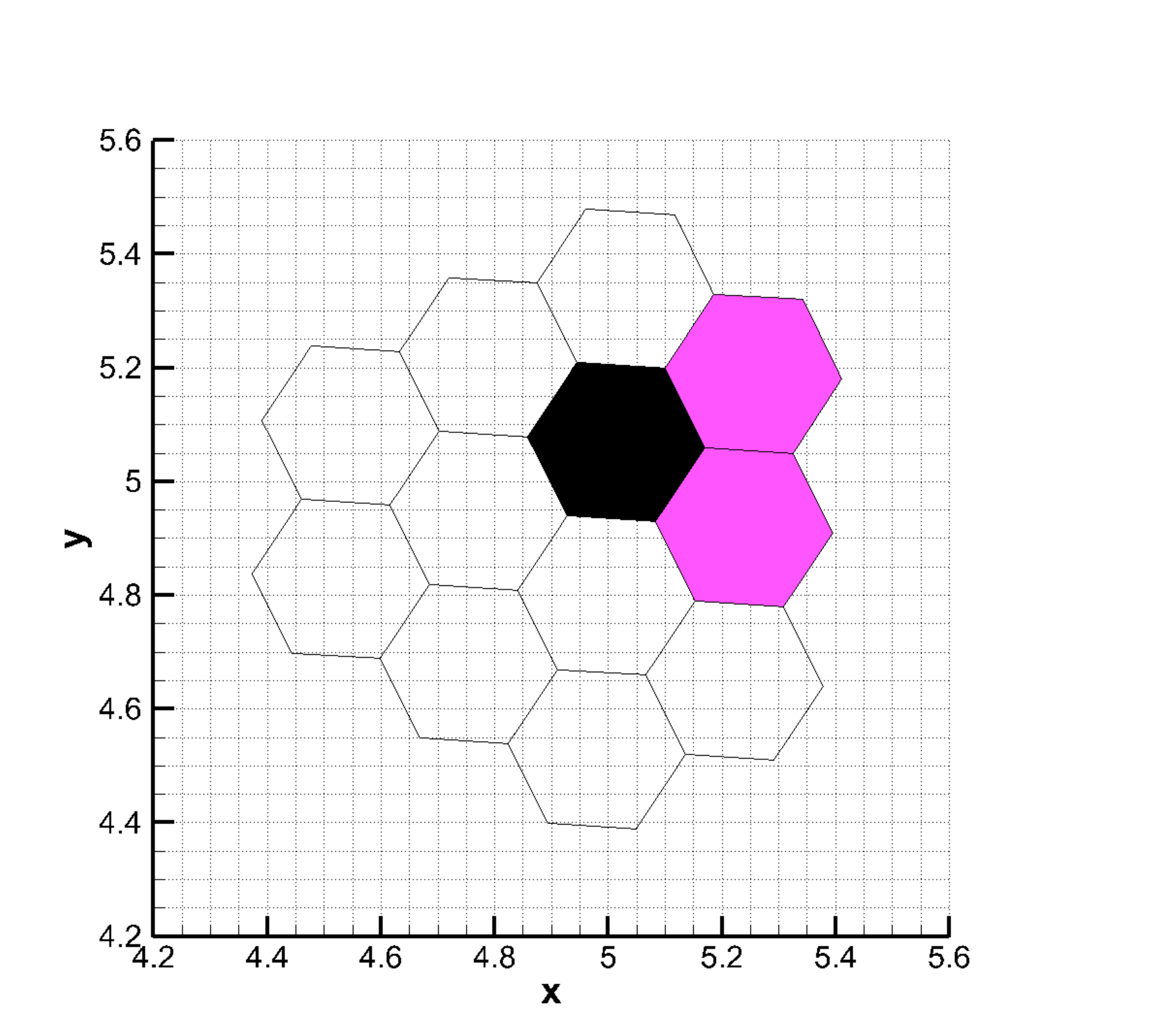} &
			\includegraphics[width=0.3\textwidth]{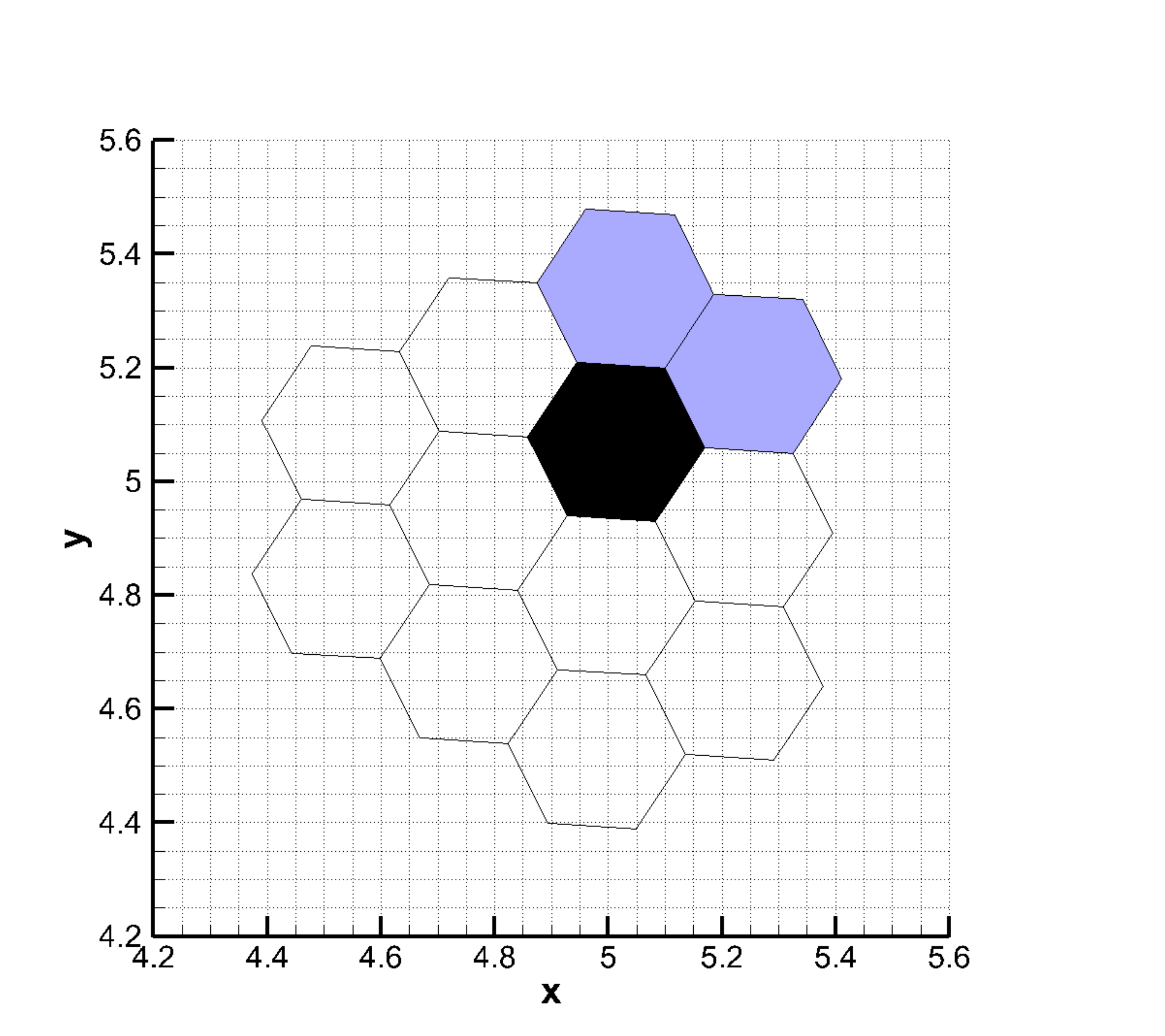}  \\
		\end{tabular}
		\caption{Reconstruction stencils for second order degree optimal polynomials. The total number of neighbors is $n_e=12$. The central $P_2$ stencil is highlighted in gray, while several one-sided $P_1$ reconstruction stencils are shown with different colours.}
		\label{fig.RecStencil}
	\end{center}
\end{figure}

Let us describe the data interpolation on this mesh. The cell averages of each component of the vector representing the distribution function at time $t^n$ is obtained by
\begin{equation}
\bar f_{k,i}^n = \frac{1}{|P_i|} \int_{P_i} f_k(x,t^n) dx, \ \forall k\in \mathcal{K}, \ x\in P_i, \ i\in[1,N_p],    
\label{eqn.cellaverage}
\end{equation}  
with $|P_i|$ the volume of cell $P_i$. We want now to construct a high order non-oscillatory polynomial representation of the distribution function $f(x,v_k,t)$. This will be called $f^w_{k,i}(x,t^n)$. We first introduce the (eventually) oscillating optimal polynomial function of arbitrary degree $M$: $p_{opt,k,i}(x)$. This has $\mathcal{M}(M)=\frac{1}{2} (M+1)(M+2)$ degree of freedom.
We consider then a central stencil $\mathcal{\tilde S}_i$ composed of $n_e=2 \cdot \mathcal{M}(M)$ cells containing the cell $P_i$
where the solution needs to be computed and its closest $(n_e-1)$ neighbors (see Figure \ref{fig.RecStencil} for a visual explanation). As suggested in \cite{barthlsq}, on general unstructured grids it is better to consider a total number of stencil elements $n_e$ larger than the exact number of degrees of freedom $\mathcal{M}(M)$ for avoiding ill posed interpolation problems which may arise due to the unstructured nature of the mesh. The polynomial $p_{opt,k,i}(x)$ relative to the stencil $\mathcal{\tilde S}_i$, is then defined by solving \begin{equation}
	\label{CWENO:Popt}
	p_{opt,k,i}(x) = \underset{{p\in\mathcal{P}_i}}{\argmin}  
	\sum_{P_j \in \mathcal{\tilde S}_i} 
	\left( \bar f^n_{k,j}-\frac{1}{|P_{j}|} \int_{P_{j}} p(x) dx \right)^2,
\end{equation}
where $\mathcal{P}_i $ is the set of all polynomials $\mathbb{P}_M$ of degree at most $M$, satisfying
\begin{equation}\label{CWENO:Popt1}
	\mathcal{P}_i = \left\{p \in \mathbb{P}_M: \bar f^n_{k,i}=\frac{1}{|P_{i}|} \int_{P_{i}} p(x) dx\right\}
	\subset \mathbb{P}_M.
\end{equation}
This means that the function $p_{opt,k,i}(x)$ is, among all the possible polynomials of degree $M$, the only one that shares the same cell average of the distribution function in the cell $P_i$, i.e. $ \bar f^n_{k,i}$, while being close in the least-square sense to the other cell averages in the stencil $\mathcal{\tilde S}_i$. In order to complete the reconstruction, one needs now a set of polynomials of degree one which permit to eliminate spurious oscillations when present. These are obtained defining a new set of stencils, all including the central cell $P_i$, and other two cells. An example of such sets is reported in Figure \ref{fig.RecStencil}.
For each stencil $S^L_i$, we indicate by $p^L_{k,i}(x)$ this linear polynomial. It remains to define one last central polynomial $p^0_{k,i}(x)$. This is computed by difference between the polynomial $p_{opt,k,i}(x)$ and the linear combination of the polynomial functions $p^L_{k,i}(x)$ of lower degree, that is
\begin{equation}
	\label{CWENO:P0}
	p^0_{k,i}(x)= \frac{1}{\lambda_0}\left(p_{opt,k,i}(x) - \sum_{L=1}^{N_{V_i}} \lambda_{L,k,i} p^L_{k,i}(x) \right) \in\mathbb{P}_M,
\end{equation}
where $\lambda_{0,k,i},\ldots,\lambda_{N_{V_i},k,i}$ are positive coefficients such that $\sum_{s=0}^{N_{V_i},k,i}\lambda_{s,k,i}=1.$
Thus, the final polynomial $f^w_{k}(x,t^n)$ is obtained by a hybridization among all stencils $s \in [0,N_{V_i}]$ as
\begin{equation}
	\label{CWENO:Prec}
	f^w_{k}(x,t^n) = \sum_{s=0}^{N_{V_i},k,i} \omega_{s,k,i} p_{k,i}^s(x),
\end{equation}
where the nonlinear weights $\omega_{s,k,i}$ are given by the standard expressions
\begin{equation}
	\label{eqn.weights}
	\omega_{s,k,i} = \frac{\tilde{\omega}_{s,k,i}}{\sum \limits_{s=0}^{N_{V_i}} \tilde{\omega}_{sk,k,i}}, 
	\qquad \textnormal{ with } \qquad 
	\tilde{\omega}_{s,k,i} = \frac{\lambda_{s,k,i}}{\left(\sigma_{s,k,i} + \epsilon \right)^r}. 
\end{equation} 
where $\epsilon=10^{-14}$ and $r=4$ are chosen according to \cite{Dumbser2007693} and $\sigma_{s,k,i}$ are standard oscillation indicators.

The polynomials $p(x)$ are expressed through the following conservative expansion
\begin{equation}
\label{eqn.recpolydef} 
p(x) = \sum \limits_{\ell=1}^\mathcal{M} \varphi_\ell(x) {\hat p}_{l,i} ,   
\end{equation}
with $\hat{p}_{l,i}$ representing the unknown expansion coefficients. The basis functions $\varphi_\ell$ are defined using a Taylor series in the space variables $x$ of components $(x_1,x_2)$ of degree $M$ and 1 for the central and the lateral stencils, respectively. This expansion is directly defined on the physical element $P_i$, expanded about its center of mass $x_{b_i}$ and normalized by the characteristic length $h_i$ of the element (the radius of the circumcircle of $P_i$)
\be
\label{eq.Dubiner_phi_spatial}
\varphi_\ell(x_1,x_2)|_{P_i} = \frac{(x_1 - x_{b_i,1})^p}{h_i^p} \, \frac{(x_2 - x_{b_i,2})^q}{h_i^q} - \int_{P_i}\frac{(x_1 - x_{b_i,1})^p}{h_i^p} \, \frac{(x_2 - x_{b_i,2})^q}{h_i^q} \, dx, 
\ee
\be
\ell = 0, \dots, \mathcal{M}, \quad (p,q) = 0,\dots, M, \ p+q\leq M, \nonumber
\ee 
In the above construction and in the reconstruction equations \eqref{CWENO:Popt} and \eqref{CWENO:Popt1}, the integrals are computed in each Voronoi polygon $P_{i}$ by summing the contribution of each sub-triangle $T\in \mathcal{T}(P_{i})$ with Gauss quadrature rule of suitable accuracy (\cite{stroud}).

Once the reconstruction procedure has been carried out, a finite volume method is used to solve the governing equations. Equation \eqref{eq:DM1_gen} is integrated on each control volume $P_i$ obtaining
\begin{equation}
	\partial_t\int_{P_i} 
	f_{k} \, dx  + \int_{P_i} v_{k} \cdot\nabla_{x}f_{k} \, dx = \int_{P_i} Q_k(f) \, dx, \ k=1,..,N, \ i=1,..,N_p,
	\label{finite_vol1}
\end{equation}
where $Q(f_k)$ represents or the BGK operator \eqref{eq:BGK} or the Boltzmann one \eqref{eq:Boltzmann}.
Using the divergence theorem, we get
\begin{equation}
	\partial_t\int_{P_i} 
	f_{k} \, dx = 
	- \sum_{j=1}^{N_{V_i}} \int_{\partial P_{i_j}} 
	L(f_k) \cdot \ {n_{ij}} \, dS +
	\int_{P_i} Q_k(f) \, dx, \quad k=1,..,N, \ i=1,..,N_p,
	\label{finite_vol3}
\end{equation}
with $n_i$ representing the unit outward normal to the element $P_i$, $\partial P_{i_j}$ denoting the face shared between element $P_i$ and its neighbor $P_j$ and $L(f_k)$ being the flux function $(v_kf_k)$. Finally, by introducing a first order in time explicit Euler discretization and using the finite volume interpretation, we get
	\begin{equation}
		\bar f_{k,i}^{n+1} = 
		\bar f_{k,i}^{n}
		- \frac{\Delta t}{|P_i|}\sum_{j=1}^{N_{V_i}} \int_{\partial P_{i_j}} 
		L(\bar f^n_{k,i}) \cdot \ {n_{ij}} \, dS +\frac{\Delta t}{|P_i|}
		\int_{P_i} Q_k(\bar f^n_{i}) \, dx, \qquad k=1,..,N, \quad i=1,..,N_p,
		\label{finite_vol5}
\end{equation}
where the boundary fluxes are evaluated relying on the Rusanov numerical flux function \cite{Rusanov:1961a}. Now, in order to obtain a high order in space scheme, we employ the CWENO reconstruction presented in the previous paragraph. This implies that $\bar f^n_{k,i}$ is replaced by $f^{n,w}_{k,i}$ for the computation of the flux functions and that $Q_k(\bar f^n_{k})$ is replaced by $Q_k^w( \bar f^n_{k})$ for the computation of the high order space approximation of the collision term. This last reconstruction needs the knowledge of the cell averages $Q_k(\bar f^n_{k})$. These will be obtained from the approximation of the collision integral through the spectral methods described in the next part. The high order in space finite volume scheme reads then
\begin{equation}\label{ho}
\bar f_{k,i}^{n+1} = 
\bar f_{k,i}^{n}- \frac{\Delta t}{|P_i|}\sum_{j=1}^{N_{V_i}} \int_{\partial P_{i_j}} 
L(f^{w,n}_{k,i}(x)) \cdot \ {n_{ij}} \, dS +\frac{\Delta t}{|P_i|}
\int_{P_i} Q^w_k(f^n_{i})(x) \, dx, \qquad k=1,..,N, \quad i=1,..,N_p.
\end{equation}

The above scheme, apart from being only first order in time, suffers from the stiffness of the equation when the Knudsen number is small. Next, we discuss high order time discretizations and a remedy to the stiffness of the kinetic equations.


\subsection{Spectral discretization of the Boltzmann operator}
\label{sec:spectral}
When the Boltzmann model \eqref{eq:Boltzmann} is employed, a fast spectral method is used to evaluate the collision term $Q(\bar f_i^n)$, and subsequently this result is projected over the discrete ordinate space giving $Q_k(\bar f_{i}^n)$ employed in \eqref{finite_vol5} and then in the high order space discretization \eqref{ho}. We proceed as follows. 
\begin{itemize}
	\item We compute an approximation of the Boltzmann operator in the center of the control volume for each cell $P_i$.
	\item Then, as for the distribution $f_{k,i}^n$, for each $k=1,..,N$, we reconstruct its value in space through the procedure of Section \ref{sec_num_approx}. This gives $Q^w_k( f_{k}^n)(x)$.
\end{itemize}
Since collisions act only at a local level, we restrict ourselves in this section on a given cell $P_i$ and we omit the dependence of the distribution function $f$ over space and time. The same computation is repeated for all volumes $P_i, \ i=1,..,N$ and times $t^n, \ n=0,..,N_t$. 
Only the dependency on the velocity variable $v$ is considered for the distribution function $f$, i.e. $f=f(v)$. In addition, since we deal with two dimensions in velocity space we also restrict to this situation in the description of the method.

We suppose, as done for the discrete ordinate method, the distribution function $f$ to have compact support on the ball $\Ball_0(R)$ of
radius $R$ centered in the origin. Then, to avoid the aliasing phenomena, at least in a single time iteration, we restrict $f(v)$ on $[-T,T]^{2}$ with $T=(3+{\sqrt 2})R/2$ and we assume $f(v)=0$ outside the ball $\Ball_0(R)$. We also extend $f(v)$ to a periodic function on the set $[-T,T]^{2}$. Let observe that $\supp (Q(f)(v)) \subset \Ball_0({\sqrt 2}R)$ and thus this choice is enough to avoid aliasing in one time step. However, the support of $f$ increases with time and thus aliasing can only be minimized by this approach. To simplify notation, let us now take $T=\pi$. This gives $R=\lambda\pi$ with $\lambda = 2/(3+\sqrt{2})$. Using one index to denote the $2$-dimensional sums, we have that the approximate function $f_N$ can be represented as the truncated Fourier
series by
\begin{equation}
f_N(v) = \sum_{k=-N/2}^{N/2} \f_k e^{i k \cdot v},
\label{eq:FU}
\end{equation}
\begin{equation}
\f_k = \frac{1}{(2\pi)^{2}}\int_{[-\pi,\pi]^{d_v}} f(v)
e^{-i k \cdot v }\,dv.
\label{eq:FC}
\end{equation}
We then obtain a spectral quadrature of our collision operator by projecting the truncated Boltzmann operator onto the space of trigonometric polynomials of degree less or equal to $N$, i.e.
\begin{equation}
{\hat Q}_k=\int_{[-\pi,\pi]^{2}}
\Q_b(f_N)
e^{-i k \cdot v}\,dv, \quad k=-N/2,\ldots,N/2. 
\label{eq:VAR}
\end{equation}
By substituting expression (\ref{eq:FU}) in (\ref{eq:VAR}) one gets
\begin{equation}
{\hat Q}_k = \sum_{\substack{l,m=-N/2\\l+m=k}}^{N/2} \f_l\,\f_m
\bb(l,m),\quad k=-N,\ldots,N,
\label{eq:CF1}
\end{equation}
where $\bb(l,m)=\B(l,m)-\B(m,m)$ are given by
\begin{equation}
\B(l,m) = \int_{\Ball_0(2\lambda\pi)}\int_{\ens{S}^{1}} 
|q| \sigma(|q|, \cos\theta) e^{-i(l\cdot q^++m\cdot q^-)}\,d\omega\,dq. \label{eq:KM}
\end{equation}
with \begin{equation}
q^{+} = \frac12(q+\vert q\vert \omega), \quad
q^{-} = \frac12(q-\vert q\vert \omega).
\label{eq:VV2}
\end{equation}
The usage of (\ref{eq:CF1}) as approximation of the collision operator requires $O(N^4)$ operations. In order to reduce the number of operations, one can express the Boltzmann operator adopting another representation \cite{Carl:EB:32}. Omitting the details, this second representation reads 
\begin{equation}
\label{defQBCarleman}
Q_b (f)= \int_{\R^{2}} \int_{\R^{2}} {\tilde B}(x,y) 
\delta(x \cdot y) 
\left[ f(v + y) \, f(v+ x) - f(v+x+y) \, f(v) \right] \, dx \,
dy,
\end{equation} 
with 
\begin{equation}\label{eq:Btilde}
\tilde{B}(|x|,|y|) =
2 \, \sigma\left(\sqrt{|x|^2+|y|^2}, \frac{|x|}{\sqrt{|x|^2+|y|^2}} \right). 
\end{equation}
This transformation yields the following new spectral quadrature formula 
\begin{equation}\label{eq:ode}
\hat{Q}_k  =
\sum_{\underset{l+m=k}{l,m=-N/2}}^{N/2} {\hat{\beta}}_F(l,m) \, \hat{f}_l \, \hat{f}_m, \ \ \
k=-N,...,N
\end{equation}
where ${\hat{\beta}}_F(l,m)=\B_F(l,m)-\B_F(m,m)$ are now 
given by
\begin{equation}
\B_F(l,m) = \int_{\Ball_0(R)} \int_{\Ball_0(R)}
\tilde{B}(x,y) \, \delta(x \cdot y) \, 
e^{i (l \cdot x+ m \cdot y)} \, dx \, dy.
\label{eq:FKM}
\end{equation}
Now, \eqref{eq:ode} can be evaluated with a convolution structure by approximating each
${\hat{\beta}}_F(l,m)$ by a sum
\[ {\hat{\beta}}_F(l,m) \simeq \sum_{p=1} ^{A} \alpha_p (l) \alpha' _p (m), \]
where $A$ represents the number of finite possible directions of collisions. This gives a sum of $A$ discrete convolutions and, consequently, the algorithm can be computed in $O(A \, N^2 \log_2 N^2)$ operations by means of
standard FFT technique (see \cite{FiMoPa:2006} for details). 

\subsection{Time Discretization}
\label{time_discr}
We discuss how to improve the time discretization and to get a fully high order space and time discretization. This is done through the use of a particular class of Implicit-Explicit (IMEX) methods \cite{Ascher,Dimarco_stiff2}. These methods are needed to handle the different scales induced by the collision and the transport operators. 
In particular, in kinetic theory one is interested in having a method which is able to handle with the same time step $\Delta t$ and the same accuracy different regimes (dense or rarefied) identified by different values of the Knudsen number $\varepsilon$. 
Our aim is precisely to be capable of capturing the fluid-limit \eqref{eq:sys1} without time step limitations due to the fast scale dynamics induced by collisions. This request is equivalent to the notion of \textit{asymptotic-preserving (AP)} schemes \cite{ACTA,degondrev,DegondAP,Jin_review} and the schemes here introduced belong to this special class. In particular, we consider the extension to the finite volume framework of the schemes discussed from the theoretical point of view in \cite{Dimarco_stiff1} and \cite{Dimarco_stiff2}.
The general formulation of such a method is
\begin{eqnarray}
	\bar F_{k,i}^{(l)} &=& \displaystyle \bar f_{k,i}^{n}-\Delta t \sum_{m=1}^{l-1} \ta_{lm} \langle L( F_{k,i}^{w,(m)})\rangle_x+\Delta t \sum_{m=1}^{\iota} a_{lm}\langle Q^w_k(F_{i}^{(m)})\rangle_x\label{eq:GIMEX} \\
	\bar f_{k,i}^{n+1} &=& \displaystyle \bar f_{k,i}^{n}-\Delta t\sum_{m=1}^{\nu}\tilde\omega_{m} \langle L( F_{k,i}^{w,(m)})\rangle_x+\Delta t\sum_{m=1}^{\iota}\omega_{m}\langle Q_k^w(F_{i}^{(m)})\rangle_x.
	\label{eq:GIMEX1}
\end{eqnarray}
In the above formula, the coefficients $\ta_{lm}$, $a_{lm}$ characterize the explicit respectively the implicit Runge-Kutta method together with the vectors $\tw =( \widetilde{\omega}_{1},..,\tw_{\iota})^{T}$, $\omega =(w_{1},..,w_{\iota})^{T}$. We refer to~\cite{Ascher,PRimex} for more details on such methods in a general framework (order conditions and properties). The functions $\bar F_{k,i}^{(l)}$ given by \eqref{eq:GIMEX} are the so-called stage values of the Runge-Kutta method which identify the cell averages of the solution at different time levels between $[t^n,t^{n+1}]$. Using these values it is possible to determine the quantities
\be \langle L( F_{k,i}^{w,(l)})\rangle_x=\frac{1}{|P_i|}\sum_{j=1}^{N_{V_i}} \int_{\partial P_{i_j}} 
L(F^{w,(l)}_{k,i}(x)) \cdot \ {n_{ij}} dS, \ee
and
\be\label{eq:ints}
\langle Q^w_k(F^{(l)}_{i})\rangle_x=\frac{1}{|P_i|}\int_{P_i} Q^w_k(F^{(l)}_{i})(x) \, dx,
\ee
where the integrals are obtained by a suitable Gauss quadrature formula. From \eqref{eq:ints}, it is clear that the direct application of such method is very difficult in practice. This, in fact, would need the inversion of a nonlinear system involving the operator $\langle Q^w_k(F^{(l)}_{i})\rangle_x$ at each stage of the Runge-Kutta time stepping. In order to find remedy to this situation, we now distinguish the case of the BGK and the case of the Boltzmann operator. 

For the case of the BGK operator, we observe that this high order integration can be conveniently written as
\begin{equation}\label{average}
	\langle Q^w_k(F^{(l)}_{i})\rangle_x=\frac{1}{|P_i|}\int_{P_i} Q^w_k( F^{(l)}_{i})(x) \, dx=\langle Q_k( \bar F^{(l)}_{i})\rangle_x+\mathcal{O}(\Delta x^2),
\end{equation}
where $\Delta x$ is the typical mesh size and $\bar F^{(l)}_{i}$ are the cell average quantities at the stage $l$. This means that the high order integral of the collision operator can be expressed by the integral of the cell centered operator plus a second order space error.
Thanks to this observation, then the time integrator can be recast as 
\begin{eqnarray}
	\bar F_{k,i}^{(l)} &=& \displaystyle \bar f_{k,i}^{n}-\Delta t \sum_{m=1}^{l-1} \ta_{lm}\left( \langle L( F_{k,i}^{w,(m)})\rangle_x+\Delta Q_k(F^{(m)}_{i})\right)+\Delta t\sum_{m=1}^{\iota} a_{lm}\langle Q_k(\bar F_{k}^{(m)})\rangle_x\label{eq:GIMEX_n} \\
	\bar f_{k,i}^{n+1} &=& \displaystyle \bar f_{k,i}^{n}-\Delta t\sum_{m=1}^{\nu}\tilde\omega_{m}\left \langle L( F_{k,i}^{w,(m)})\rangle_x+\Delta Q_k(F^{(m)}_{i})\right)+\Delta t\sum_{m=1}^{\iota}\omega_{m}\langle Q_k(\bar F_{i}^{(m)})\rangle_x.
	\label{eq:GIMEX_n2}
\end{eqnarray}
where the quantities $\Delta Q_k(F^{(m)}_{i})$ are given by
\begin{equation}
	\Delta Q_k(F^{(m)}_{i})=\langle Q^w_k(F^{(m)}_{i})\rangle_x-\frac{1}{|P_i|}\int_{P_i} Q_k(\bar F^{(m)}_{i}) \, dx, =\langle Q^w_k(F^{(m)}_{i})\rangle_x-\langle Q_k(\bar F_{i}^{(m)})\rangle_x,
\end{equation}
i.e. they represent the difference between the high order evaluation of the source term and the respective cell average. This permits a direct evaluation of the implicit terms without resorting to the approximate solution of nonlinear systems. To better understand this, let us remark that the stage evaluation of the cell average of the distribution function \eqref{eq:GIMEX_n} can be rewritten as
\begin{equation}
	\bar F_{k,i}^{(l)} = \displaystyle \bar f_{k,i}^{n}-\Delta t \sum_{m=1}^{l-1} \ta_{lm}\left( \langle L(\bar F_{k,i}^{w,(m)})\rangle_x+\Delta Q_k(F^{(m)}_{i})\right)+\Delta t\sum_{m=1}^{l-1} a_{lm}\langle Q_k( \bar F_{i}^{(m)})\rangle_x+\frac{\Delta t}{|P_i|}\int_{P_i} \frac{\nu}{\varepsilon} \left(\E_{k}[\bar U_i^{(l)}]-\bar F^{(l)}_{k,i}\right) \, dx,
	\label{eq:GIMEX3}
\end{equation}
where the only implicit term is now diagonal factor $\frac{\nu}{\varepsilon} \left(\E_{k}[\bar U_i^{(l)}]-\bar F^{(l)}_{k,i}\right)$ since the implicit methods which are considered in this work are all diagonally implicit. Let us notice now that $\bar F^{(l)}_{k,i}$ can easily be moved to the left hand side of equation \eqref{eq:GIMEX3}, while $\E_{k}[\bar U_i^{(l)}]$ depends only on the first three moments of the distribution function $\bar U_i^{(l)}$ and these are easily obtained by integrating in velocity space the cell centered values of the distribution function again from \eqref{eq:GIMEX3}.
Thus, to have a completely determined method one needs only to compute 
the moments $\bar U_i^{(l)}$. This is done by
integrating equation (\ref{eq:GIMEX3}) in velocity space which is explicitly obtained, since collision terms disappear as a result of the integration, from formula
\begin{eqnarray}
	\bar U_i^{(l)} &=& \bar U_i^{n}-\Delta t \sum_{m=1}^{l-1} \ta_{lm}  L_U^w(\bar U_i^{n}),
\end{eqnarray}
where $ L_U^w(\bar U_i^{n})$ is a high order evaluation of the macroscopic fluxes for respectively density, momentum and energy.
As a consequence, $\bar U_i^{(l)}$ and thus $\E_{k}[\bar U_i^{(l)}]$ can be explicitly evaluated and then the scheme (\ref{eq:GIMEX_n})-(\ref{eq:GIMEX_n2}) is \textit{explicitly} solvable. This concludes the presentation of the scheme for the BGK case.

In the case of the Boltzmann operator \eqref{eq:Boltzmann}, one cannot employ directly an equivalent strategy to the one above described for the BGK case. In fact, even if equation \eqref{average} still holds true, it is not possible to invert the cell centred integral of $Q(f)$, $\langle Q_k( \bar F^{(l)}_{i})\rangle_x$, in this case. Consequently, the scheme (\ref{eq:GIMEX_n})-(\ref{eq:GIMEX_n2}) is not explicitly solvable and another solution has to be found. This problem can be circumvented relying on a penalization strategy where the penalization term is the cell centred BGK operator.
In fact, in stiff regimes the solution $f$ is close to the Maxwellian equilibrium $M_f$ and then the Boltzmann and the BGK models furnish closer and closer solutions as the fluid regime is approached. Thus, one can rewrite equation \eqref{eq:B} as
\be
\partial_t
f+v\cdot\nabla_x
f=\frac{1}{\varepsilon}G_P(f)+\frac{1}{\varepsilon}Q_{BGK}(f),\label{eq:Bref2}
\ee
where $Q_{BGK}$ is the BGK operator and $G_p(f)=Q_b-Q_{BGK}$ with $Q_b$ the Boltzmann collision operator. Now, the idea is to use a time discretization in which $Q_{BGK}$ is implicit while $G_p(f)$ is explicit. 
The modified IMEX Runge-Kutta schemes become then
\begin{eqnarray}
\bar F_{k,i}^{(l)} &=& \displaystyle \bar f_{k,i}^{n}-\Delta t \sum_{m=1}^{l-1} \ta_{lm}\left(\frac{1}{\varepsilon} G_{k,p}( F_i^{(m)})+ \langle L(\bar F_{k,i}^{w,(m)})\rangle_x\right)+\Delta t\sum_{m=1}^{\iota} a_{lm}\langle Q_{k,BGK}(\bar F_{i}^{(m)})\rangle_x\label{eq:GIMEXb} \\
\bar f_{k,i}^{n+1} &=& \displaystyle \bar f_{k,i}^{n}-\Delta t
\sum_{m=1}^{\nu}\tilde \omega_{m}\left(\frac{1}{\varepsilon} G_{k,p}( F_{i}^{(m)})+ \langle L(\bar F_{k,i}^{w,(m)})\rangle_x\right)+\Delta
t\sum_{m=1}^{\iota}\omega_{m}\langle Q_{k,BGK}(\bar F_{i}^{(m)})\rangle_x,
\label{eq:GIMEX1b}
\end{eqnarray}
where 
\begin{equation}
	G_{k,P}( F^{(m)}_{i})=\langle Q_{k,b}^w( F_{i}^{(m)})\rangle_x-\langle Q_{k,BGK}(\bar F_{i}^{(m)})\rangle_x.
\end{equation}
As for the case of the BGK operator the only implicit term is $\frac{\nu}{\varepsilon} \left(\E_{k}[\bar U_i^{(l)}]-\bar F^{(l)}_{k,i}\right)$ which can be computed following the schemes \eqref{eq:GIMEX_n}-\eqref{eq:GIMEX_n2}. This concludes the presentation of the time integration methods.

Finally, thanks to the choice done about the numerical time integration, in all the numerical simulation discussed in the next section, the time step is fixed according to
\be
\Delta t =\frac{1}{2}  \left(\frac{\Delta x}{\max_\mathcal{K}(|v_k|)} \right),
\label{eq:Time2}
\ee
which means it is limited only by the hyperbolic transport term.
In practice, we consider in the next section three time integration schemes. The first one is the standard first order implicit-explicit Euler scheme for which we do not report the Butcher tableau. The second one is the second order ARS(2,2,2) \cite{Ascher} scheme 
{\small
	\[
	\begin{array}{c|ccc}
		0 & 0 & 0 & 0   \\
		\gamma   & \gamma & 0 & 0 \\
		1   & \delta & 1-\delta & 0\\
		\hline
		& \delta &  1-\delta & 0
	\end{array}\qquad
	\begin{array}{c|ccc}
		0 & 0 & 0 & 0 \\
		\gamma   & 0 & \gamma & 0 \\
		1   & 0 & 1-\gamma & \gamma\\
		\hline
		& 0 & 1-\gamma & \gamma
	\end{array}
	\]
}
with $\gamma=1-1/\sqrt{2}$ and $\delta=1-1/(2\gamma)$, while the third one is the third order BPR(3,4,3) \cite{BPR}
{\small
	\[
	\begin{array}{c|cccccc}
		0 & 0 & 0 & 0 & 0 & 0  \\
		1   & 1 & 0 & 0 & 0 & 0 \\
		2/3   & 4/9 & 2/9 & 0 & 0 & 0\\
		1   & 1/4 & 0 & 3/4 & 0 & 0\\
		1   & 1/4 & 0 & 3/4 & 0 & 0\\
		\hline
		& 1/4 & 0 & 3/4 & 0 & 0\\
	\end{array}\qquad
	\begin{array}{c|cccccc}
		0 & 0 & 0 & 0 & 0 & 0   \\
		1   & 1/2 & 1/2 & 0 & 0& 0 \\
		2/3   & 5/18 & -1/9 & 1/2 & 0& 0 \\
		1   & 1/2 & 0 & 0 & 1/2& 0 \\
		1   & 1/4 & 0 & 3/4 & -1/2& 1/2 \\
		\hline
		& 1/4 & 0 & 3/4 & -1/2& 1/2 \\
	\end{array}
	\]
}
The stability properties of these schemes, the capability of describing different collisional regimes and the consistency with the underlying compressible Euler model have been discussed in \cite{Dimarco_stiff2} in the case of finite difference method for the Boltzmann and the BGK equations while in the case of finite volume methods for the sole BGK model in \cite{BosDim}.

\section{Numerical results}
\label{sec.validation}
In this section, we measure the capability of the presented method to describe fluid flows under different regimes. Our analysis starts by a set of numerical convergence tests in velocity, in space, in time and in space-time. These convergence studies are carried out in different regimes identified by the value of the Knudsen number. In addition, we also report the computational efforts involved in the resolution of the Boltzmann and the BGK models. The second part is devoted to some standard benchmark studies in rarefied gas dynamics. The last part is dedicated to applications: we study the flow around a NACA 0012 airfoil. This study is possible only thanks to the arbitrarily shaped finite volume approach used, which permits to describe the generic geometry considered. For all computations the gas polytropic index is set to $\gamma=2$.
%
\subsection{Numerical convergence studies}
\begin{figure}[!htbp]
	\begin{center}
		\includegraphics[width=0.5\textwidth]{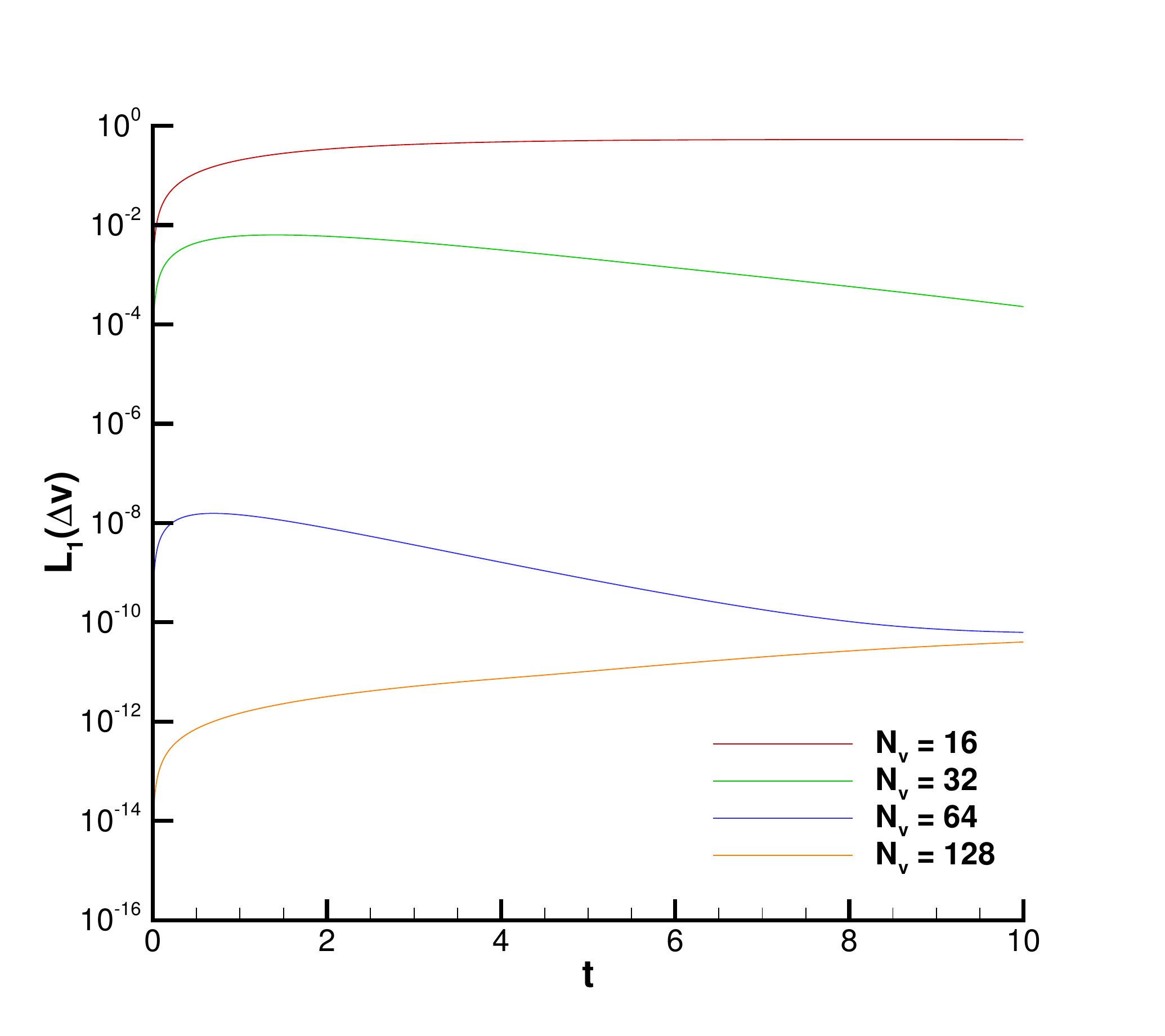}
		\caption{Time evolution of the $L_1$ error in velocity space computed according to \eqref{error} up to $t=10$.}
		\label{fig.vel-error}
	\end{center}
\end{figure}

This part is dedicated to the convergence tests. The first test validates the correct implementation of the spectral method for the discretization of the Boltzmann operator. Here we resolve the Boltzmann homogeneous equation in the 2D velocity space. The following analytical solution \cite{Bobylev1975} 
\begin{equation}
	f(v,t) = \frac{\exp{(-v^2/2S)}}{2\pi S^2} \Bigg[2S - 1 + \frac{1-S}{2S}v^2\Bigg]	\qquad	\text{with}	\qquad	S=1-\frac{\exp{(-t/8)}}{2},	
\end{equation}
is approached by the spectral scheme. We adopt a domain $\mathcal{V}=[-12;12]\times [-12;12]$, for which we consider different values velocity points $N_v$ corresponding also to the number of modes of the Fourier expansion. The different configurations for the velocity grid are reported in Table \ref{tab.conv_spectral}. 
The errors in the $L_1$ and $L_2$ norms are obtained with the following formula
\begin{equation}
	L_p (\Delta v, t) = \Bigg(\frac{\sum_{k=1}^{N_v} |f_{\Delta v}(v_k,t) - f_{ref}(v_k,t)|^p}{\sum_{k=1}^{N_v} |f_{ref}(v_k,t)|^p}\Bigg)^{1/p},
	\label{error}
\end{equation}
with $\Delta v$ representing the characteristic mesh size in the velocity space. The quantities $f_{\Delta v}$ and $f_{ref}$ represent the numerical and reference solutions of the distribution function computed on the velocity points, respectively. The $L_{\infty}$ norm is also computed, for which we simply measure the maximum distance between exact and computed solutions. Table \ref{tab.conv_spectral} shows the results of such analysis at the final time $t_{f}=0.1$. As expected, the convergence rate grows very fast and thus errors become closer to machine precision. Figure \ref{fig.vel-error} reports the $L_1$ errors over time up to $t_f=10$.

\begin{table}[!htbp]  
	\caption{Numerical convergence results for the homogeneous Boltzmann equation using the Spectral method presented in \ref{sec:spectral} on a sequence of refined Cartesian meshes of size $\Delta v$ in the velocity space. The error is measured for the distribution function in the $L_1$, $L_2$ and $L_\infty$ norms using the definition \eqref{error} at time $t_{f}=0.1$.}  
	\begin{center} 
		\begin{small}
			\renewcommand{\arraystretch}{1.0}
			\begin{tabular}{cc|cccccc}
				\hline
				$N_v$ & $\Delta v$ & $f_{L_1}$ & $\mathcal{O}(f_{L_1})$ & $f_{L_2}$ & $\mathcal{O}(f_{L_2})$ &  $f_{L_\infty}$ & $f_{L_\infty}$ \\ 
				\hline
				16  & 1.60E-00 & 2.392E-02 & -    & 1.516E-02 & -    & 1.660E-03 & -    \\ 
				32  & 7.74E-01 & 1.183E-03 & 4.1  & 5.235E-04 & 4.6  & 4.252E-05 & 5.0  \\ 
				64  & 3.81E-01 & 5.493E-09 & 17.3 & 1.702E-09 & 17.8 & 1.159E-10 & 18.1 \\ 
				128	& 1.89E-01 & 1.425E-13 & 15.1 & 8.916E-14 & 14.1 & 1.096E-14 & 13.2 \\ 
			\end{tabular}
		\end{small}
	\end{center}
	\label{tab.conv_spectral}
\end{table}

\begin{figure}[!htbp]
	\begin{center}
		\begin{tabular}{cc} 
			\includegraphics[width=0.47\textwidth]{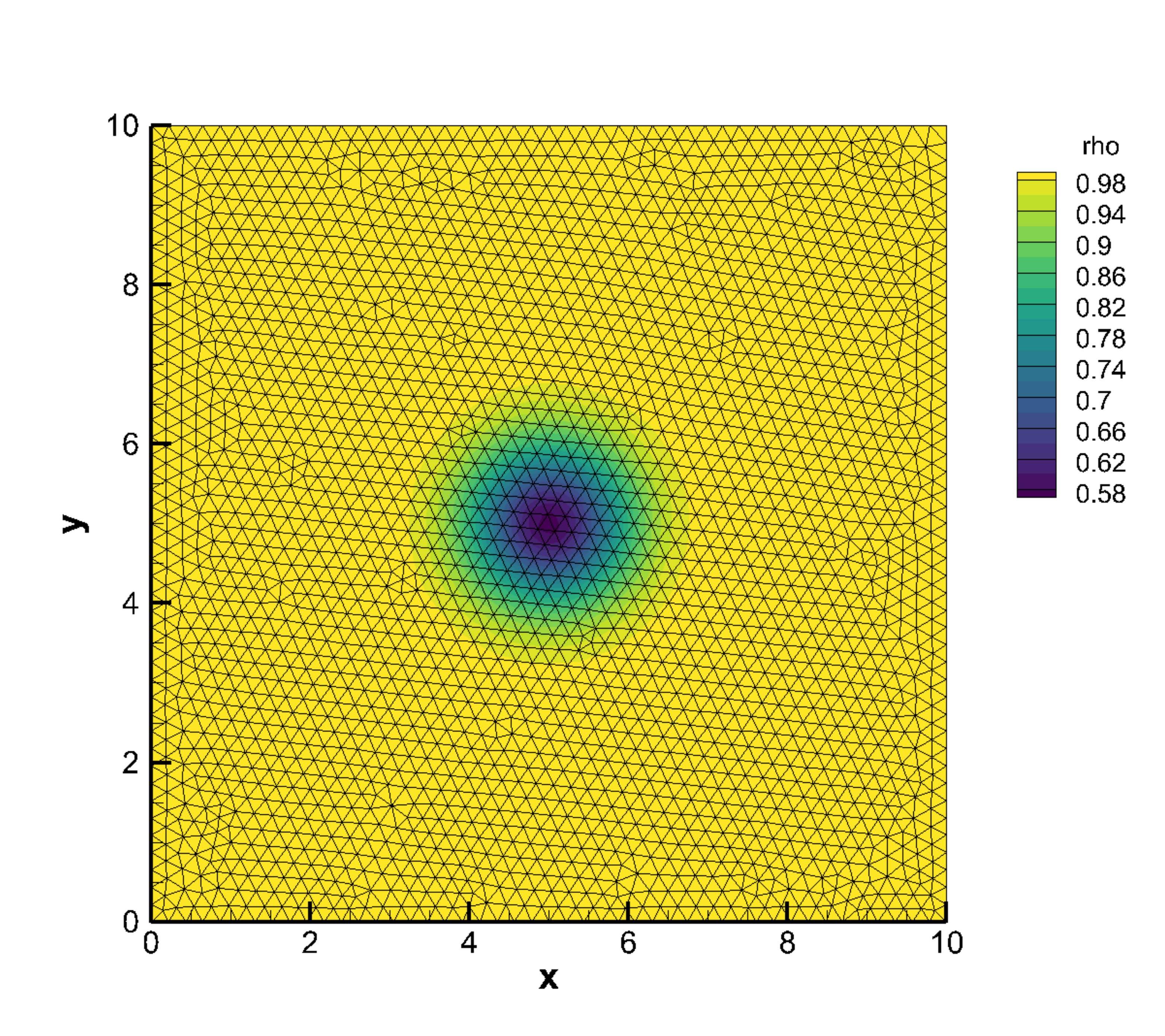}  &           
			\includegraphics[width=0.47\textwidth]{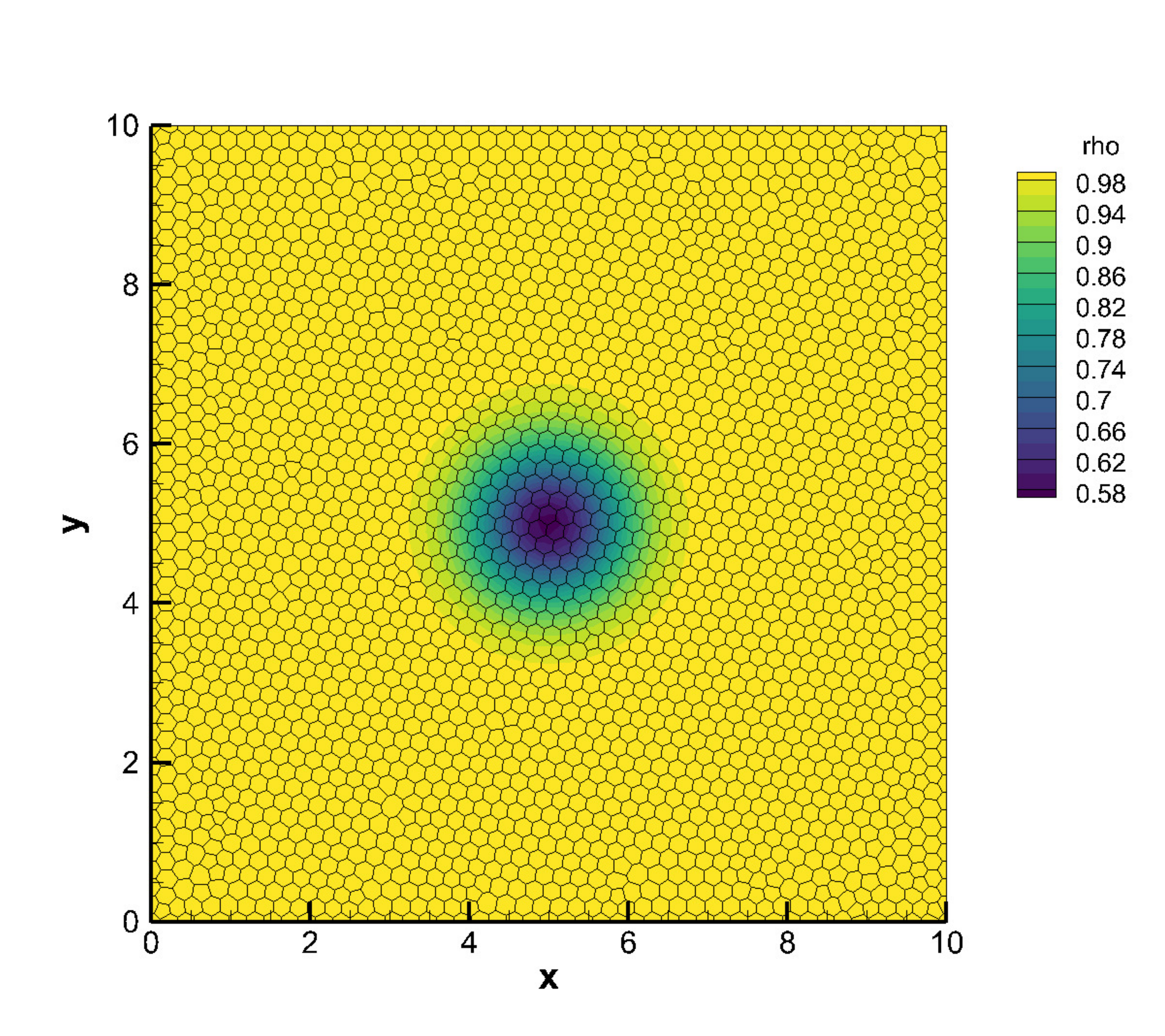} \\
		\end{tabular}
		\caption{Initial density distribution of the smooth isentropic vortex on an unstructured triangular mesh (left) and on a Voronoi tessellation (right) with characteristic mesh size $h=0.25$.}
		\label{fig.SVmesh}
	\end{center}
\end{figure}

For the second part of convergence studies, we consider as test case an isentropic vortex in physical space \cite{HuShuTri}. The velocity space is now fixed to $\mathcal{V}=[-10;10]\times [-10;10]$ and it is paved with a Cartesian grid of $900$ equal elements. Periodic boundary conditions are imposed on each side of the box. The initial condition is given by
\begin{equation}
\mathcal{U} = (\rho,u_x,u_y,T) = (1+\delta \rho, 1+\delta u_x, 1+\delta u_y, \delta T),
\label{eq.ConvEul-IC}
\end{equation}
where the perturbations for temperature $\delta T$, density $\delta \rho$ and velocity $(\delta u_x, \delta u_y)$ are
\begin{eqnarray}
\label{ShuVortDelta}
\delta T    = -\frac{(\gamma-1)\beta^2}{8\gamma\pi^2}e^{1-r^2},\quad \delta \rho = (1+\delta T)^{\frac{1}{\gamma-1}}-1,\quad 
\left(\begin{array}{c} \delta u_x \\ \delta u_y \end{array}\right) = \frac{\beta}{2\pi}e^{\frac{1-r^2}{2}} \left(\begin{array}{c} -(y-5) \\ \phantom{-}(x-5) \end{array}\right), \nonumber 
\end{eqnarray}
with $\beta=5$. 
Figure \ref{fig.SVmesh} shows the initial density distribution of this smooth isentropic vortex on a triangular and on a polygonal unstructured mesh. 
The space corresponds to the box $\Omega=[0;10] \times [0;10]$ discretized with a sequence of refined unstructured meshes of characteristic size $h(\Omega)=\big(\sum_{i=1}^{N_P} h_i\big)/N_P$ with $h_i=\sqrt{|P_i|}$.
The errors are measured in the $L_1$ and $L_2$ norms for respectively the density and the temperature as follows:
\begin{eqnarray}
L_1 &=& \int \limits_{\Omega} \left| \mathcal{U}_{ref}(x,t_f) - \mathcal{U}_h(x,t_f) \right| dx \, dy, \\
L_2 &=& \sqrt{\int \limits_{\Omega} \left( \mathcal{U}_{ref}(x,t_f) - \mathcal{U}_h(x,t_f) \right)^2 dx \, dy}.
\end{eqnarray}	
Here, $\mathcal{U}_h(x,t_f)$ represents the high order reconstructed solution for the \textit{macroscopic} quantities obtained with our scheme, while $\mathcal{U}_{ref}(x,t_f)$ is a prescribed \textit{reference solution}. We perform several convergence measures: convergence of the time discretization and convergence of the space-time discretization. More in details:
\begin{itemize}
	\item In the limit $\varepsilon \to 0$, since the analytical solution is known, we consider a \textit{space-time} error: $\mathcal{U}_{ref}(\x,t)=\mathcal{U}(x-t,0)$.
	\item If $\varepsilon \ne 0$, the reference solution is instead the one computed on the \textit{same} phase-space mesh using a very small time step $\Delta t_e\ll \Delta t$, i.e. $\mathcal{U}_{ref}(x,t_f)=\mathcal{U}_{h_t}(x,t_f)$. In this case, a time convergence is proposed.
\end{itemize}
The validation of the CWENO reconstruction procedure has been already presented up fourth order in the previous work \cite{BosDim} and we omit it here.
The results of this test are shown in Table \ref{tab.conv_0} for $\varepsilon=0$ while in Table \ref{tab.conv_1e-3} for $\varepsilon=10^{-3}$, in Table \ref{tab.conv_1e-2} for $\varepsilon=10^{-2}$ and finally in Table \ref{tab.conv_1} for $\varepsilon=1$. The numerical results confirm that the CWENO-IMEX schemes applied to the Boltzmann model achieve the theoretical order of accuracy in velocity, in space and in time.
\begin{table}[!htbp]  
	\caption{Numerical convergence results for the Boltzmann model at time $t_f=0.1$ with $\varepsilon=0$ on a sequence of refined polygonal meshes of size $h(\Omega)$. The errors are measured in $L_1$ and $L_2$ norm and refer to the variables $\rho$ (density) and $T$ (temperature).}  
	\begin{center} 
		\begin{small}
			\renewcommand{\arraystretch}{1.0}
			\begin{tabular}{c|cccc}
				\multicolumn{5}{c}{CWENO-IMEX $\mathcal{O}2$ $Kn=0$} \\
				\hline
				$h(\Omega)$ & $\rho_{L_1}$ & $\mathcal{O}(\rho_{L_1})$ & $\rho_{L_2}$ & $\mathcal{O}(\rho_{L_2})$ \\ 
				\hline
				3.84E-01 & 4.523E-02 & -   & 1.585E-02 & -   \\ 
				1.96E-01 & 8.759E-03 & 2.5 & 2.996E-03 & 2.5 \\ 
				1.32E-01 & 3.793E-03 & 2.1 & 1.218E-03 & 2.3 \\ 
				9.90E-02 & 2.078E-03 & 2.1 & 6.853E-04 & 2.0 \\ 
				\hline
				$h(\Omega)$ & $T_{L_1}$ & $\mathcal{O}(T_{L_1})$ & $T_{L_2}$ & $\mathcal{O}(T_{L_2})$ \\ 
				\hline
				3.84E-01 & 1.589E-01 & -   & 6.601E-02 & -   \\ 
				1.96E-01 & 2.770E-01 & 2.6 & 3.779E-02 & 2.5 \\ 
				1.32E-01 & 1.075E-02 & 2.4 & 1.907E-02 & 2.4 \\ 
				9.90E-02 & 5.707E-03 & 2.2 & 1.136E-02 & 2.2 \\
				\multicolumn{5}{c}{} \\
				\multicolumn{5}{c}{CWENO-IMEX $\mathcal{O}3$ $Kn=0$} \\
				\hline
				$h(\Omega)$ & $\rho_{L_1}$ & $\mathcal{O}(\rho_{L_1})$ & $\rho_{L_2}$ & $\mathcal{O}(\rho_{L_2})$ \\ 
				\hline
				3.84E-01 & 4.971E-02 & -   & 1.717E-02 & -   \\ 
				1.96E-01 & 6.951E-03 & 2.9 & 2.425E-03 & 2.9 \\ 
				1.32E-01 & 2.044E-03 & 3.1 & 6.989E-04 & 3.1 \\ 
				9.90E-02 & 9.111E-04 & 2.8 & 3.102E-04 & 2.8 \\ 
				\hline
				$h(\Omega)$ & $T_{L_1}$ & $\mathcal{O}(T_{L_1})$ & $T_{L_2}$ & $\mathcal{O}(T_{L_2})$ \\ 
				\hline
				3.84E-01 & 1.753E-01 & -   & 5.725E-02 & -   \\ 
				1.96E-01 & 2.601E-02 & 2.8 & 8.721E-03 & 2.8 \\ 
				1.32E-01 & 7.632E-03 & 3.1 & 2.564E-03 & 3.1 \\ 
				9.90E-02 & 3.315E-03 & 2.9 & 1.136E-03 & 2.8 \\
			\end{tabular}
		\end{small}
	\end{center}
	\label{tab.conv_0}
\end{table}

\begin{table}[!htbp]  
	\caption{Numerical convergence results for the Boltzmann model at time $t_f=9 \cdot 10^{-3}$ with $\varepsilon=10^ {-3}$ on a triangular mesh with characteristic mesh size $h(\Omega)=1/3$. The errors are measured in $L_1$ and $L_2$ norm and refer to the variables $\rho$ (density) and $T$ (temperature).}  
	\begin{center} 
		\begin{small}
			\renewcommand{\arraystretch}{1.0}
			\begin{tabular}{c|cccc}
				\multicolumn{5}{c}{CWENO-IMEX $\mathcal{O}2$ $\varepsilon=10^{-3}$} \\
				\hline
				$\dt$ & $\rho_{L_1}$ & $\mathcal{O}(\rho_{L_1})$ & $\rho_{L_2}$ & $\mathcal{O}(\rho_{L_2})$ \\ 
				\hline
				1.80E-03 & -         & -   & -         & -   \\ 
				9.00E-04 & 1.423E-06 & -   & 6.693E-07 & -   \\ 
				4.50E-04 & 3.794E-07 & 1.9 & 1.784E-07 & 1.9 \\ 
				2.25E-04 & 9.803E-08 & 2.0 & 4.610E-08 & 2.0 \\
				\hline
				$h(\Omega)$ & $T_{L_1}$ & $\mathcal{O}(T_{L_1})$ & $T_{L_2}$ & $\mathcal{O}(T_{L_2})$ \\ 
				\hline
				1.80E-03 & -         & -   & -         & -   \\ 
				9.00E-04 & 5.647E-06 & -   & 2.249E-06 & -   \\ 
				4.50E-04 & 1.511E-06 & 1.9 & 6.003E-06 & 1.9 \\ 
				2.25E-04 & 3.911E-07 & 1.9 & 1.552E-07 & 2.0 \\
				\multicolumn{5}{c}{} \\
				\multicolumn{5}{c}{CWENO-IMEX $\mathcal{O}3$ $\varepsilon=10^{-3}$} \\
				\hline
				$h(\Omega)$ & $\rho_{L_1}$ & $\mathcal{O}(\rho_{L_1})$ & $\rho_{L_2}$ & $\mathcal{O}(\rho_{L_2})$ \\ 
				\hline
				2.25E-04 & -         & -   & -         & -   \\ 
				1.12E-05 & 7.033E-09 & -   & 3.641E-09 & -   \\  
				5.62E-05 & 1.051E-09 & 2.7 & 5.422E-10 & 2.7 \\
				2.81E-05 & 1.566E-10 & 2.8 & 8.013E-11 & 2.8 \\ 
				\hline
				$h(\Omega)$ & $T_{L_1}$ & $\mathcal{O}(T_{L_1})$ & $T_{L_2}$ & $\mathcal{O}(T_{L_2})$ \\ 
				\hline
				2.25E-04 & -         & -   & -         & -   \\ 
				1.12E-05 & 2.631E-08 & -   & 1.182E-08 & -   \\ 
				5.62E-05 & 3.908E-09 & 2.8 & 1.734E-09 & 2.8 \\ 
				2.81E-05 & 5.782E-10 & 2.8 & 2.544E-10 & 2.8 \\ 
			\end{tabular}
		\end{small}
	\end{center}
	\label{tab.conv_1e-3}
\end{table}

\begin{table}[!htbp]  
	\caption{Numerical convergence results for the Boltzmann model at time $t_f=9 \cdot 10^{-3}$ with $\varepsilon=10^ {-2}$ on a triangular mesh with characteristic mesh size $h(\Omega)=1/3$. The errors are measured in $L_1$ and $L_2$ norm and refer to the variables $\rho$ (density) and $T$ (temperature).}  
	\begin{center} 
		\begin{small}
			\renewcommand{\arraystretch}{1.0}
			\begin{tabular}{c|cccc}
				\multicolumn{5}{c}{CWENO-IMEX $\mathcal{O}2$ $\varepsilon=10^{-2}$} \\
				\hline
				$\Delta t$ & $\rho_{L_1}$ & $\mathcal{O}(\rho_{L_1})$ & $\rho_{L_2}$ & $\mathcal{O}(\rho_{L_2})$ \\ 
				\hline
				1.80E-03 &         - & -   &        -  & -   \\ 
				9.00E-04 & 2.609E-06 & -   & 7.473E-07 & -   \\ 
				4.50E-04 & 6.438E-07 & 2.0 & 1.861E-07 & 2.0 \\ 
				2.25E-04 & 1.600E-07 & 2.0 & 4.650E-08 & 2.0 \\ 
				\hline
				$\Delta t$ & $T_{L_1}$ & $\mathcal{O}(T_{L_1})$ & $T_{L_2}$ & $\mathcal{O}(T_{L_2})$ \\ 
				\hline
				1.80E-03 &         - & -   &        -  & -   \\ 
				9.00E-04 & 1.422E-05 & -   & 4.180E-06 & -   \\ 
				4.50E-04 & 3.537E-06 & 2.0 & 1.039E-06 & 2.0 \\ 
				2.25E-04 & 8.831E-07 & 2.0 & 2.593E-07 & 2.0 \\ 
				\multicolumn{5}{c}{} \\
				\multicolumn{5}{c}{CWENO-IMEX $\mathcal{O}3$ $\varepsilon=10^{-2}$} \\
				\hline
				$\Delta t$ & $\rho_{L_1}$ & $\mathcal{O}(\rho_{L_1})$ & $\rho_{L_2}$ & $\mathcal{O}(\rho_{L_2})$ \\ 
				\hline
				1.80E-03 &         - & -   &        -  & -   \\ 
				9.00E-04 & 5.084E-08 & -   & 1.957E-08 & -   \\ 
				4.50E-04 & 6.799E-09 & 2.9 & 2.665E-09 & 2.9 \\ 
				2.25E-04 & 8.825E-10 & 2.9 & 3.491E-10 & 2.9 \\ 
				\hline
				$\Delta t$ & $T_{L_1}$ & $\mathcal{O}(T_{L_1})$ & $T_{L_2}$ & $\mathcal{O}(T_{L_2})$ \\ 
				\hline
				1.80E-03 &         - & -   &        -  & -   \\ 
				9.00E-04 & 2.994E-07 & -   & 9.614E-08 & -   \\ 
				4.50E-04 & 3.874E-08 & 3.0 & 1.256E-08 & 2.9 \\ 
				2.25E-04 & 4.930E-09 & 3.0 & 1.608E-09 & 3.0 \\
			\end{tabular}
		\end{small}
	\end{center}
	\label{tab.conv_1e-2}
\end{table}

\begin{table}[!htbp]  
	\caption{Numerical convergence results for the Boltzmann model at time $t_f=9 \cdot 10^{-3}$ with $\varepsilon=10^ {0}$ on a triangular mesh with characteristic mesh size $h(\Omega)=1/3$. The errors are measured in $L_1$ and $L_2$ norm and refer to the variables $\rho$ (density) and $T$ (temperature).}  
	\begin{center} 
		\begin{small}
			\renewcommand{\arraystretch}{1.0}
			\begin{tabular}{c|cccc}
				\multicolumn{5}{c}{CWENO-IMEX $\mathcal{O}2$ $\varepsilon=10^{0}$} \\
				\hline
				$h(\Omega)$ & $\rho_{L_1}$ & $\mathcal{O}(\rho_{L_1})$ & $\rho_{L_2}$ & $\mathcal{O}(\rho_{L_2})$ \\ 
				\hline
				1.80E-03 &         - & -   &        -  & -   \\ 
				9.00E-04 & 2.231E-06 & -   & 6.224E-07 & -   \\ 
				4.50E-04 & 5.407E-07 & 2.0 & 1.508E-07 & 2.0 \\ 
				2.25E-04 & 1.331E-07 & 2.0 & 3.711E-08 & 2.0 \\ 
				\hline
				$h(\Omega)$ & $T_{L_1}$ & $\mathcal{O}(T_{L_1})$ & $T_{L_2}$ & $\mathcal{O}(T_{L_2})$ \\ 
				\hline
				1.80E-03 &         - & -   &        -  & -   \\ 
				9.00E-04 & 1.245E-05 & -   & 3.762E-06 & -   \\ 
				4.50E-04 & 3.007E-06 & 2.0 & 9.076E-07 & 2.1 \\ 
				2.25E-04 & 7.389E-07 & 2.0 & 2.229E-07 & 2.0 \\ 
				\multicolumn{5}{c}{} \\
				\multicolumn{5}{c}{CWENO-IMEX $\mathcal{O}3$ $\varepsilon=10^{0}$} \\
				\hline
				$h(\Omega)$ & $\rho_{L_1}$ & $\mathcal{O}(\rho_{L_1})$ & $\rho_{L_2}$ & $\mathcal{O}(\rho_{L_2})$ \\ 
				\hline
				1.80E-03 &         - & -   &        -  & -   \\ 
				9.00E-04 & 2.258E-08 & -   & 7.657E-09 & -   \\ 
				4.50E-04 & 2.737E-09 & 3.0 & 9.307E-10 & 3.0 \\ 
				2.25E-04 & 3.369E-10 & 3.0 & 1.148E-10 & 3.0 \\ 
				\hline
				$h(\Omega)$ & $T_{L_1}$ & $\mathcal{O}(T_{L_1})$ & $T_{L_2}$ & $\mathcal{O}(T_{L_2})$ \\ 
				\hline
				1.80E-03 &         - & -   &        -  & -   \\ 
				9.00E-04 & 1.415E-07 & -   & 4.930E-08 & -   \\ 
				4.50E-04 & 1.708E-08 & 3.1 & 5.963E-09 & 3.0 \\ 
				2.25E-04 & 2.098E-09 & 3.0 & 7.334E-10 & 3.0 \\
			\end{tabular}
		\end{small}
	\end{center}
	\label{tab.conv_1}
\end{table}

We conclude this part by reporting in Table \ref{tab.comp_time} the computational times needed to run the simulations detailed in the rest of the section. In particular, we aim at comparing the costs involved in the computation of respectively the Boltzmann and the BGK models. The results show that, in the case of two-dimensional velocity and physical space simulations, the two models are almost equivalently expensive. Thus, the main part of the computation is devoted to the determination of the numerical fluxes in the hyperbolic part. We can conclude that in the tested cases the spectral method is indeed so fast that it does not engrave on the global simulation time.
\begin{table}[!htbp]  
	\caption{Computational time for some test cases for Boltzmann and BGK model with different Knudsen numbers. Absolute time of each simulation measured in seconds [s] and percentage of the computational time needed to perform the evaluation of the Boltzmann collision operator w.r.t. to the BGK computational time $\tau_{coll}=\left(\frac{T_{Boltz}-T_{BGK}}{T_{BGK}}\cdot 100\right)$ (Lax = Lax problem, EP2D = Explosion problem, DMR = Double Mach Reflection problem).}  
	\begin{center} 
		\begin{small}
			\renewcommand{\arraystretch}{1.0}
			\begin{tabular}{cc|cc|c}
				\hline
				Test & Knudsen & $T_{Boltz}$ [s] & $T_{BGK}$ [s] & $\tau_{coll}$ [\%] \\ 
				\hline
				Lax  & $5 \cdot 10^{-3}$ & 5.0028E+05 & 4.8807E+05 & 2.502 \\ 
				Lax  & $5 \cdot 10^{-4}$ & 5.1038E+05 & 4.9976E+05 & 2.126 \\
				Lax  & $5 \cdot 10^{-5}$ & 5.1142E+05 & 5.0108E+05 & 2.064 \\
				EP2D & $5 \cdot 10^{-3}$ & 8.4881E+06 & 8.3637E+06 & 1.488 \\ 
				EP2D & $5 \cdot 10^{-4}$ & 8.4893E+06 & 8.3828E+06 & 1.270 \\
				EP2D & $5 \cdot 10^{-5}$ & 8.4900E+06 & 8.3968E+06 & 1.110 \\
				DMR  & $5 \cdot 10^{-3}$ & 1.1573E+07 & 1.1566E+07 & 0.060 \\ 
				DMR  & $5 \cdot 10^{-4}$ & 1.1146E+07 & 1.1070E+07 & 0.690 \\
				DMR  & $5 \cdot 10^{-5}$ & 1.1134E+07 & 1.1046E+07 & 0.797 \\
			\end{tabular}
		\end{small}
	\end{center}
	\label{tab.comp_time}
\end{table}

%
\subsection{Lax problem} 
\label{sec.Lax}
\begin{figure}[!htbp]
	\begin{center}
		\begin{tabular}{ccc} 
			\includegraphics[width=0.33\textwidth]{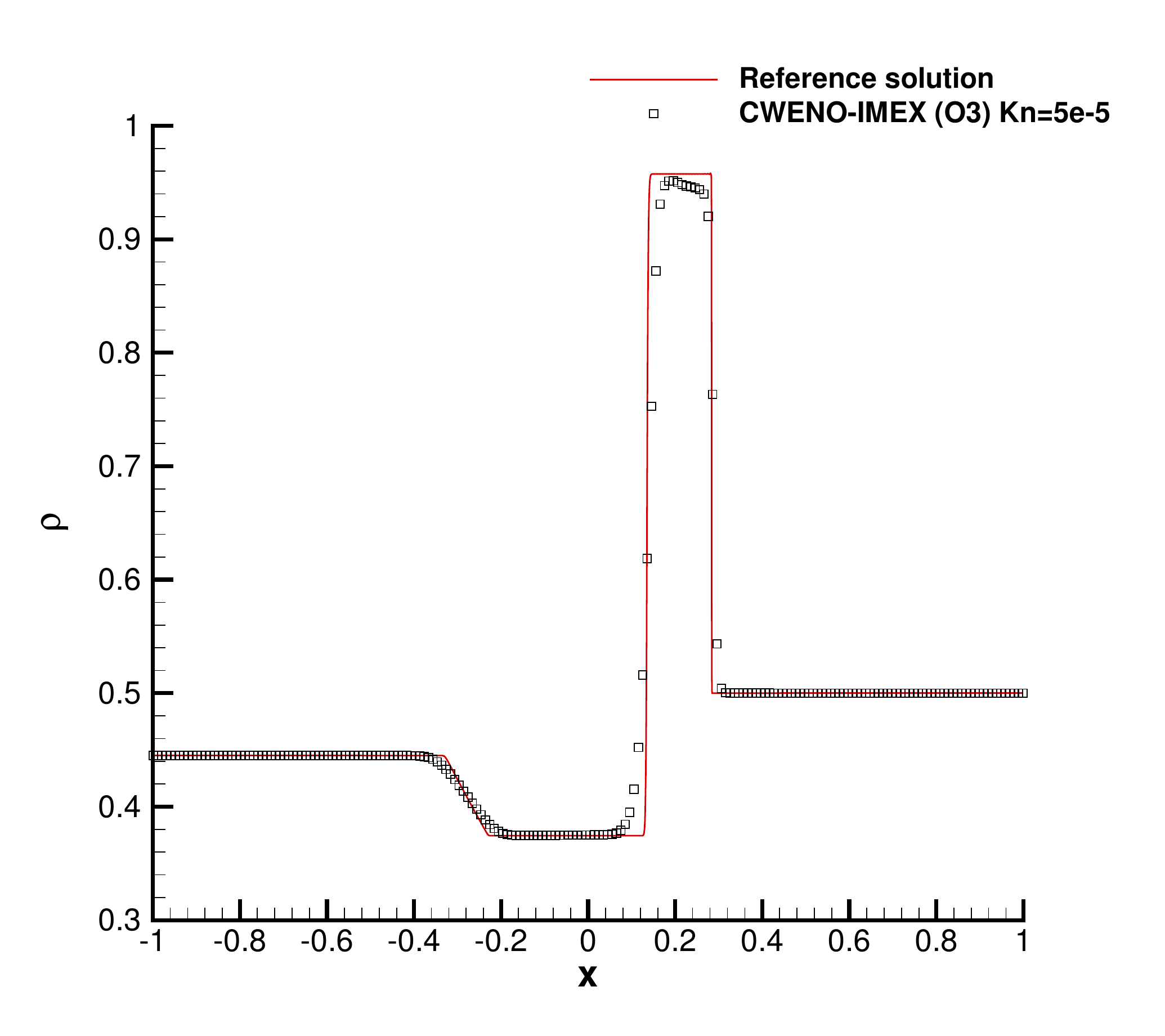}  &           
			\includegraphics[width=0.33\textwidth]{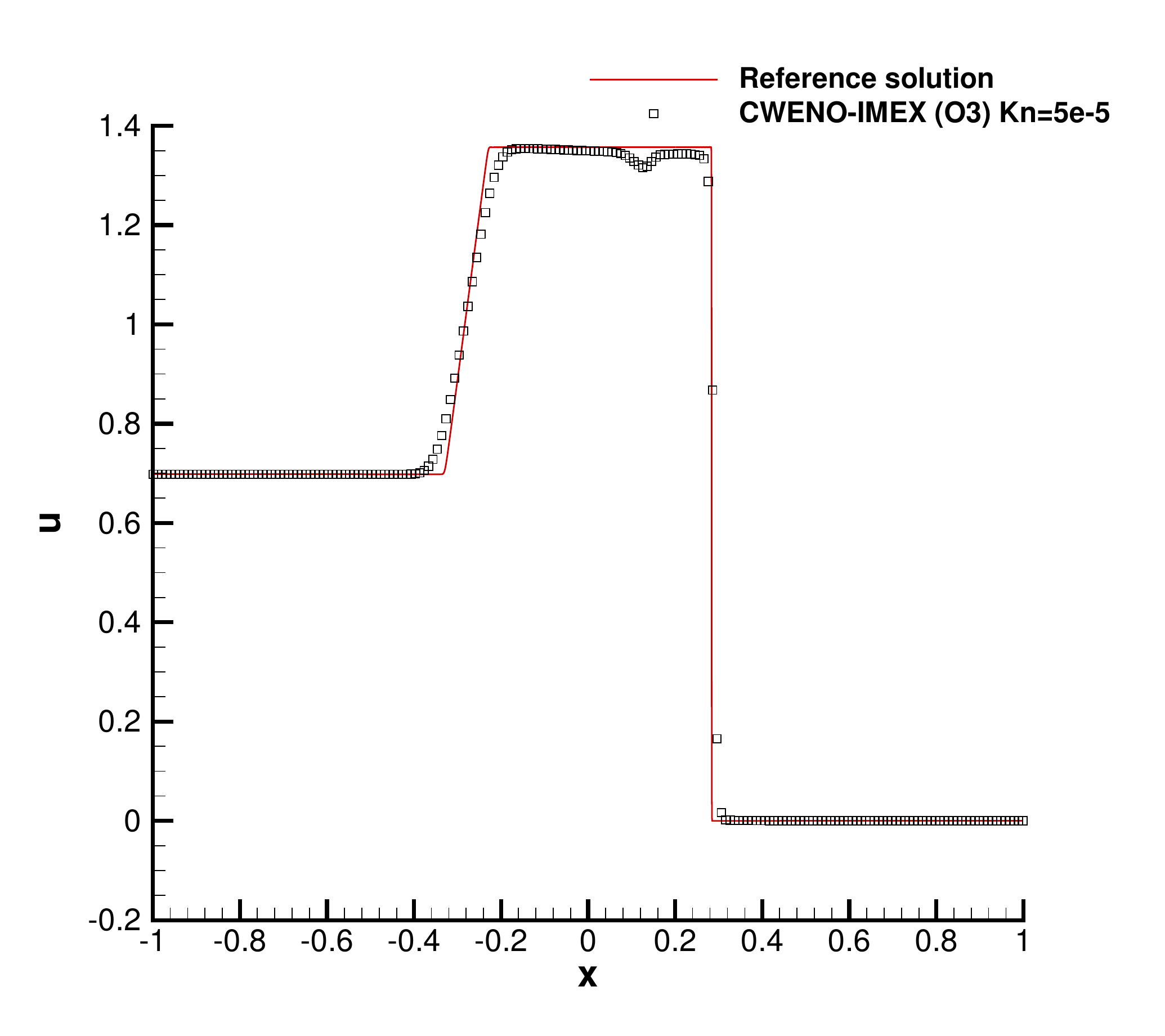} &
			\includegraphics[width=0.33\textwidth]{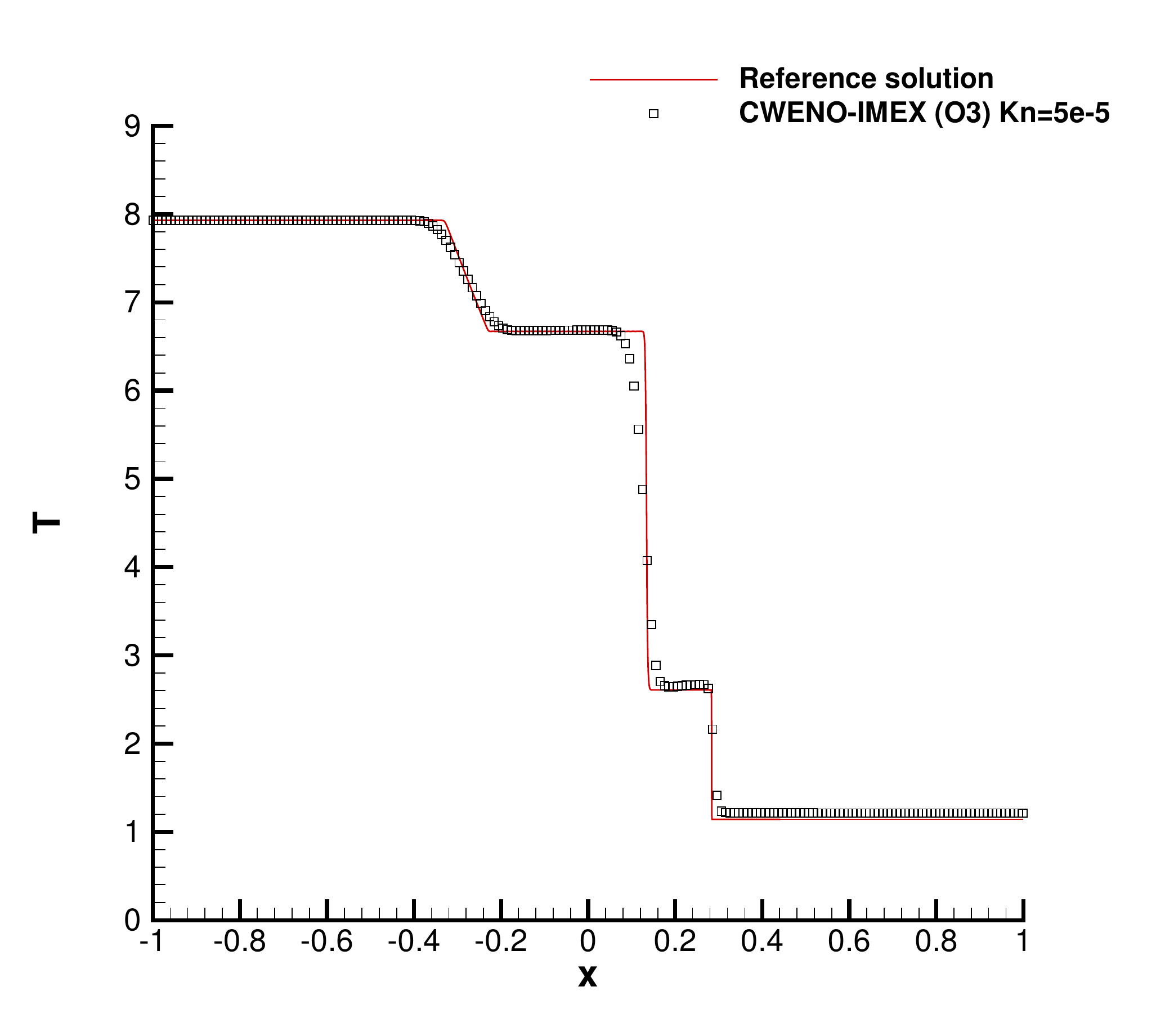}  \\ 
			\includegraphics[width=0.33\textwidth]{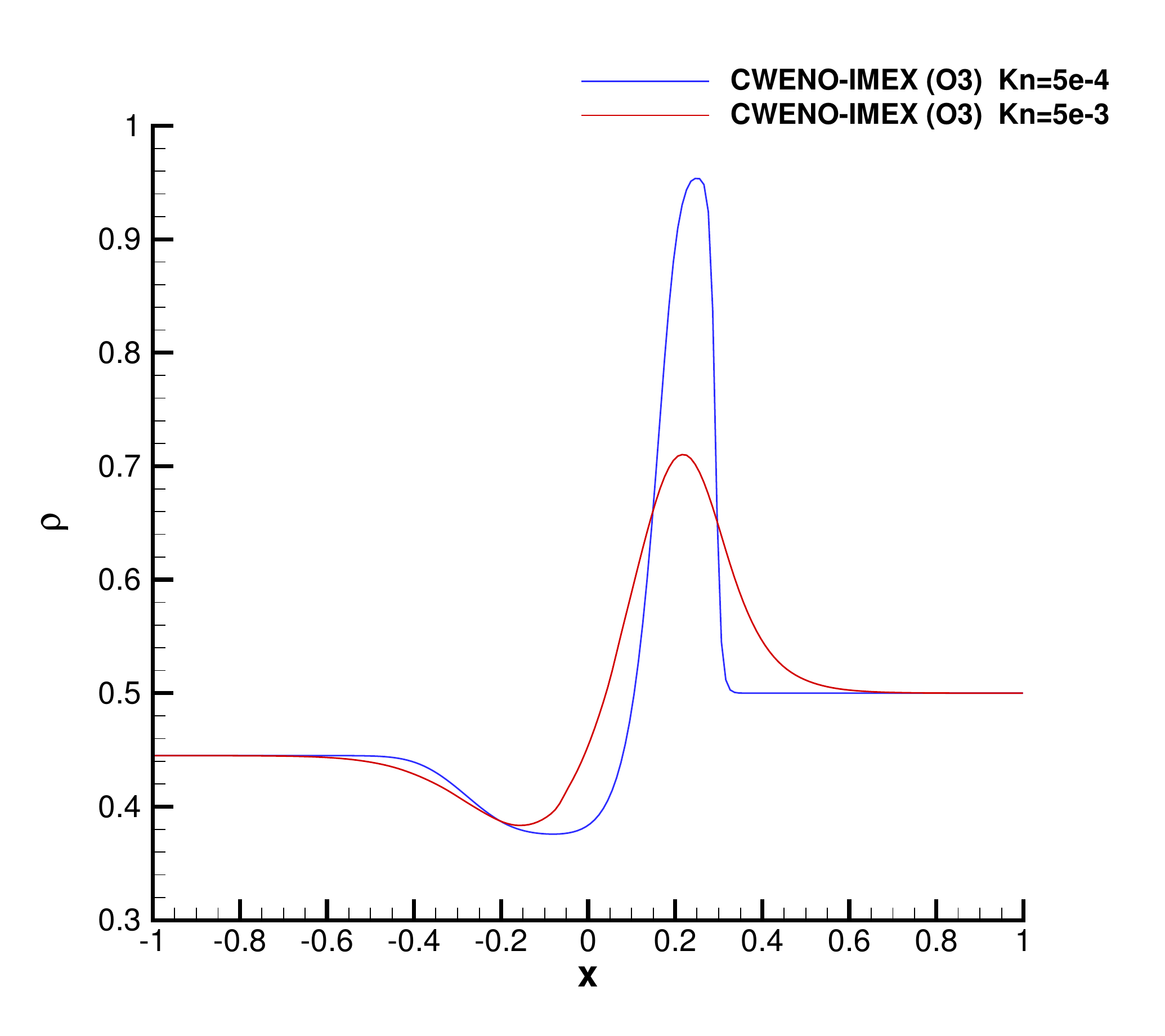}  &           
			\includegraphics[width=0.33\textwidth]{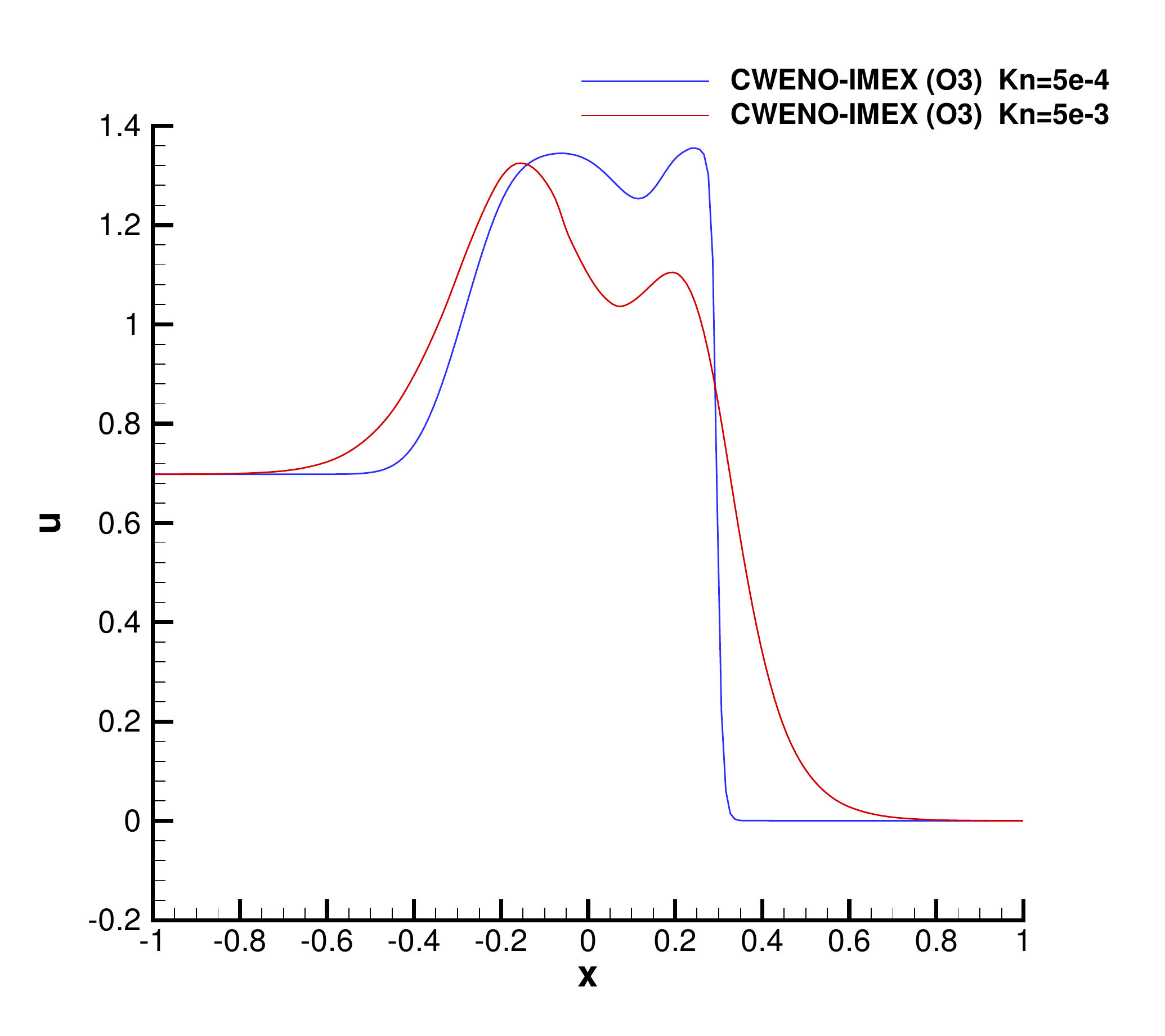} &
			\includegraphics[width=0.33\textwidth]{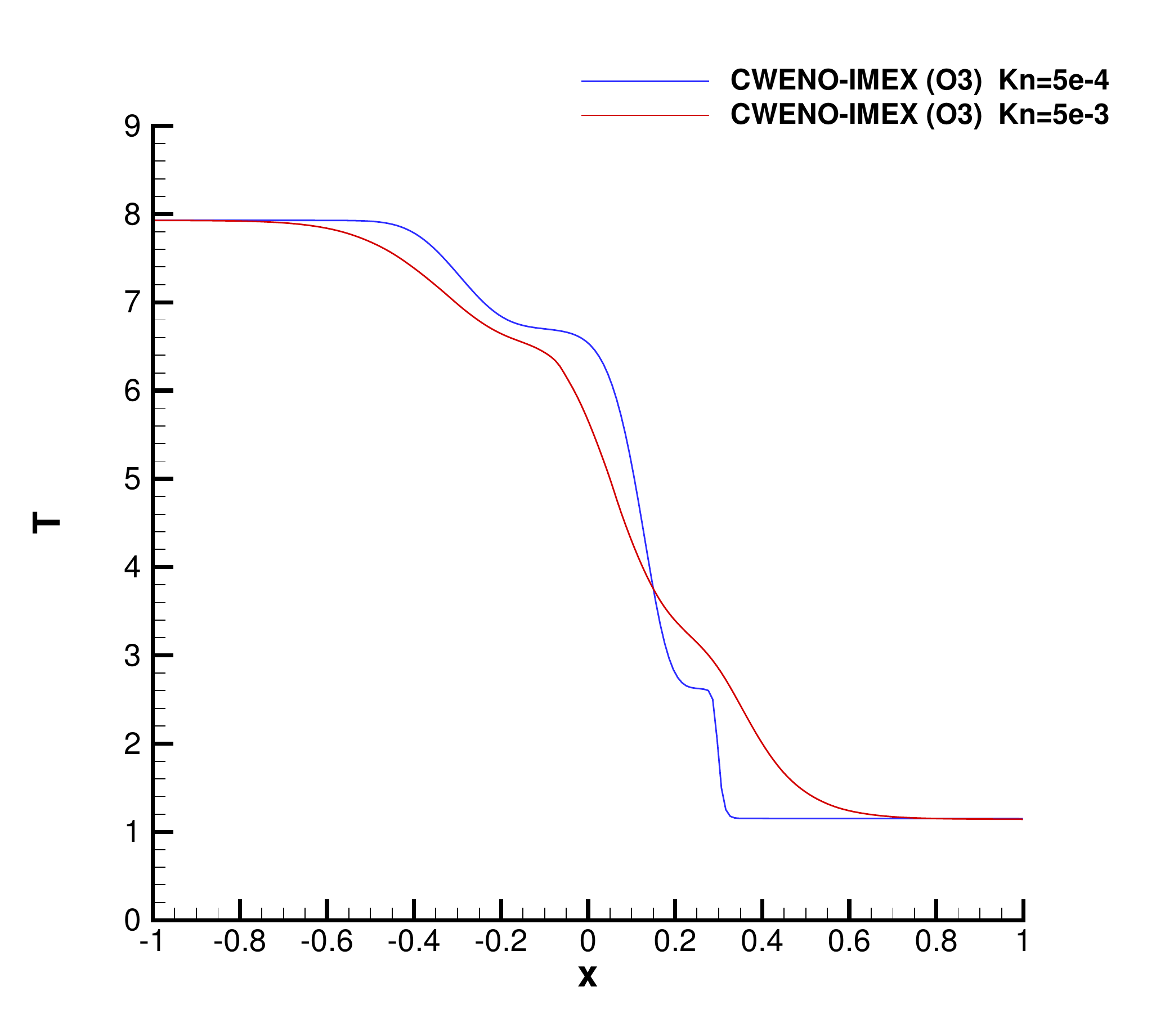}  \\  
			\includegraphics[width=0.33\textwidth]{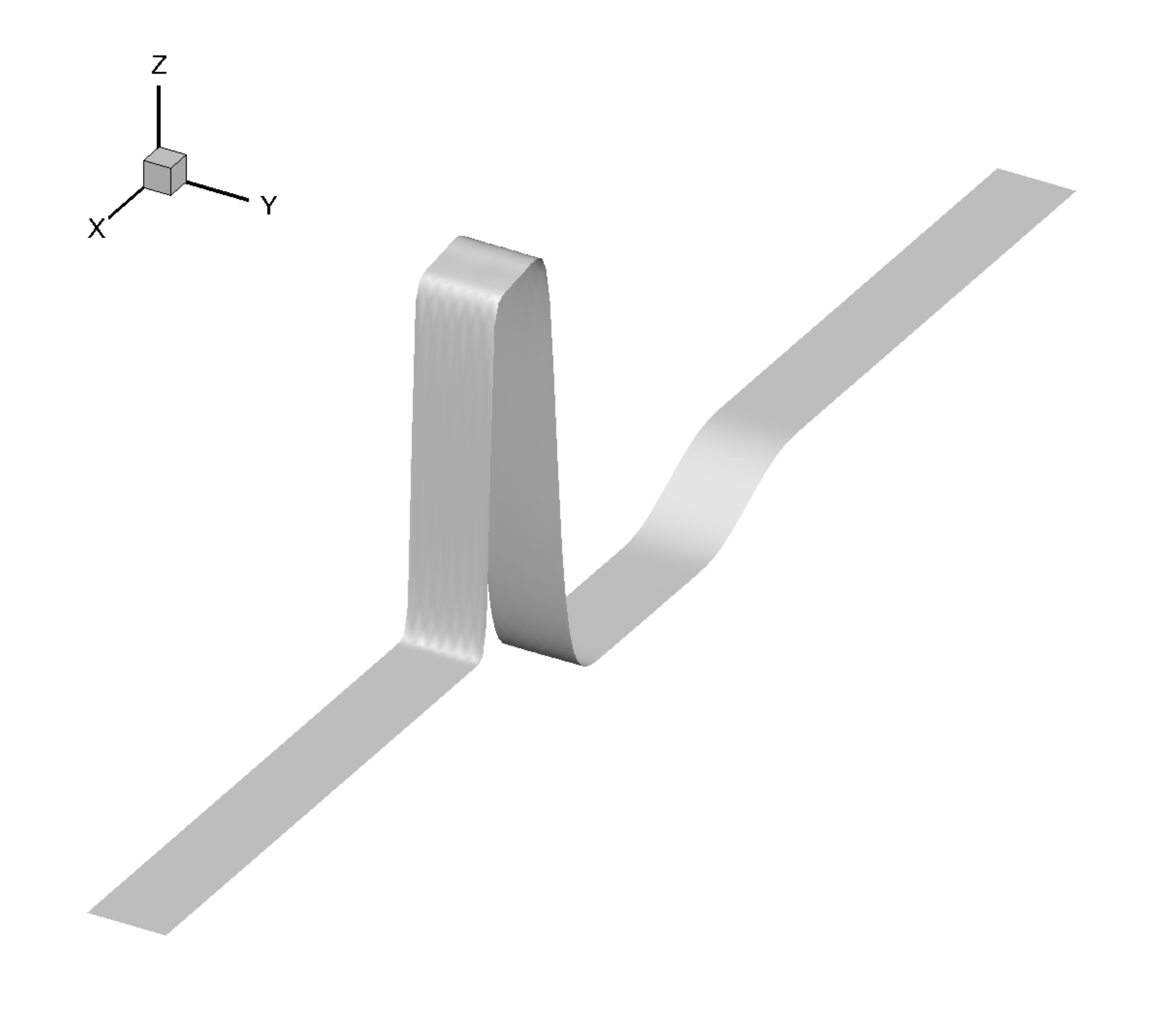}  &           
			\includegraphics[width=0.33\textwidth]{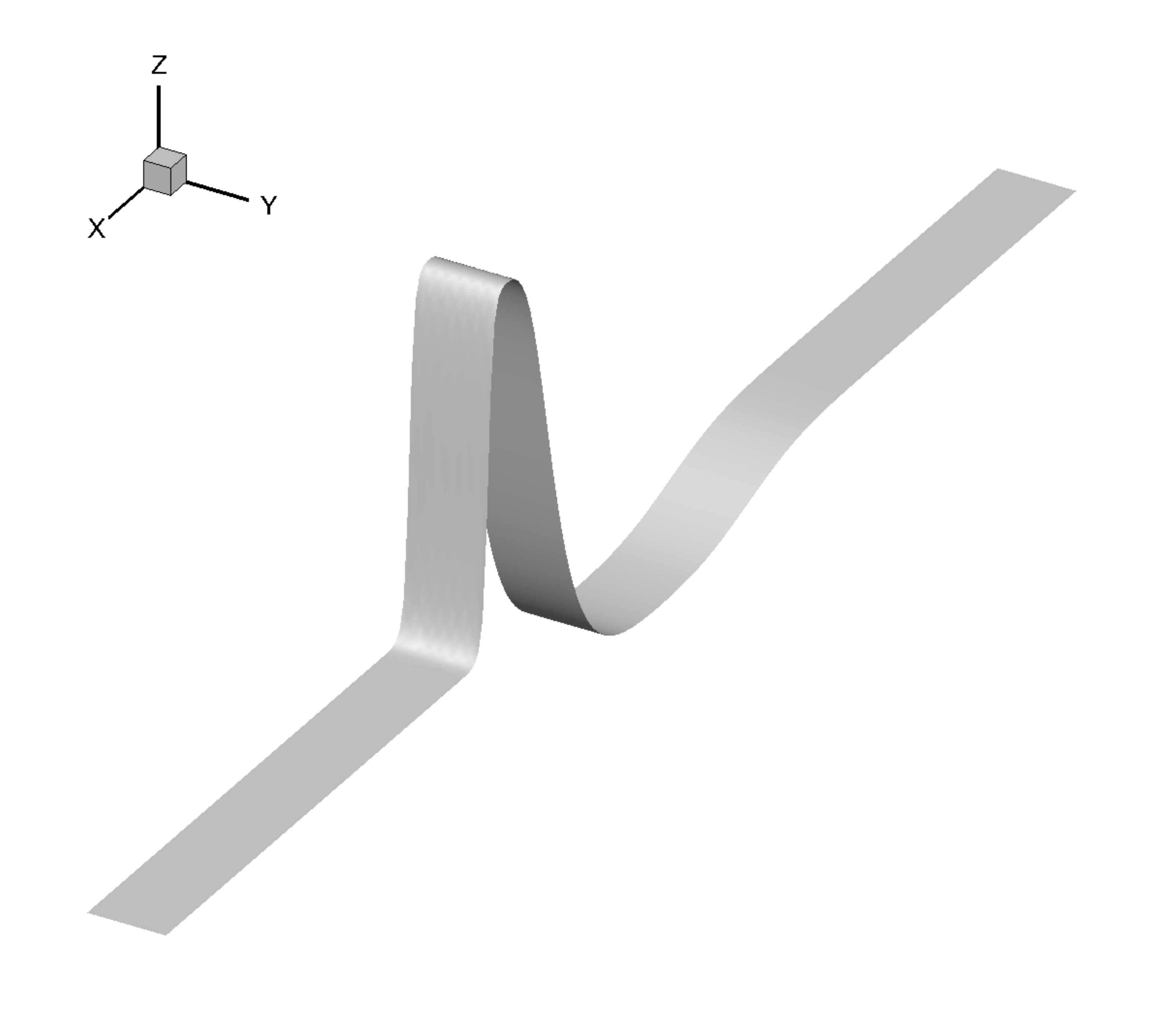} &
			\includegraphics[width=0.33\textwidth]{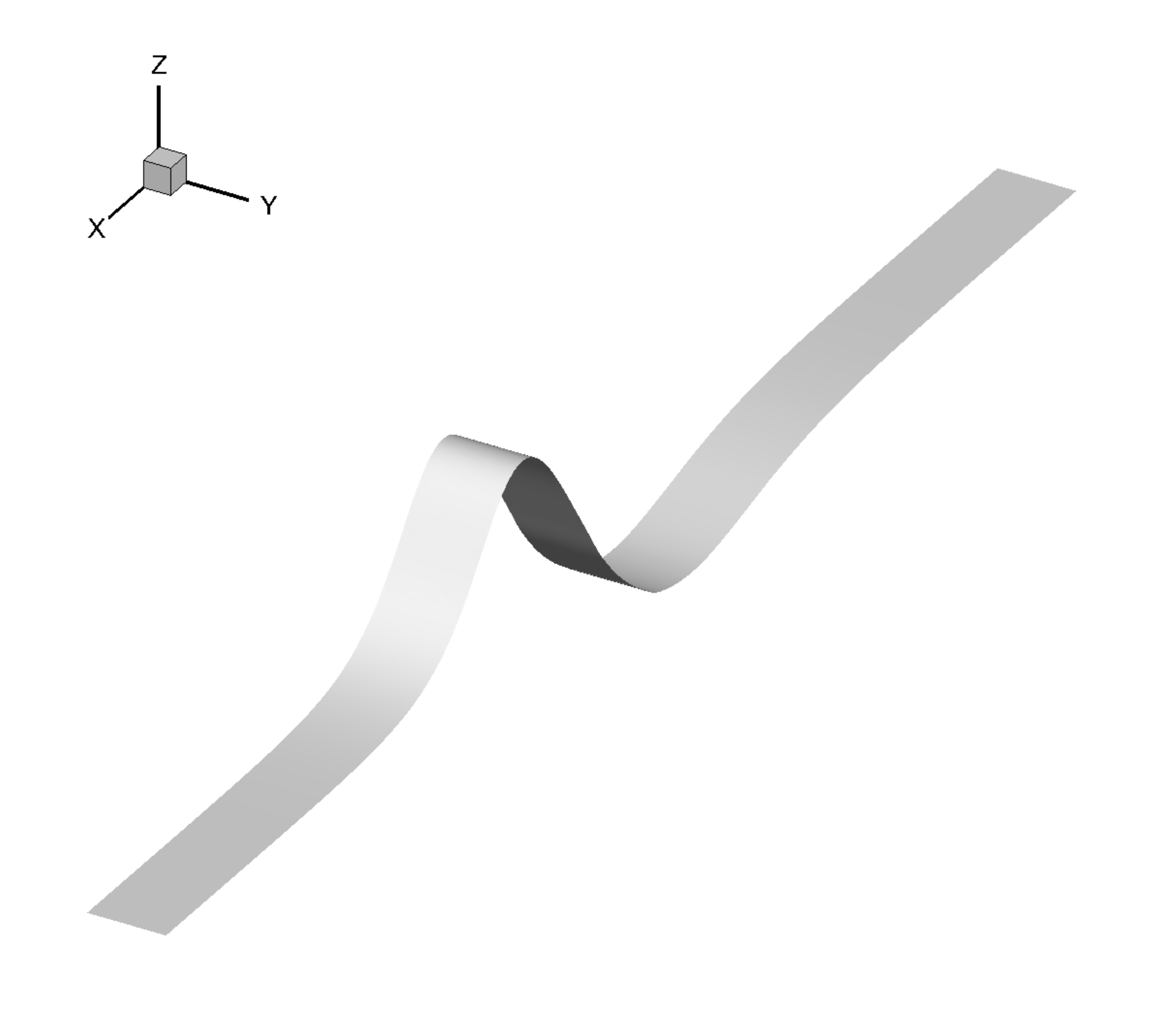}  \\           
		\end{tabular}
		\caption{Lax shock tube problem at time $t_f=0.1$ with third order CWENO-IMEX scheme for the Boltzmann model. 1D cut along the $x$-axis through the third order numerical results for density (left), horizontal velocity (middle) and temperature (right). Top: comparison between reference solution with Euler equations and numerical solution with $\varepsilon=5 \cdot 10^{-5}$. Middle: comparison between numerical solution with $\varepsilon=5 \cdot 10^{-4}$ (blue solid line) and $\varepsilon=5 \cdot 10^{-3}$ (red solid line). Bottom: three-dimensional view of density profile for the Boltzmann model with $\varepsilon=5 \cdot 10^{-5}$ (left), $\varepsilon=5 \cdot 10^{-4}$ (middle) and $\varepsilon=5 \cdot 10^{-3}$ (right). }
		\label{fig.Lax-Boltz-Kn}
	\end{center}
\end{figure}

We consider now a classical Lax shock tube Riemann problem. The essential feature one can observe in rarefied flows with respect to dense fluids is the presence of a physical diffusion which mitigates the waves making the solution smoother. From the numerical side, another important aspect consists in analysing the capability of the methods working on unstructured meshes to keep the one dimensional structure of the solution even if the element edges are not aligned with the fluid motion. The computational domain is the box $\Omega=[-1;1] \times [-0.05;0.05]$ discretized with $[200 \times 10]$ triangular control volumes with characteristic mesh size of $h(\Omega)=0.01$. Periodic boundary conditions are set in $y-$direction, while Dirichlet boundaries are imposed in the $x-$direction. The velocity space $\mathcal{V}=[-15;15]\times [-15;15]$ counts a total number of $32^2=1024$ regular Cartesian control volumes. The initial condition is 
\begin{equation}
\mathcal{U}_L = \left(0.445,0.698,0,7.928 \right), \qquad \mathcal{U}_R = \left(0.5,0,0,1.142\right),
\end{equation}
with the two states separated at $x=0$. Figure \ref{fig.Lax-Boltz-Kn} shows a one dimensional cut along the $x-$axis for the density, the horizontal velocity and the temperature profiles for the Boltzmann model with $\varepsilon=10^{-5}$ at the top and for $\varepsilon=10^{-4}$ and $\varepsilon=10^{-3}$ in the middle at the final time of the simulation $t_f=0.1$. The third order CWENO-IMEX method is employed using $200$ cells in space.
The numerical results are compared against the exact solution of the Euler equations of compressible gas dynamics, highlighting the correct behavior of our scheme as $\varepsilon \to 0$. A comparison between our approach and a Direct Simulation Monte Carlo method (DSMC) solving the Boltzmann equation \eqref{eq:Boltzmann} is proposed in the case in which $\varepsilon=0.005$ and $\varepsilon=5\cdot 10^{-4}$. This method works on a regular Cartesian mesh of $250^2$ elements and employs around 62 millions of particles. The details can be found in \cite{pareschi_dsmc}. Also in this case, we observe a good matching between our method and the DSMC solution. Finally, a three-dimensional view is shown at the bottom for the density profile at different Knudsen numbers highlighting the capability of the method to keep the one dimensional structure of the solution in presence of an unstructured mesh.
\begin{figure}[!htbp]
	\begin{center}
		\begin{tabular}{ccc} 
			\includegraphics[width=0.33\textwidth]{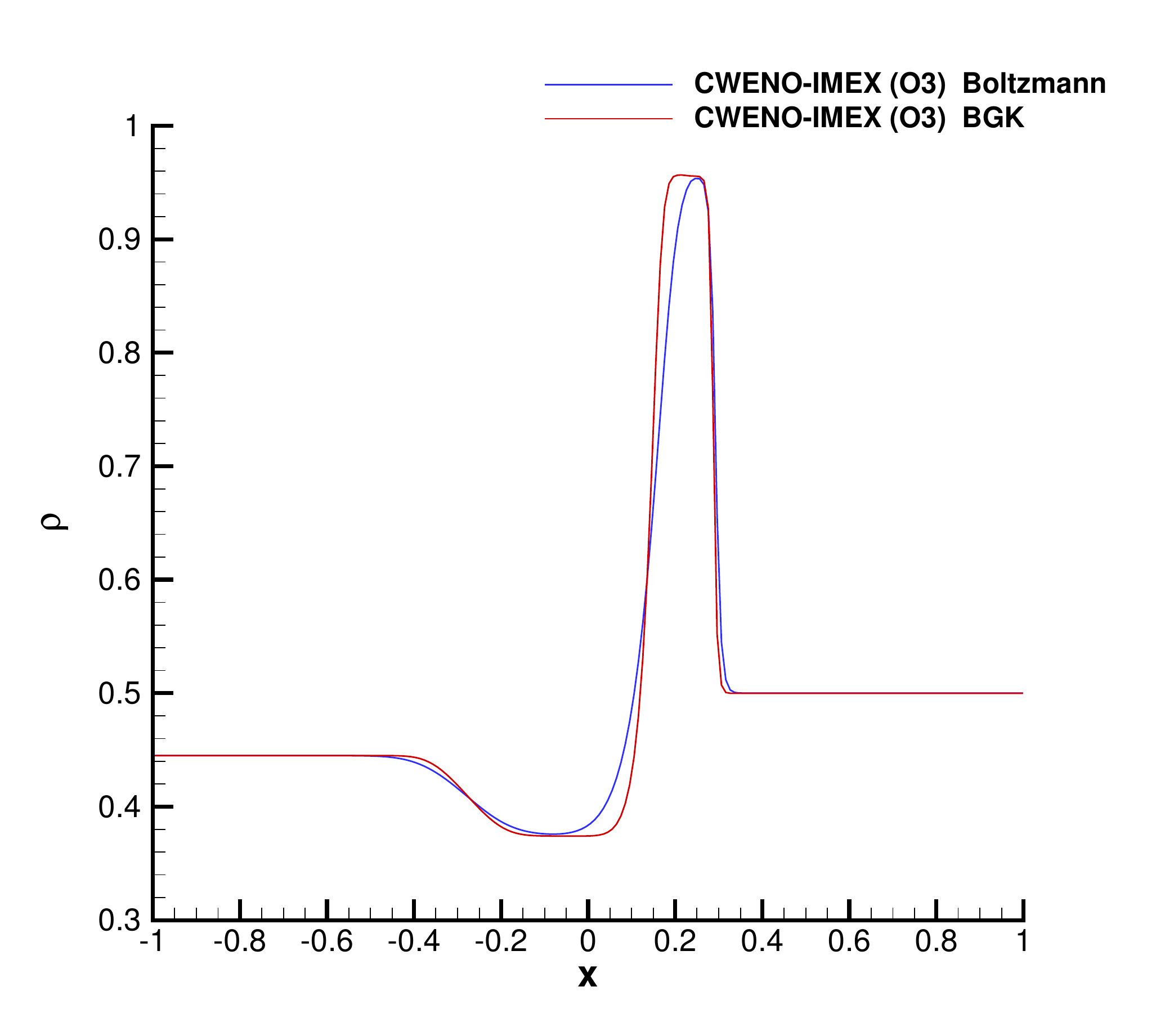}  &           
			\includegraphics[width=0.33\textwidth]{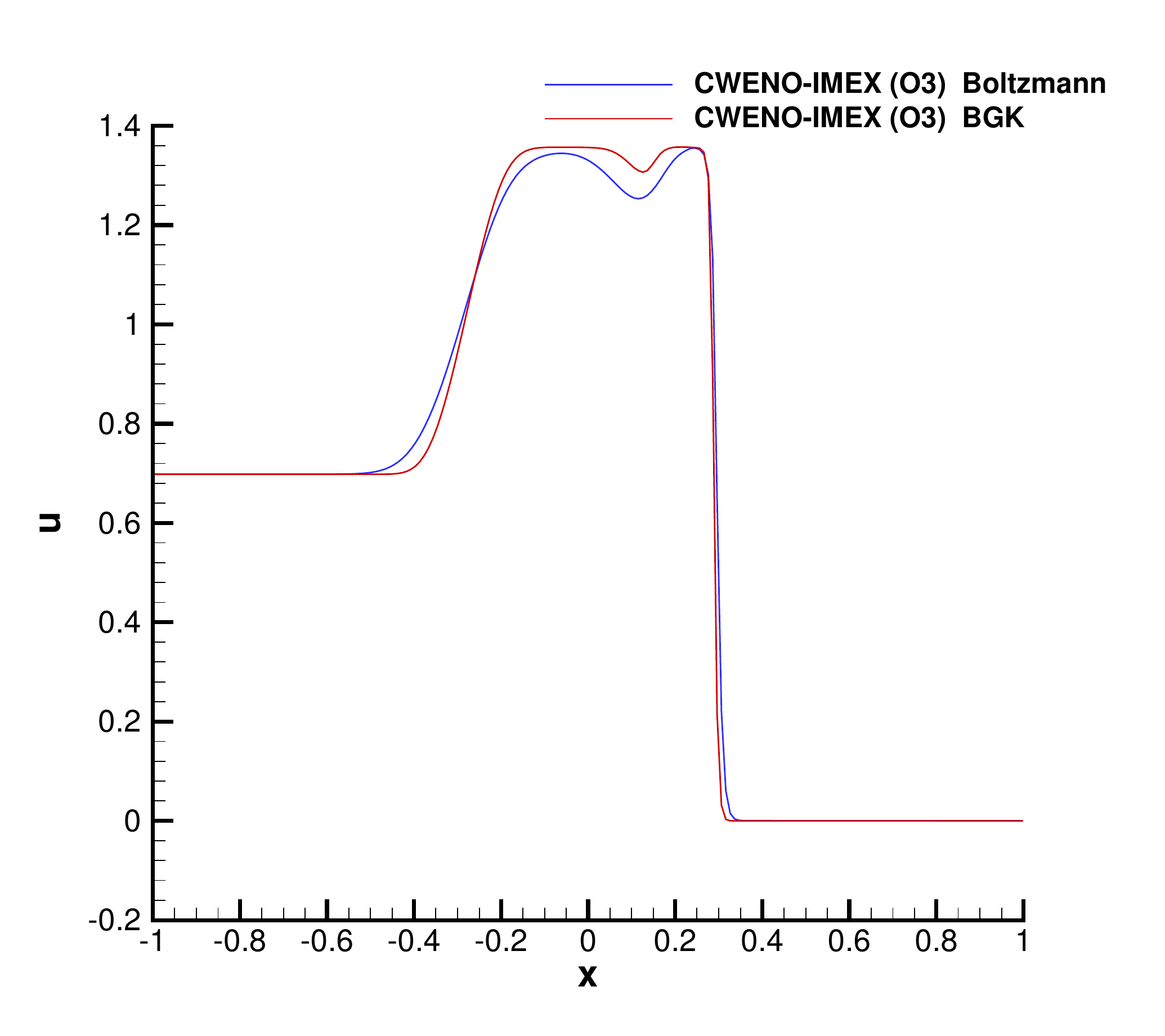} &
			\includegraphics[width=0.33\textwidth]{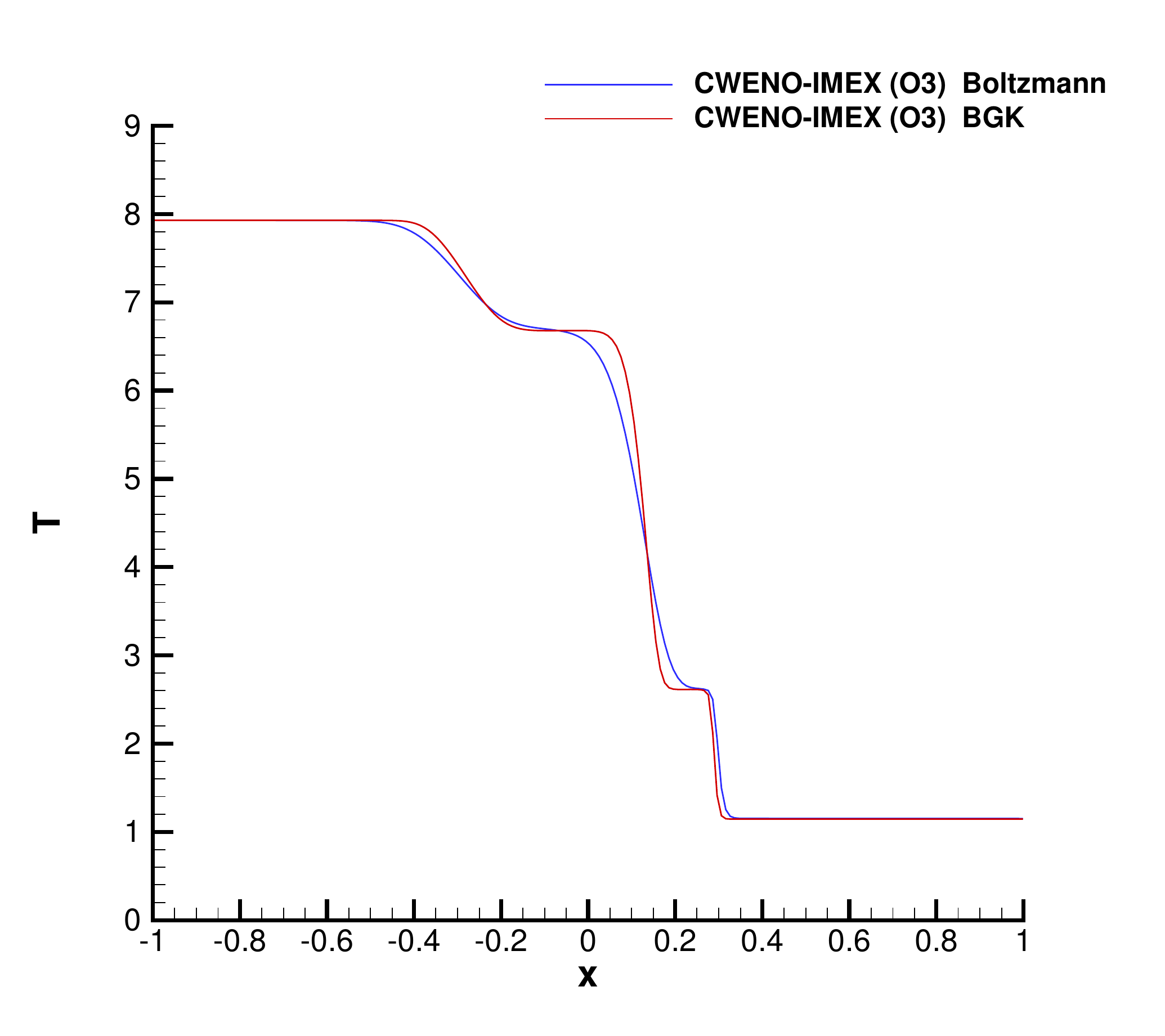}  \\  
			\includegraphics[width=0.33\textwidth]{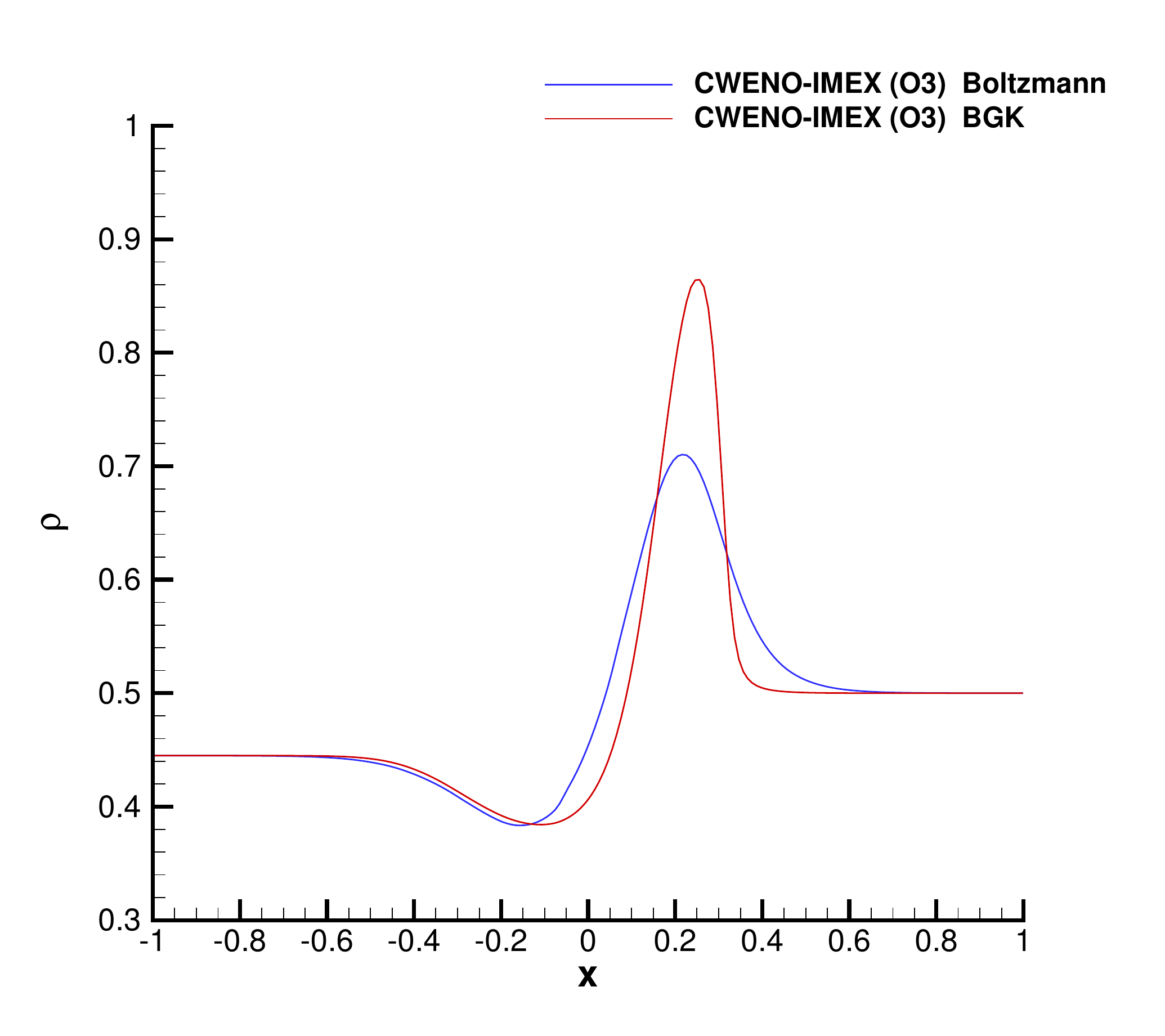}  &           
			\includegraphics[width=0.33\textwidth]{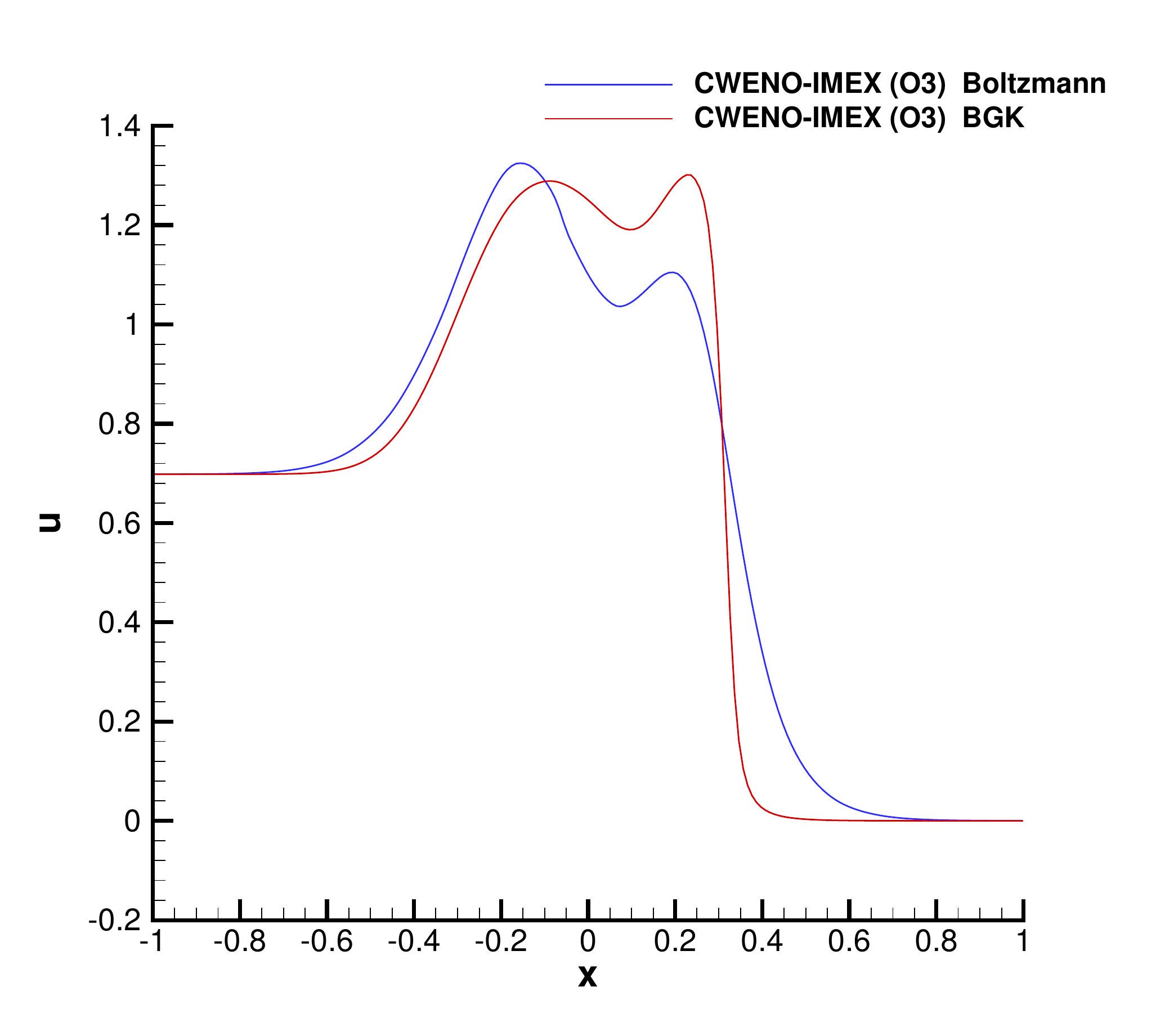} &
			\includegraphics[width=0.33\textwidth]{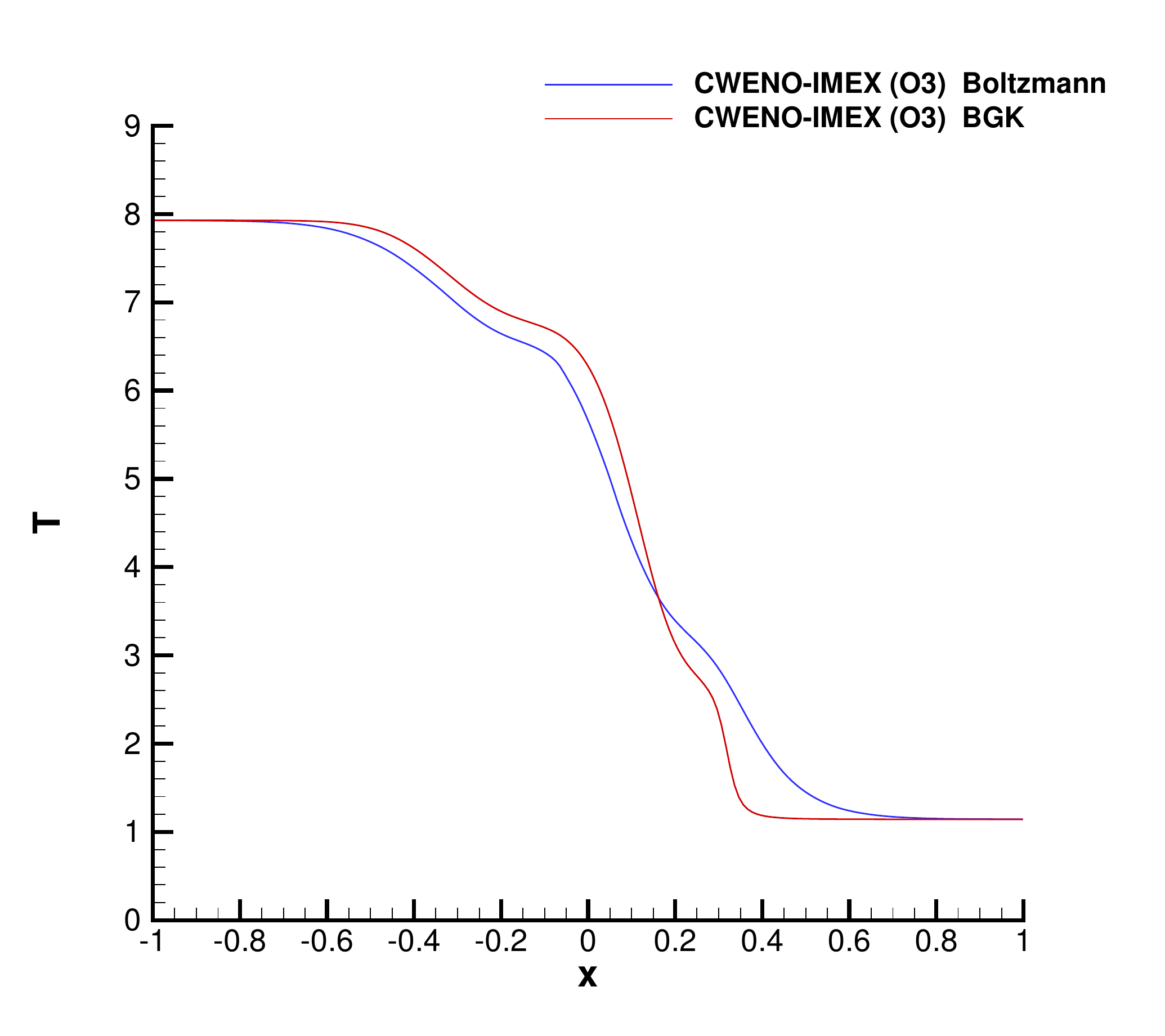}  \\           
		\end{tabular}
		\caption{Lax shock tube problem with $\varepsilon=5 \cdot 10^{-4}$ (top row) and $\varepsilon=5 \cdot 10^{-3}$ (bottom row) for Boltzmann (blue solid line) and BGK (red solid line) model at time $t_f=0.1$. 1D cut along the $x$-axis through the third order numerical results for density (left), horizontal velocity (middle) and temperature (right).}
		\label{fig.Lax-BoltzBGK}
	\end{center}
\end{figure}
Figure \ref{fig.Lax-BoltzBGK} instead shows a comparison between the Boltzmann and the BGK models for different values of the Knudsen number. As expected \cite{ACTA} the BGK model tends to over relax the solution towards the fluid limit. This is particularly clear when the Knudsen number is quite large $\varepsilon=5\cdot 10^{-3}$ and it is in line with the theoretical observations which claim that the BGK model is adapted for the description of rarefied flows only when the collisional rate is such that the fluid is close to equilibrium.

\begin{figure}[!htbp]
	\begin{center}
		\begin{tabular}{ccc} 
			\includegraphics[width=0.33\textwidth]{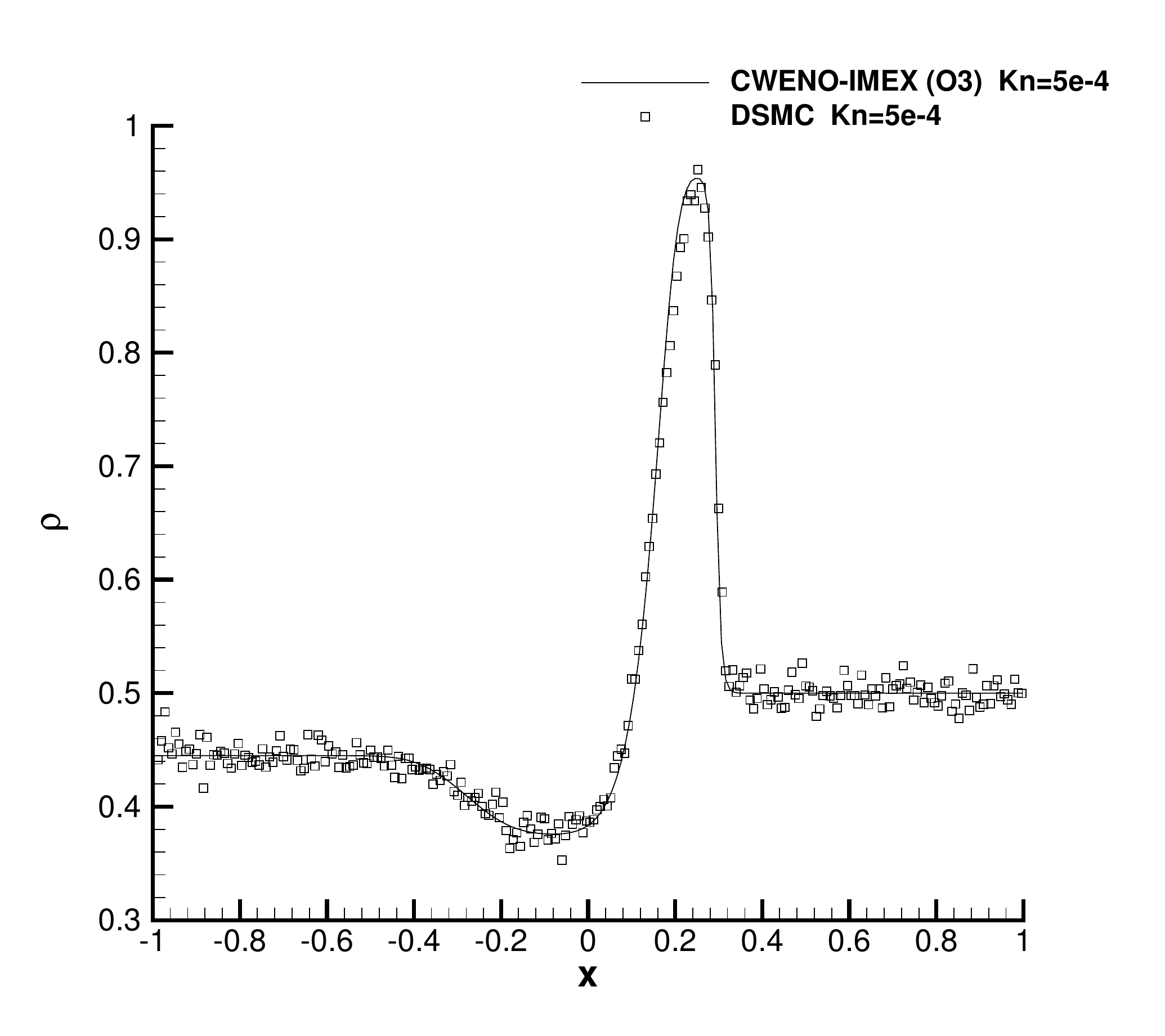}  &           
			\includegraphics[width=0.33\textwidth]{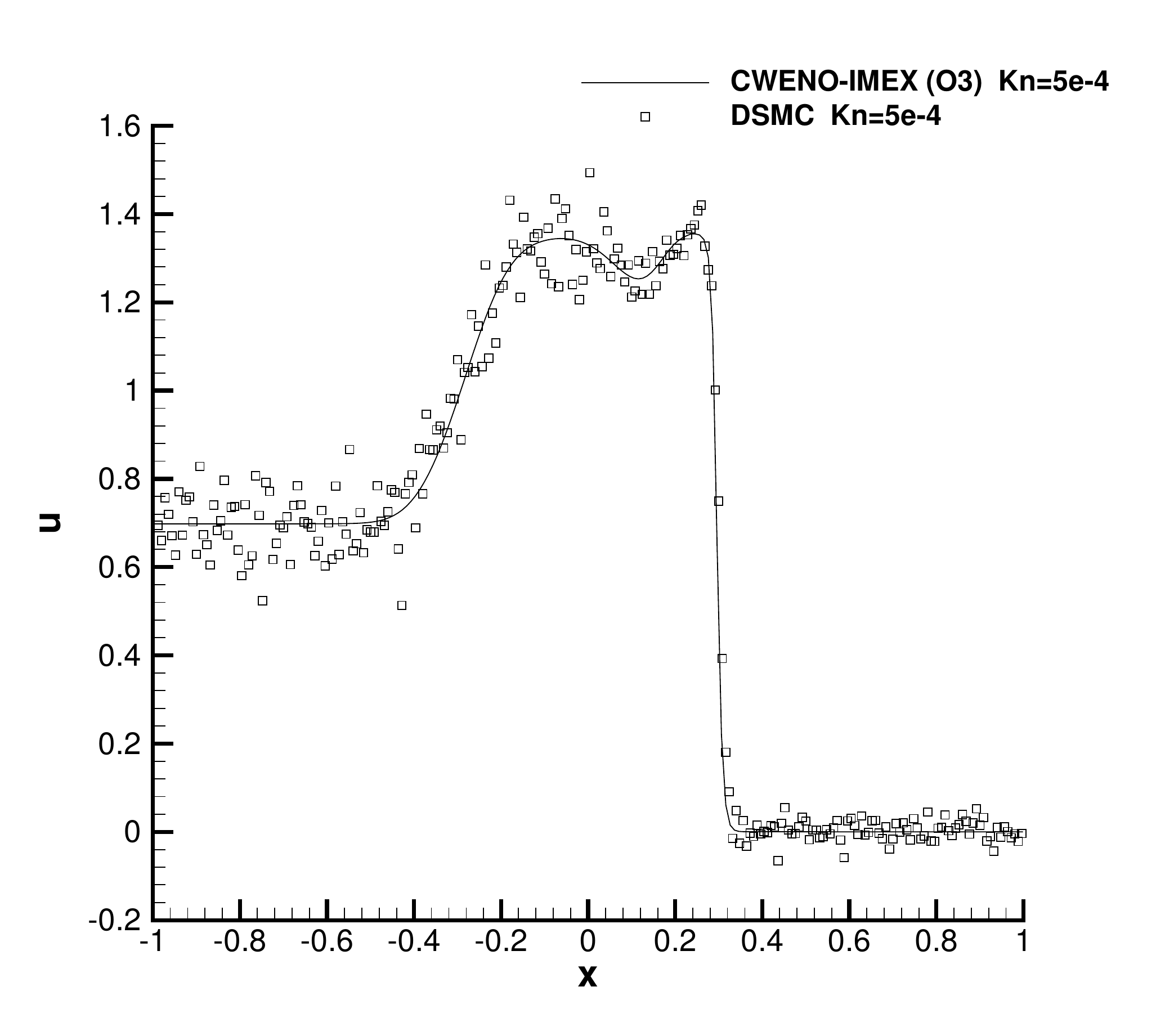} &
			\includegraphics[width=0.33\textwidth]{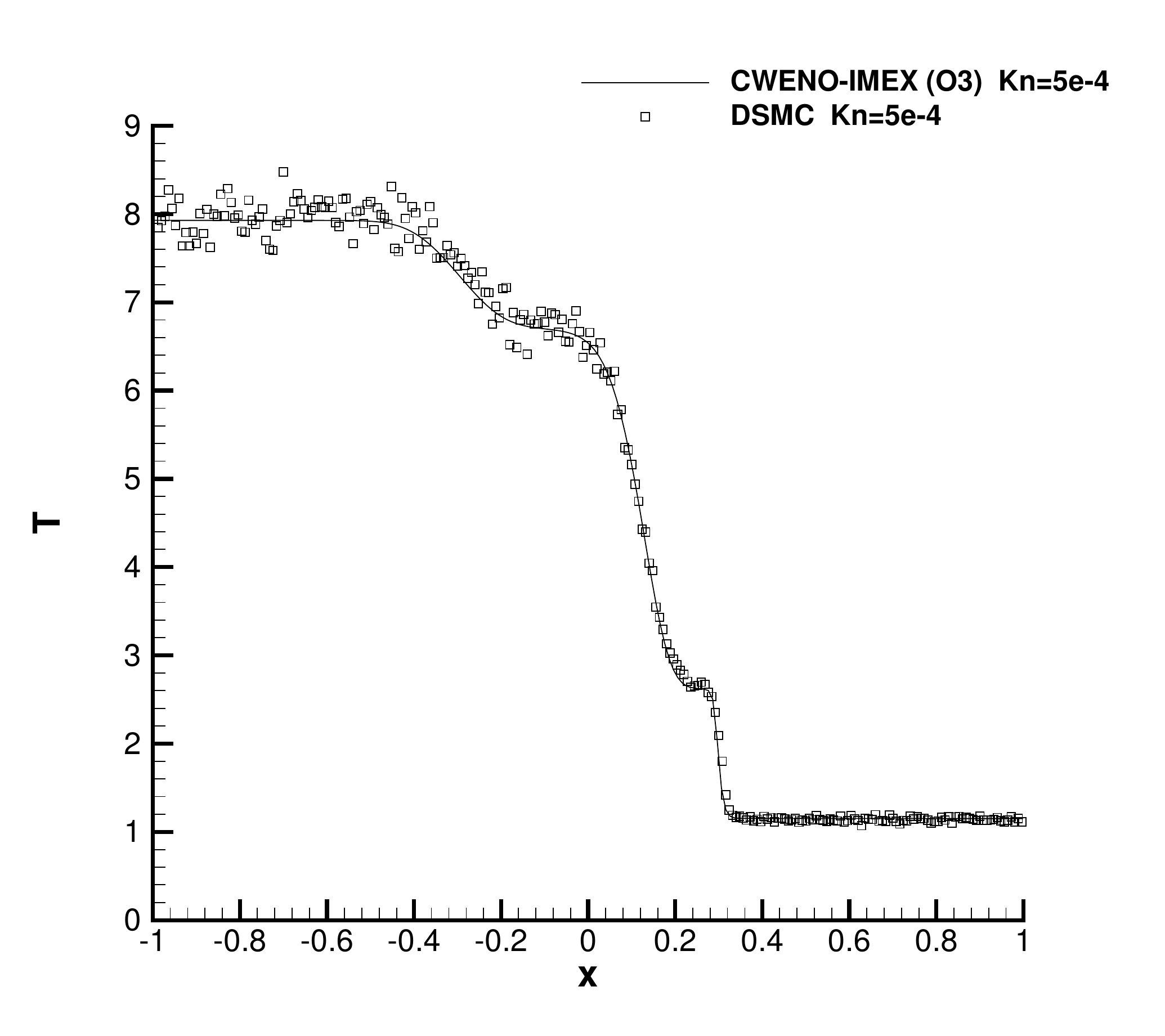}  \\  
			\includegraphics[width=0.33\textwidth]{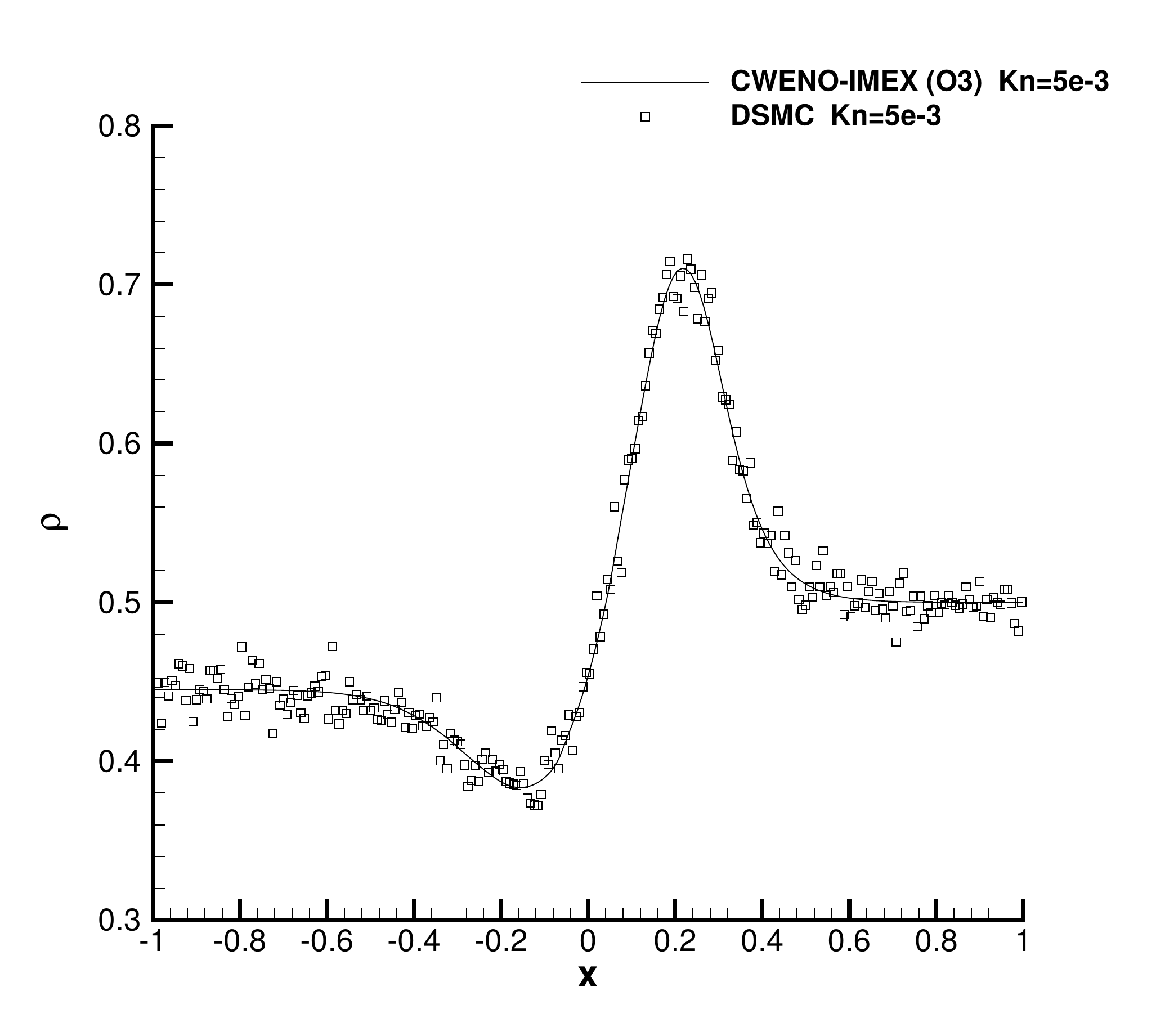}  &           
			\includegraphics[width=0.33\textwidth]{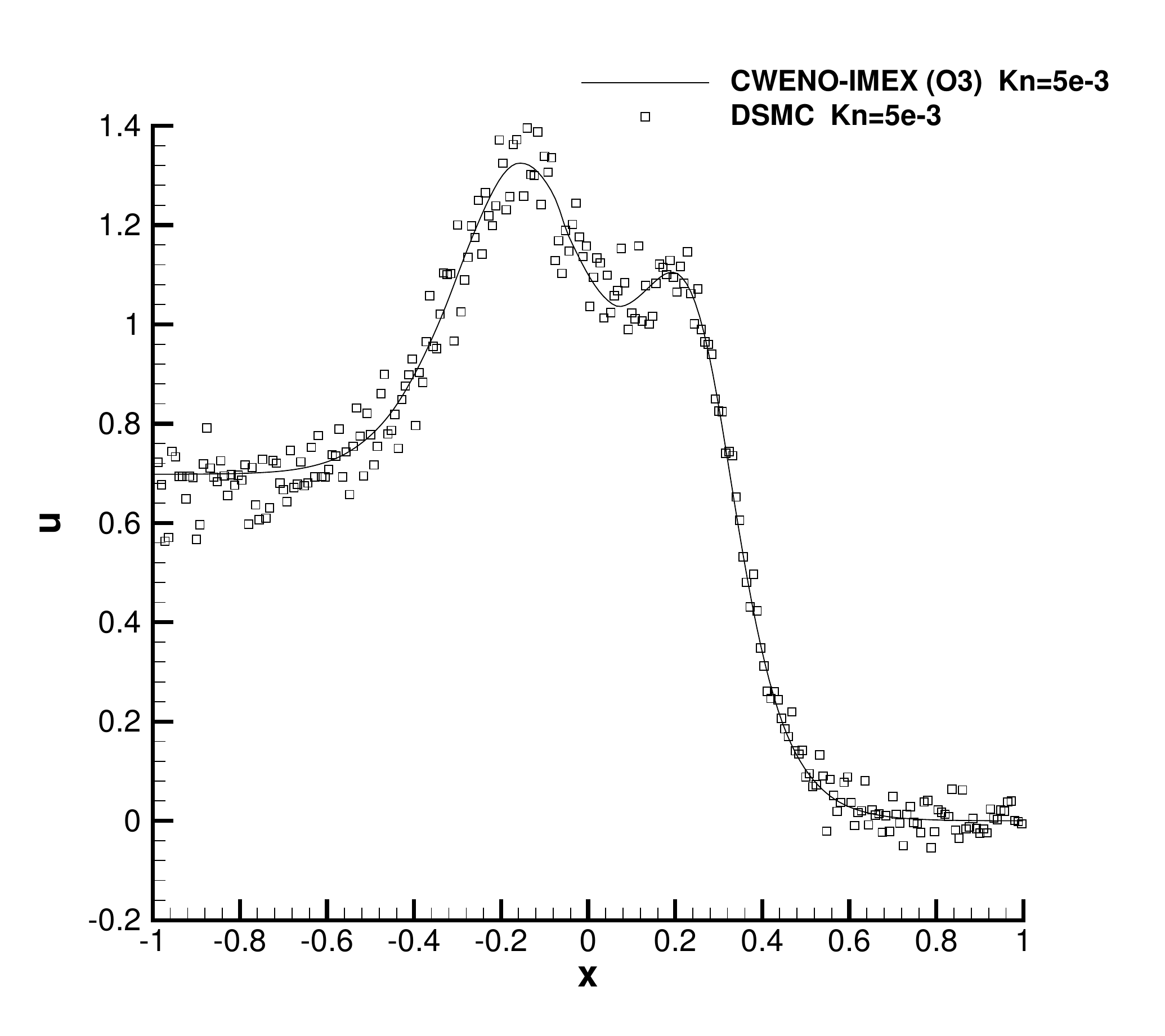} &
			\includegraphics[width=0.33\textwidth]{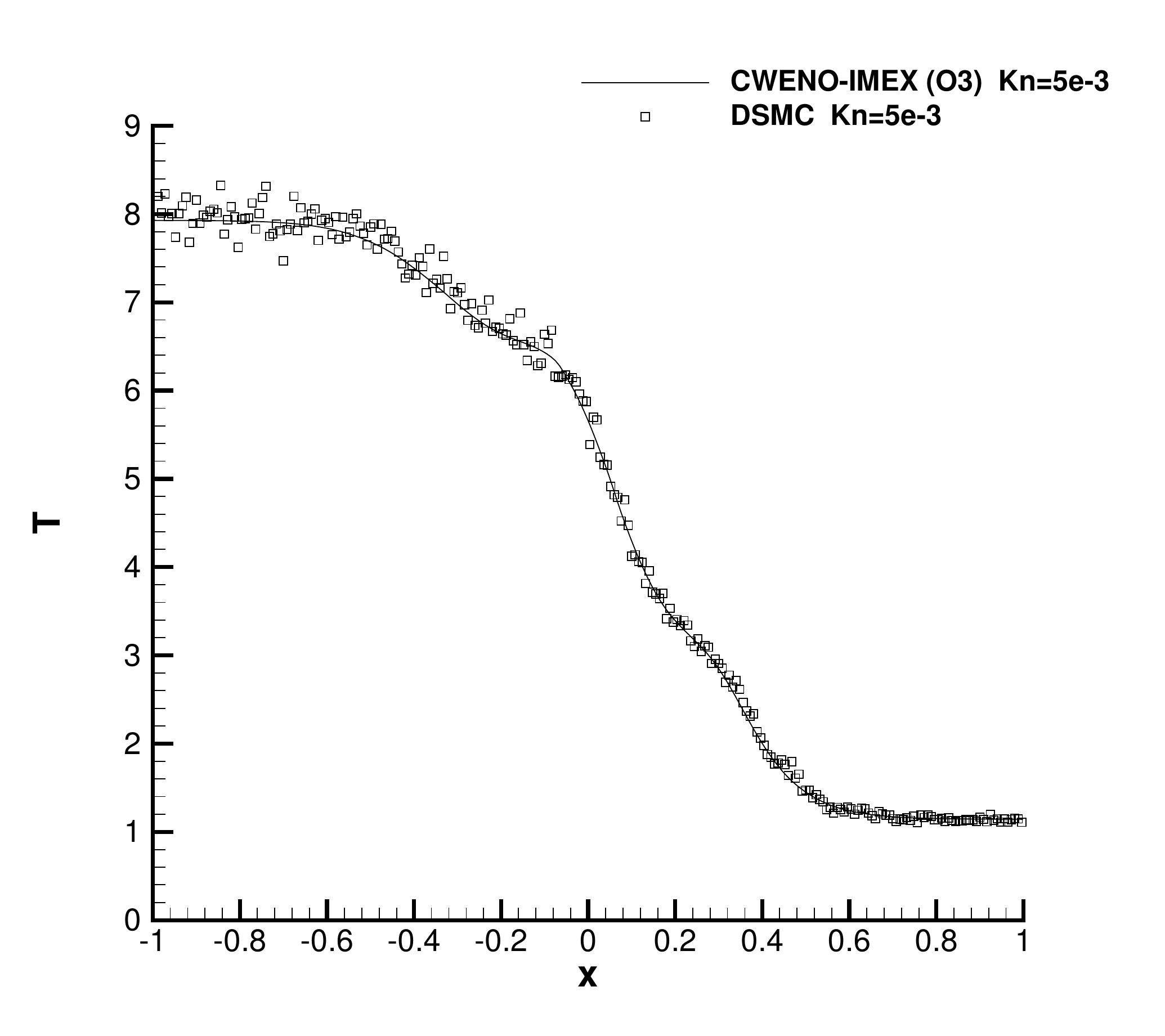}
		\end{tabular}
		\caption{Lax shock tube problem with $\varepsilon=5 \cdot 10^{-4}$ (top row) and $\varepsilon=5 \cdot 10^{-3}$ (bottom row) at time $t_f=0.1$. 1D cut along the $x$-axis for density (left), horizontal velocity (middle) and temperature (right). The solid line represents the solution obtained with the third order CWENO-IMEX schemes, while scatter plots represent the results computed with a DSMC simulation on a grid composed of $250^2$ Cartesian cells with approximately $62 \cdot 10^6$ particles.}
		\label{fig.Lax-BoltzDSMC}
	\end{center}
\end{figure}

Finally, the third order results for the Boltzmann model are compared against the solution obtained by employing a DSMC method in Figure \ref{fig.Lax-BoltzDSMC}. For the DSMC simulation a Cartesian grid of $250^2$ elements is used with around 62 millions of particles. An overall very good agreement can be appreciated between the solvers, which confirms what already observed for the case of the fluid limit.

%
\subsection{Explosion problem}
\begin{figure}[!htbp]
	\begin{center}
		\begin{tabular}{cc} 
			\includegraphics[width=0.47\textwidth]{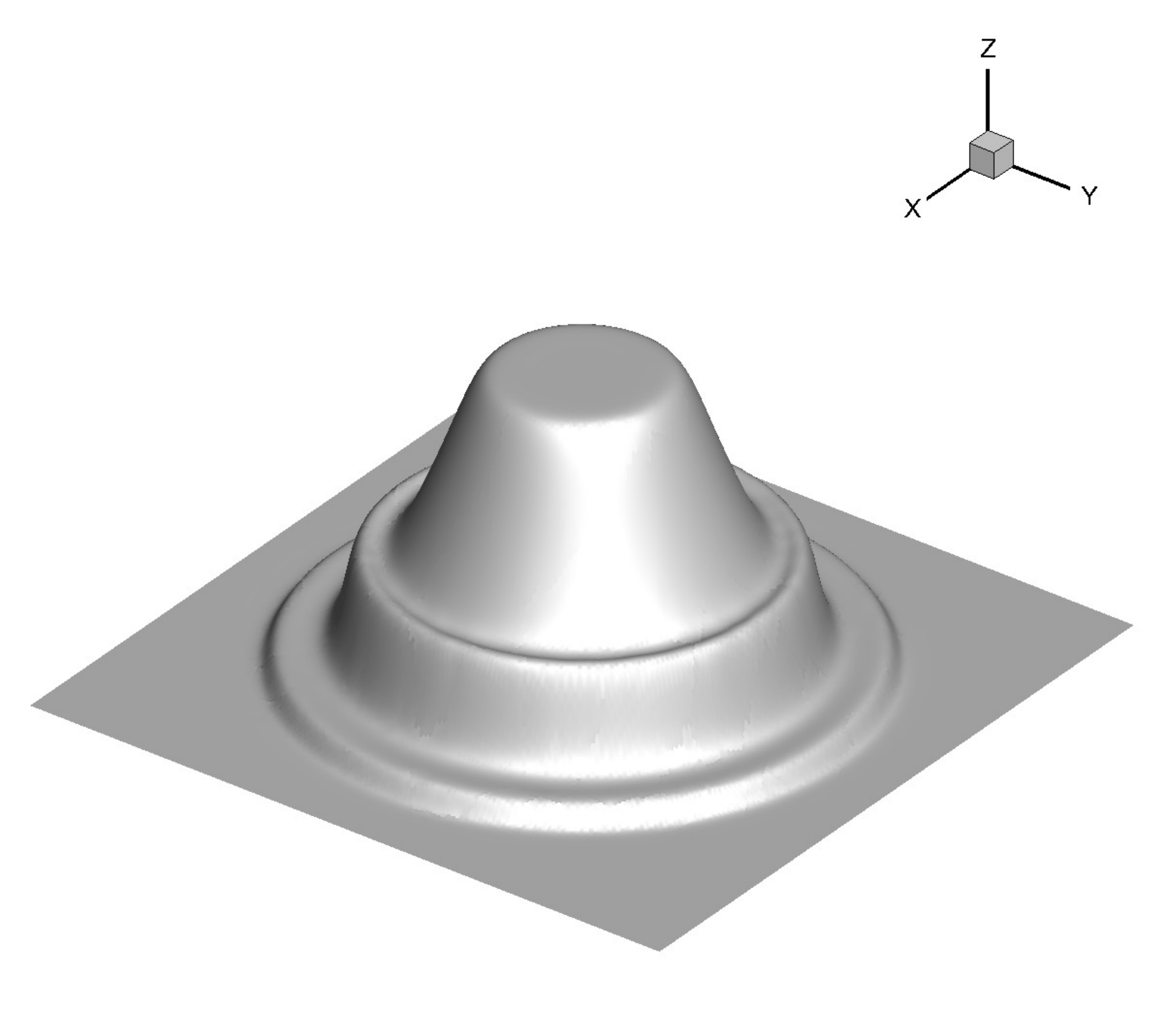}  &           
			\includegraphics[width=0.47\textwidth]{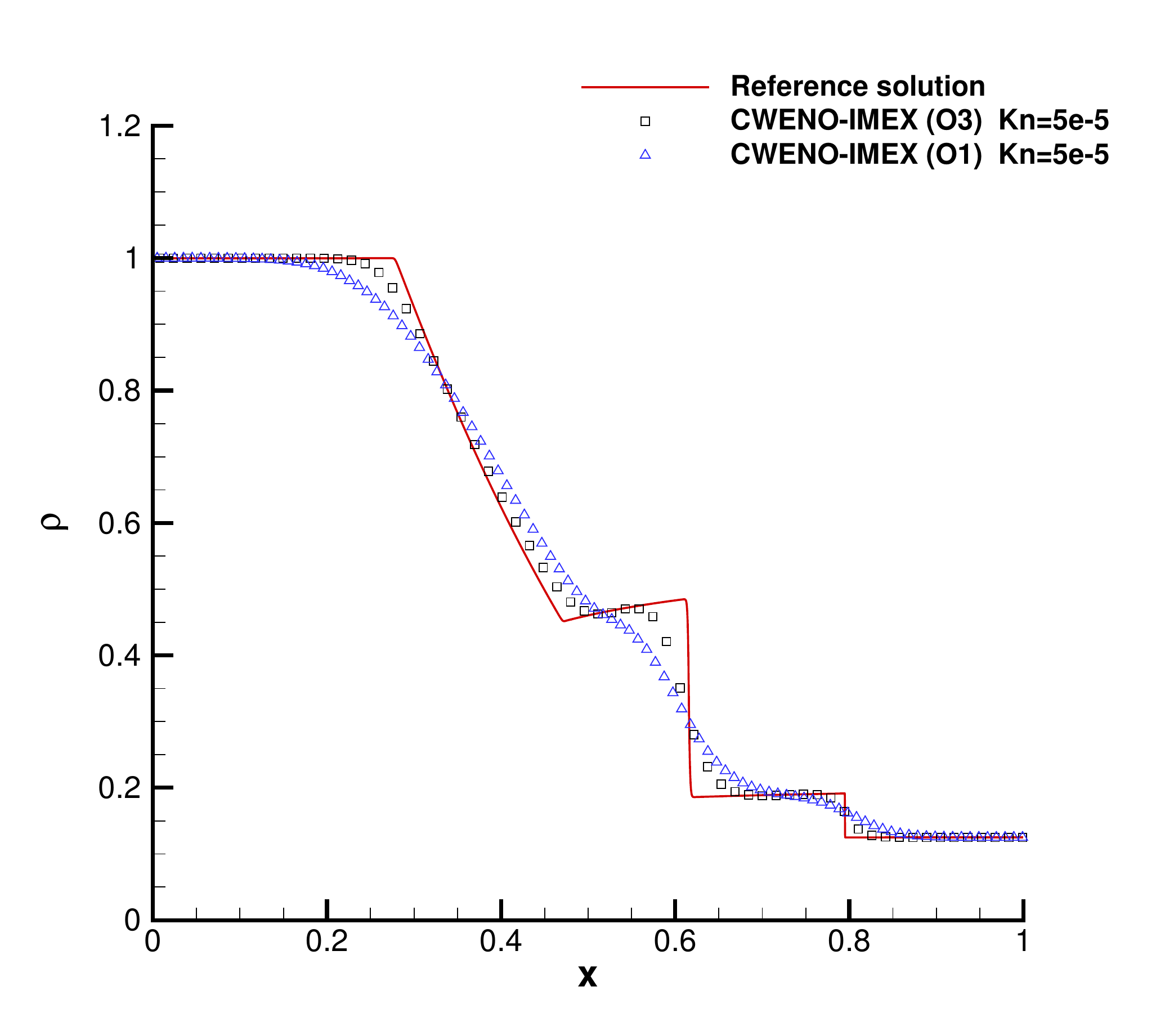} \\
			\includegraphics[width=0.47\textwidth]{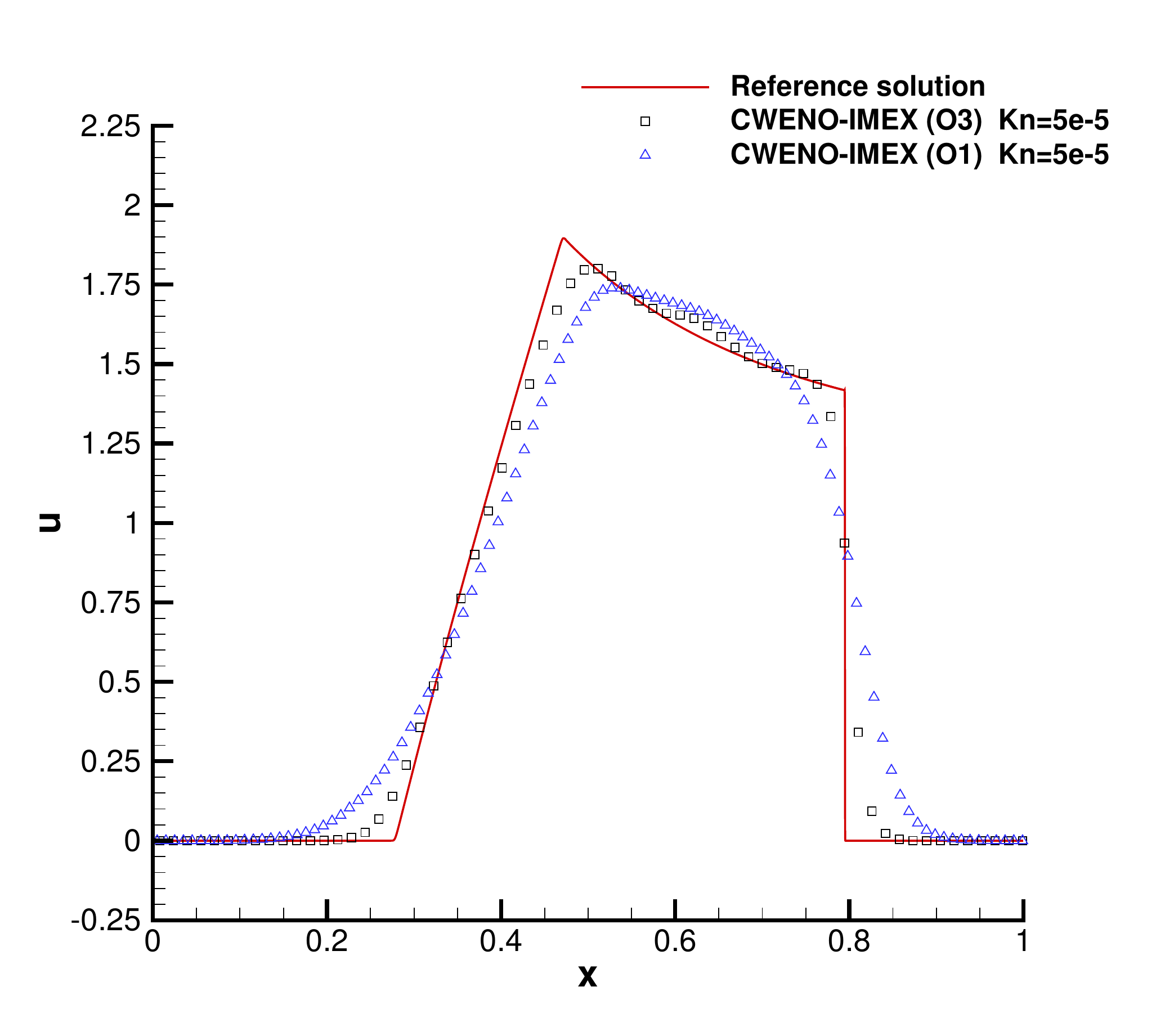}  &           
			\includegraphics[width=0.47\textwidth]{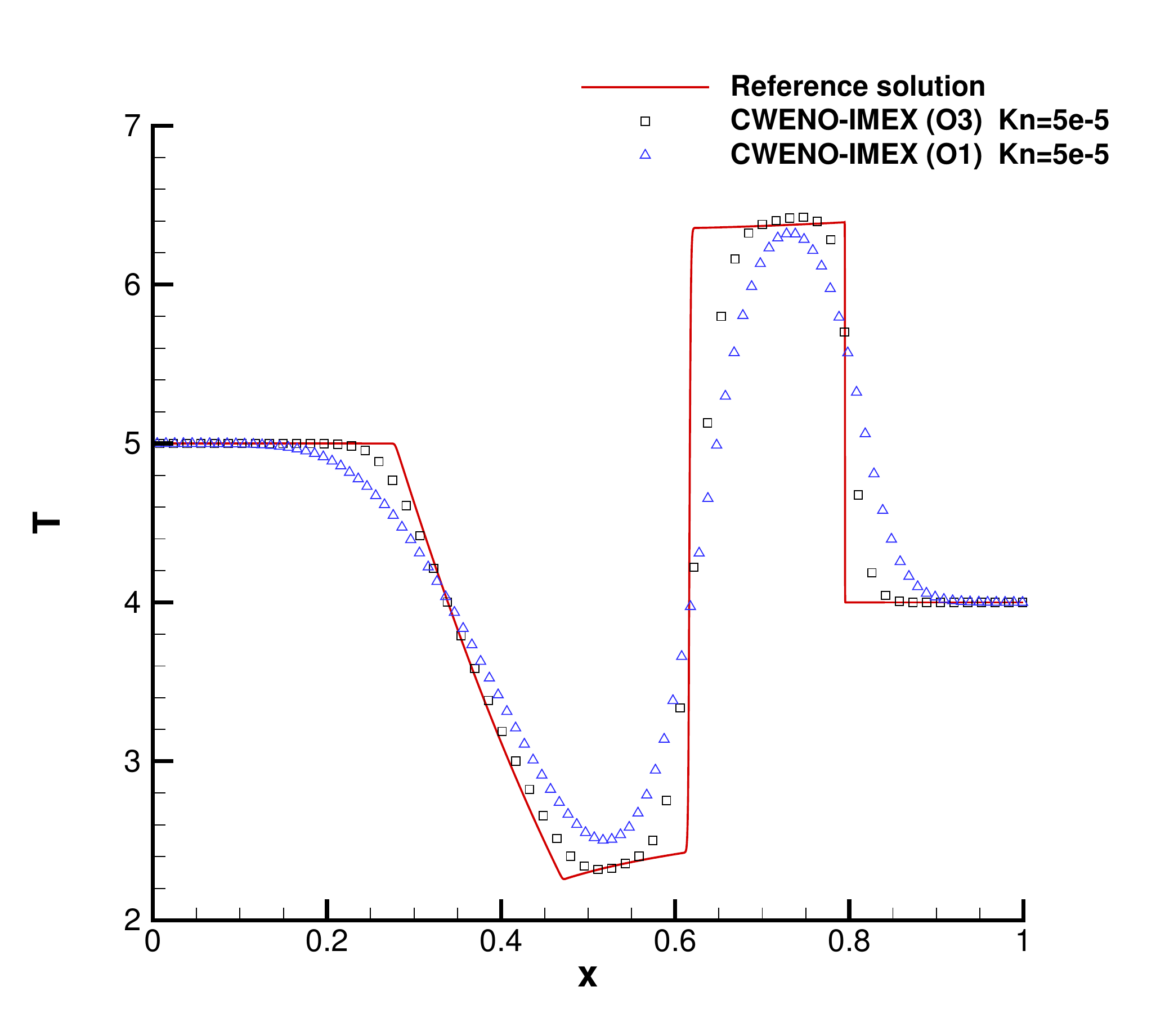} \\
		\end{tabular}
		\caption{Explosion problem with $\varepsilon=5 \cdot 10^{-5}$ at time $t_f=0.07$.  Three dimensional view of density profile for Boltzmann model together with a 1D cut along the $x$-axis through the third order numerical results and comparison with exact solution of the Euler equations for density, horizontal velocity and temperature.}
		\label{fig.EP-5e-5}
	\end{center}
\end{figure}

\begin{figure}[!htbp]
	\begin{center}
		\begin{tabular}{ccc} 
			\includegraphics[width=0.33\textwidth]{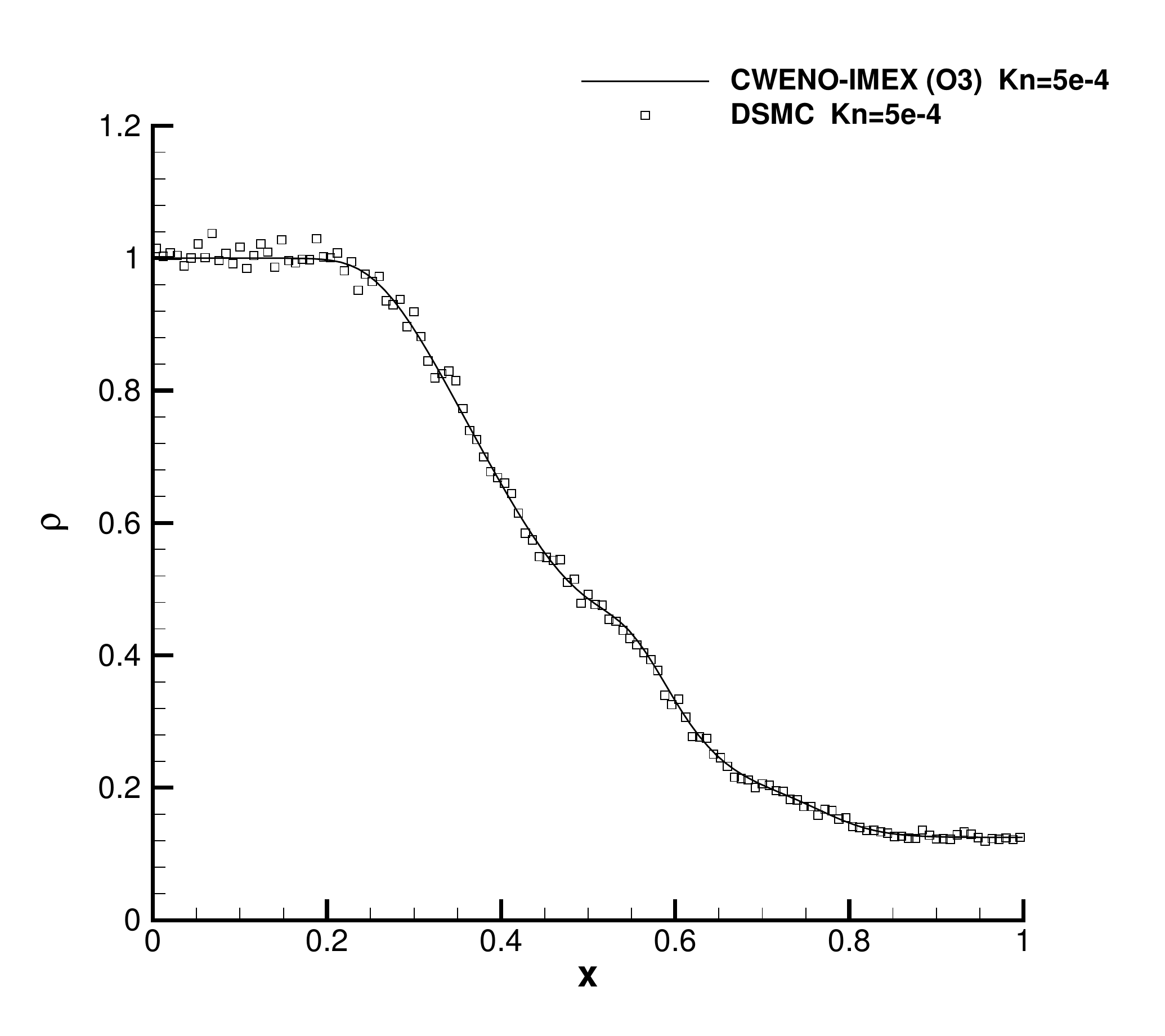}  &          
			\includegraphics[width=0.33\textwidth]{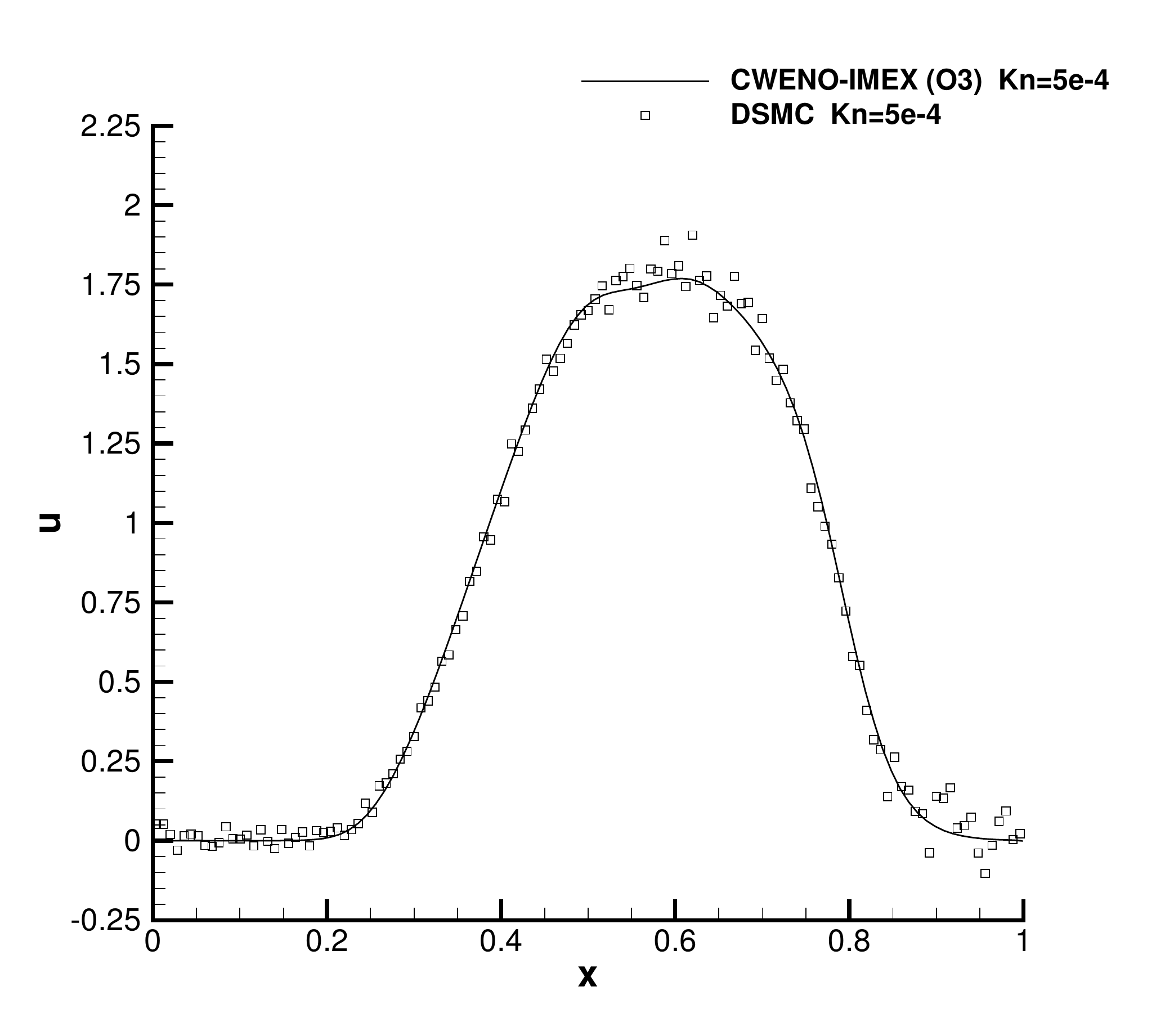} &
			\includegraphics[width=0.33\textwidth]{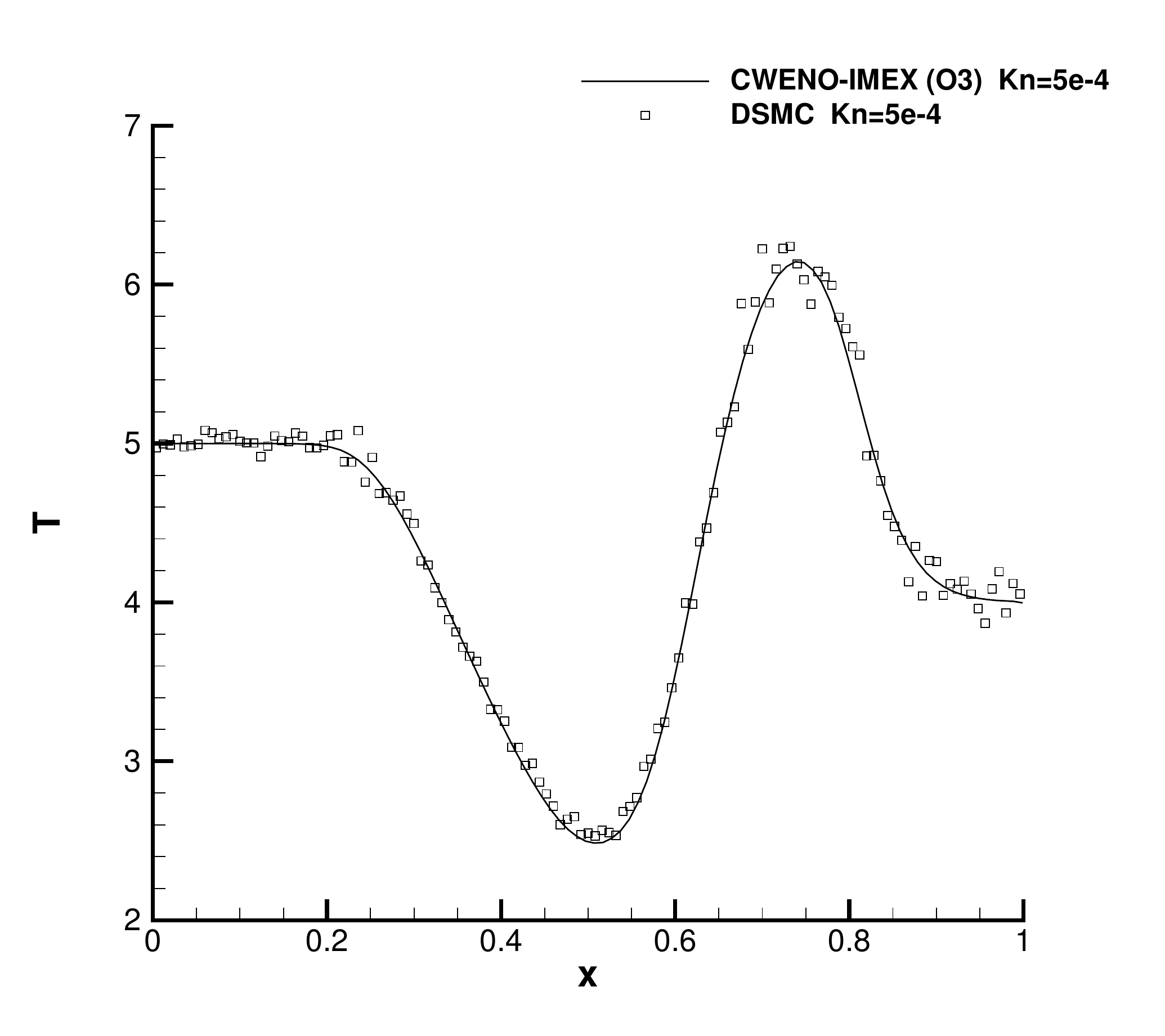}  \\           
			\includegraphics[width=0.33\textwidth]{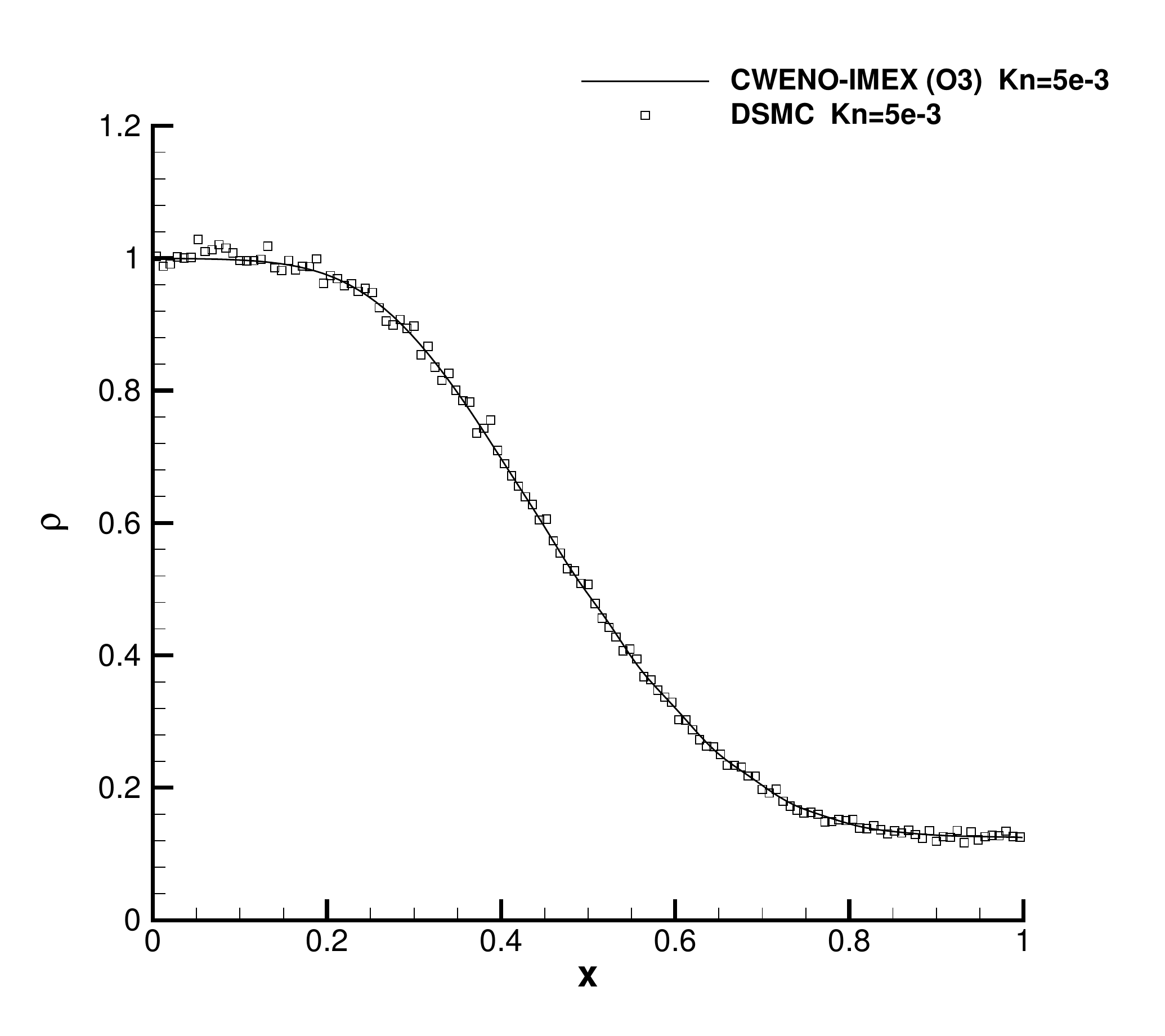}  &          
			\includegraphics[width=0.33\textwidth]{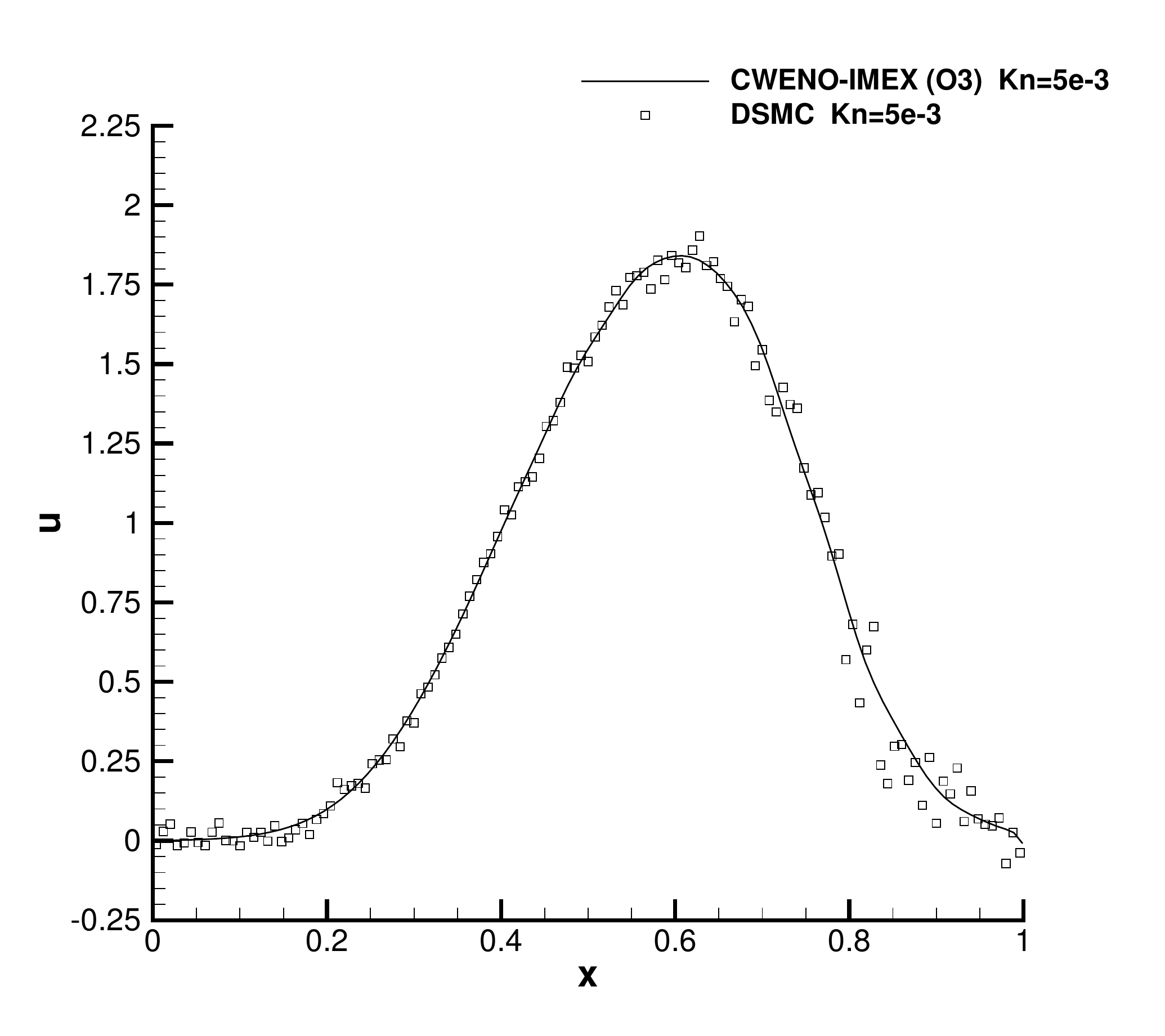} &
			\includegraphics[width=0.33\textwidth]{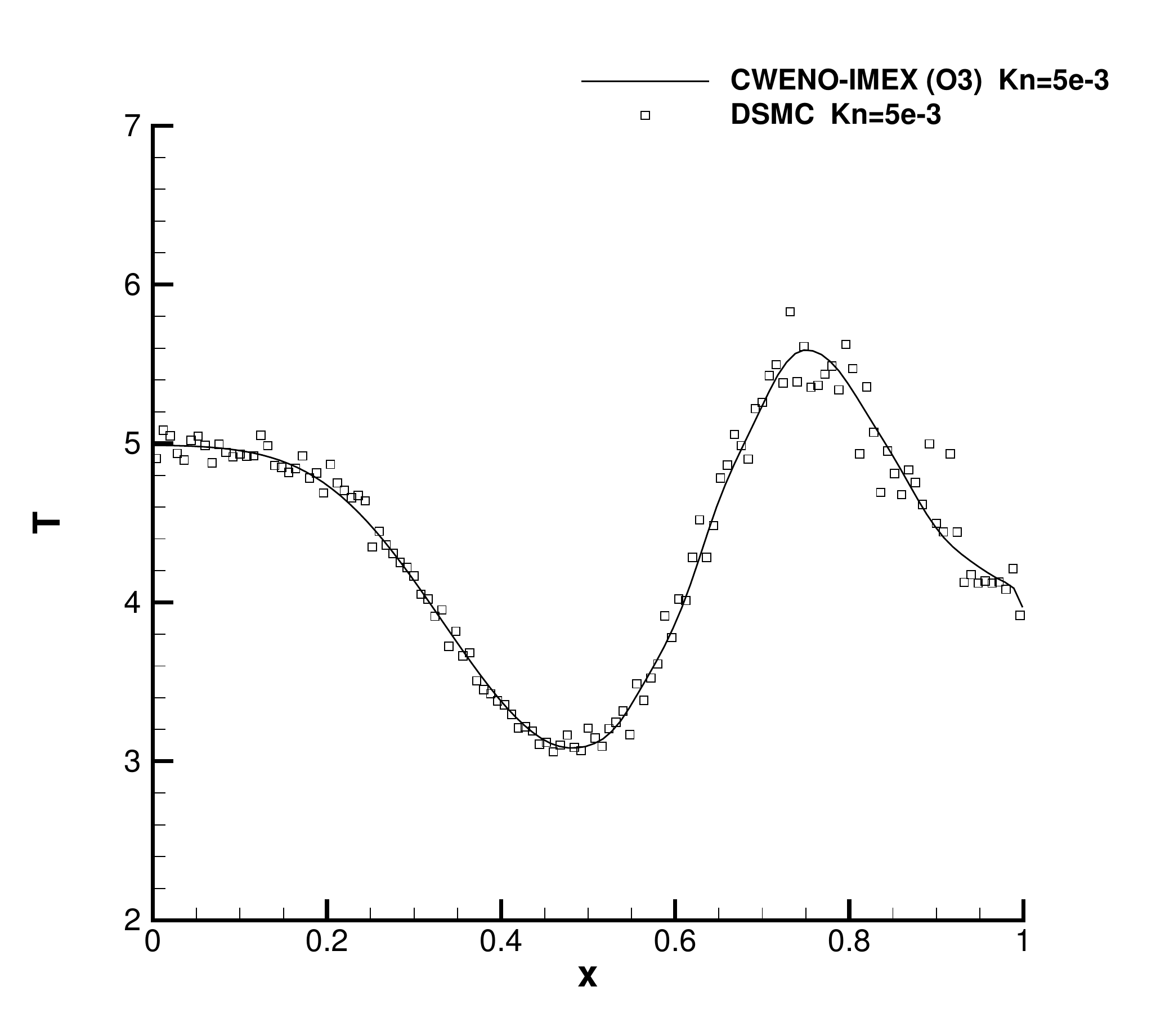}  \\
		\end{tabular}
	    \begin{tabular}{cc} 
	    	\includegraphics[width=0.47\textwidth]{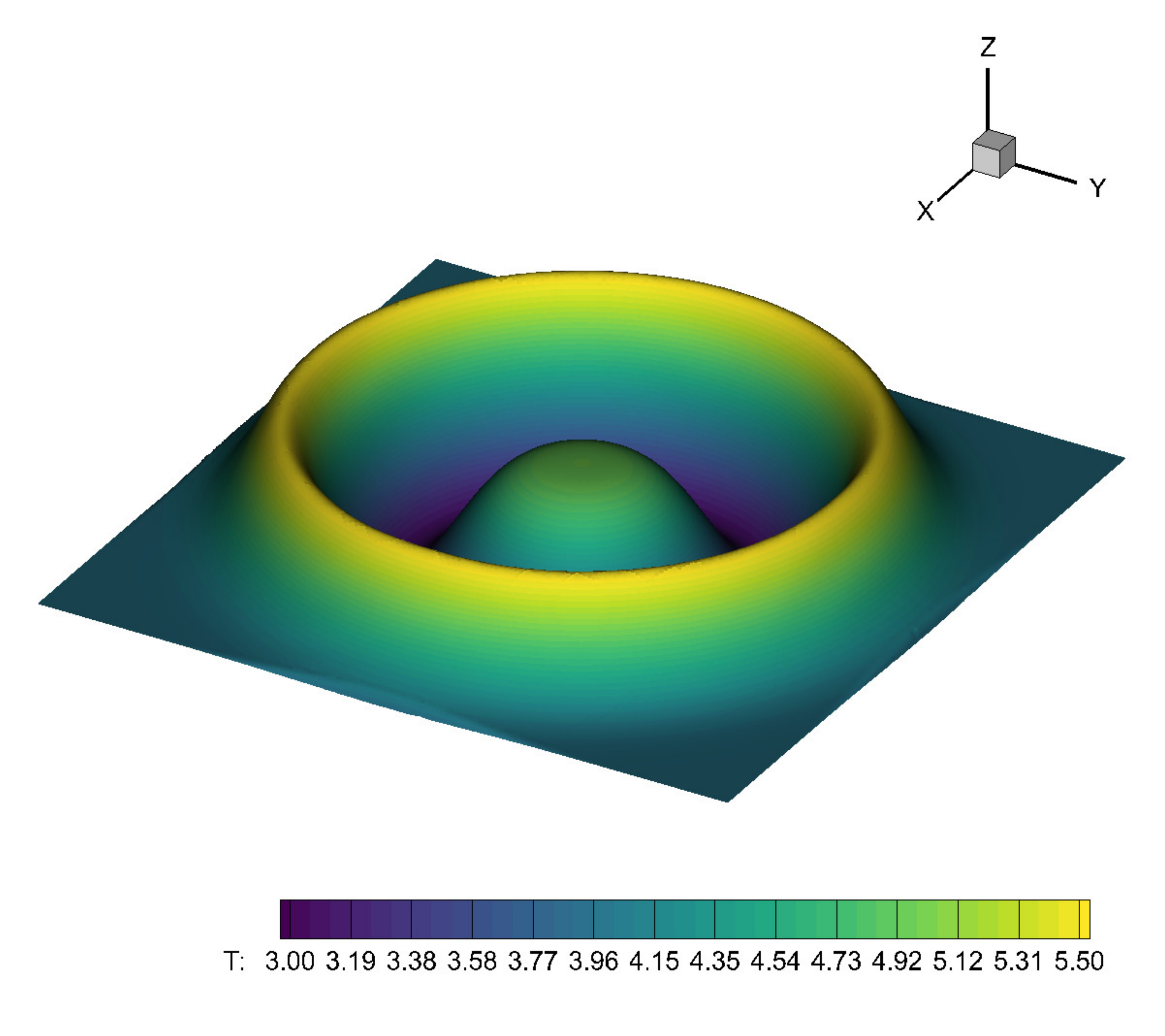}  &          
	    	\includegraphics[width=0.47\textwidth]{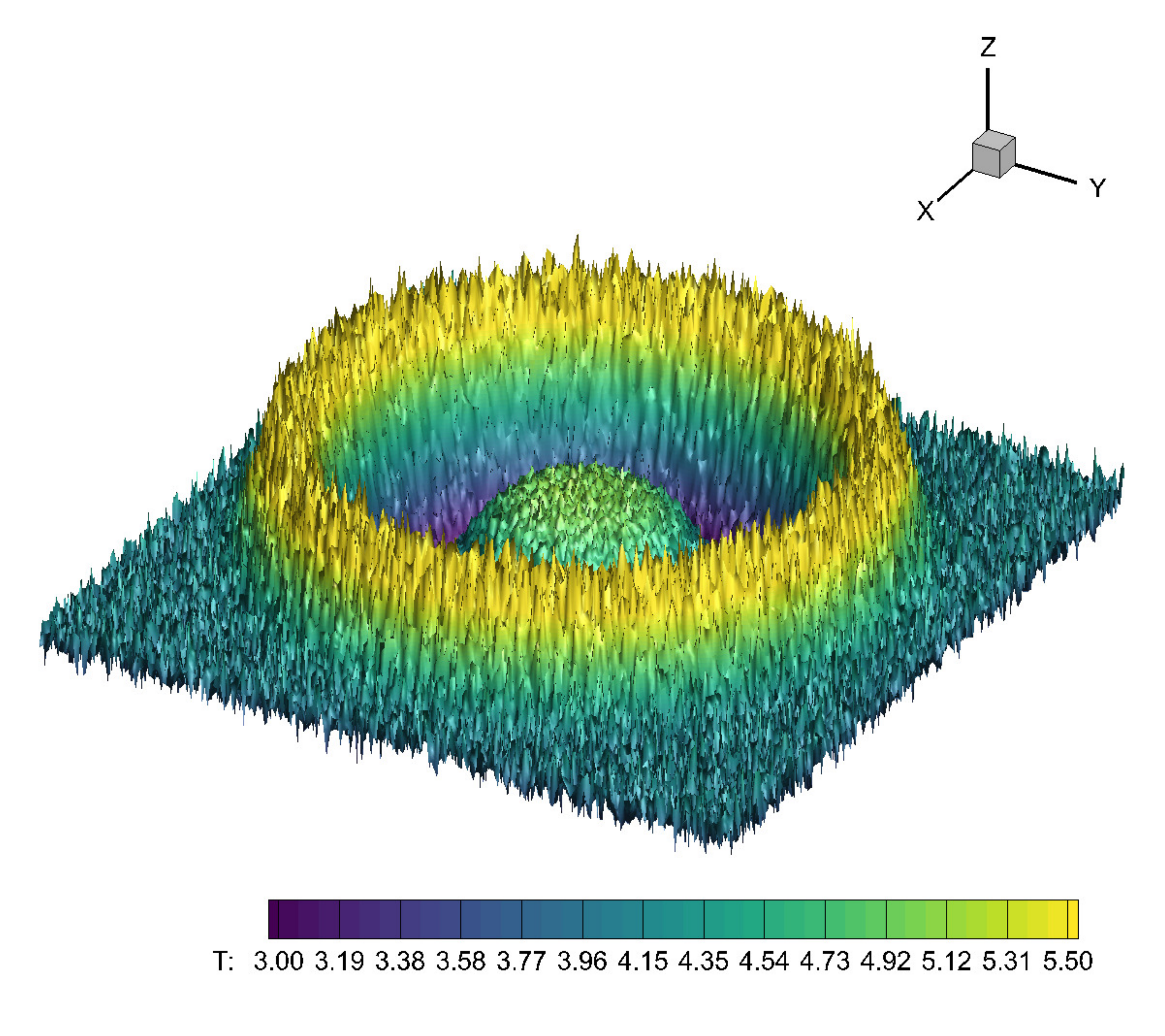}  \\
	    \end{tabular}	
		\caption{Explosion problem with $\varepsilon=5 \cdot 10^{-4}$ (top row) and $\varepsilon=5 \cdot 10^{-3}$ (middle row) at time $t_f=0.07$. 1D cut along the $x$-axis for density (left), horizontal velocity (middle) and temperature (right). The solid line represents the solution obtained with the third order CWENO-IMEX schemes, while scatter plots represent the results computed with a DSMC simulation on a grid composed of $250^2$ Cartesian cells with approximately $62 \cdot 10^6$ particles. Bottom: 3D view of the temperature distribution with 40 contour levels within the interval $[3;5.5]$ for the CWENO-IMEX (left) and DSMC (right) solver.}
		\label{fig.EP-Kn}
	\end{center}
\end{figure}

The setup of the explosion problem involves a square computational domain of dimension $\Omega=[-1;1]\times[-1;1]$, which is discretized by $N_P \approx 16'500$ Voronoi elements. The initial condition is given in terms of two different states $(\mathcal{U}_L,\mathcal{U}_R)$, separated by a discontinuity at distance $R_d=0.5$ from the origin: 
\begin{equation}
	\begin{cases} 
		\mathcal{U}_L = \left(1,0,0,1 \right) & r \leq R_d, \\
		\mathcal{U}_R = \left(0.125,0,0,0.8\right) & r > R_d,
	\end{cases}
\end{equation}
with $r=\sqrt{x^2+y^2}$ the radial position. Dirichlet boundary conditions are imposed everywhere and the final time of the simulation is chosen to be $t_f=0.07$. The velocity space is composed of $32^2=1024$ Cartesian cells ranging in the interval $[-20;20]^2$. Figure \ref{fig.EP-5e-5} depicts the third order numerical results obtained with the CWENO-IMEX scheme with Knudsen number $\varepsilon=5 \cdot 10^{-5}$ which approaches the fluid limit. A comparison against the reference solution of the Euler equations of compressible gas dynamics \cite{Lagrange2D,ToroBook} is proposed, showing an overall good agreement considering the very low number of mesh elements employed. The same simulation is also run with a first order finite volume scheme which exhibits much more numerical dissipation, as expected. Let us notice that such smoother profile could even be wrongly interpreted as a different physical regime of the gas flow, mimicking the behavior of the fluid at larger Knudsen numbers. This underlines the importance of the development of high order numerical discretizations, in order to correctly capture and simulate the physical behavior of the observed phenomena. Finally, Figure \ref{fig.EP-Kn} reports the numerical distribution of density, horizontal velocity and temperature with Knudsen numbers $\varepsilon=5 \cdot 10^{-4}$ and $\varepsilon=5 \cdot 10^{-3}$. For these regimes, a comparison with a DSMC method is proposed. This method works on a regular Cartesian mesh of $250^2$ elements and employs around 62 millions of particles. As for the case of the fluid limit a good agreement is observed between the two solvers.

%
\subsection{Double Mach reflection problem}
\begin{figure}[!htbp]
	\begin{center}
		\begin{tabular}{cc} 
			\includegraphics[width=0.47\textwidth]{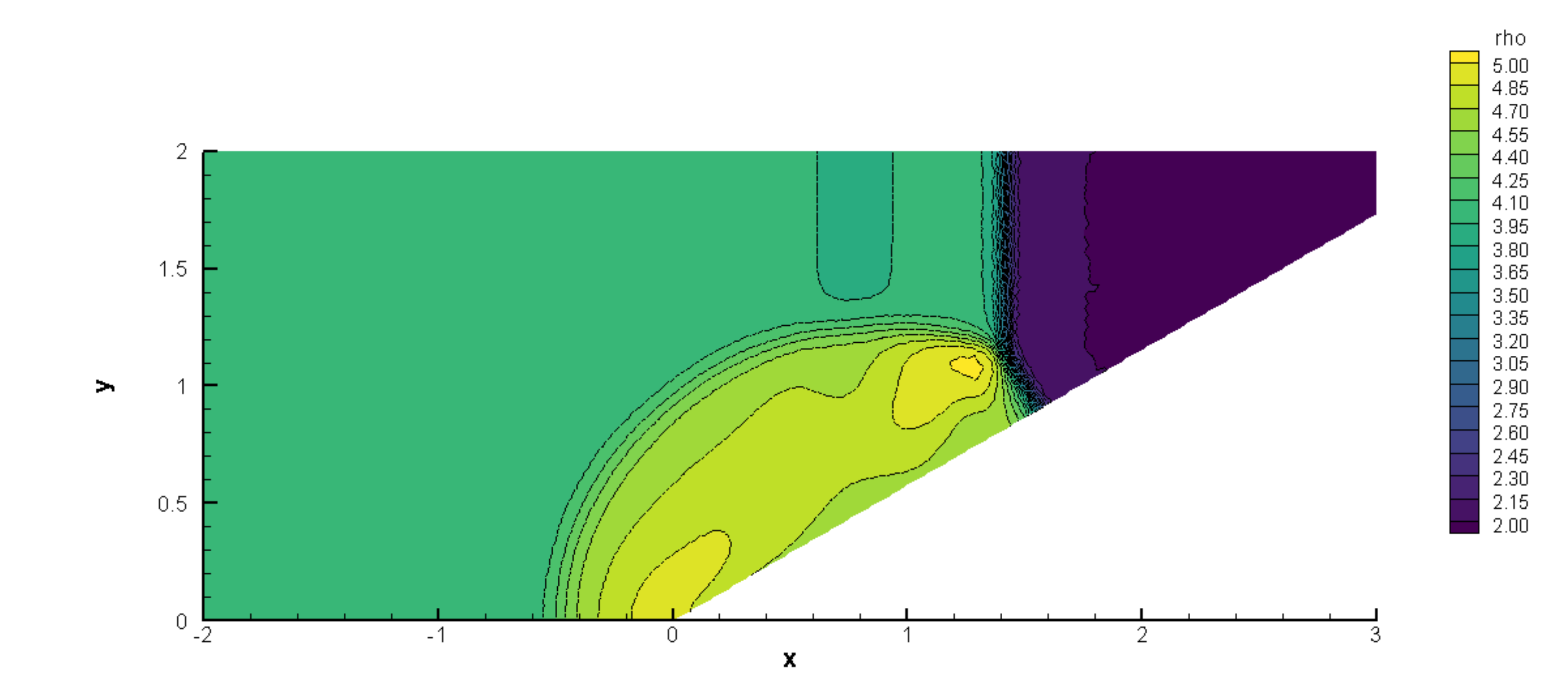}  &          
			\includegraphics[width=0.47\textwidth]{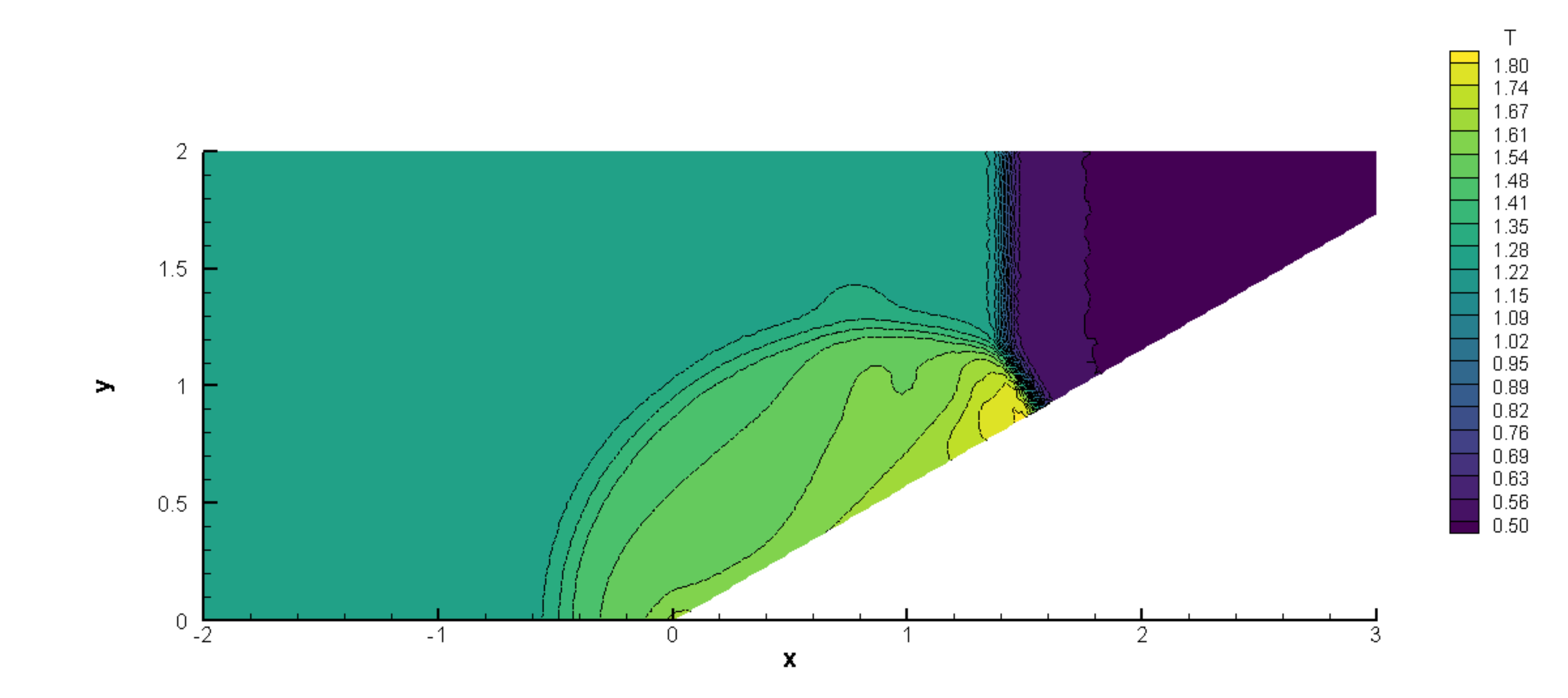}  \\
			\includegraphics[width=0.47\textwidth]{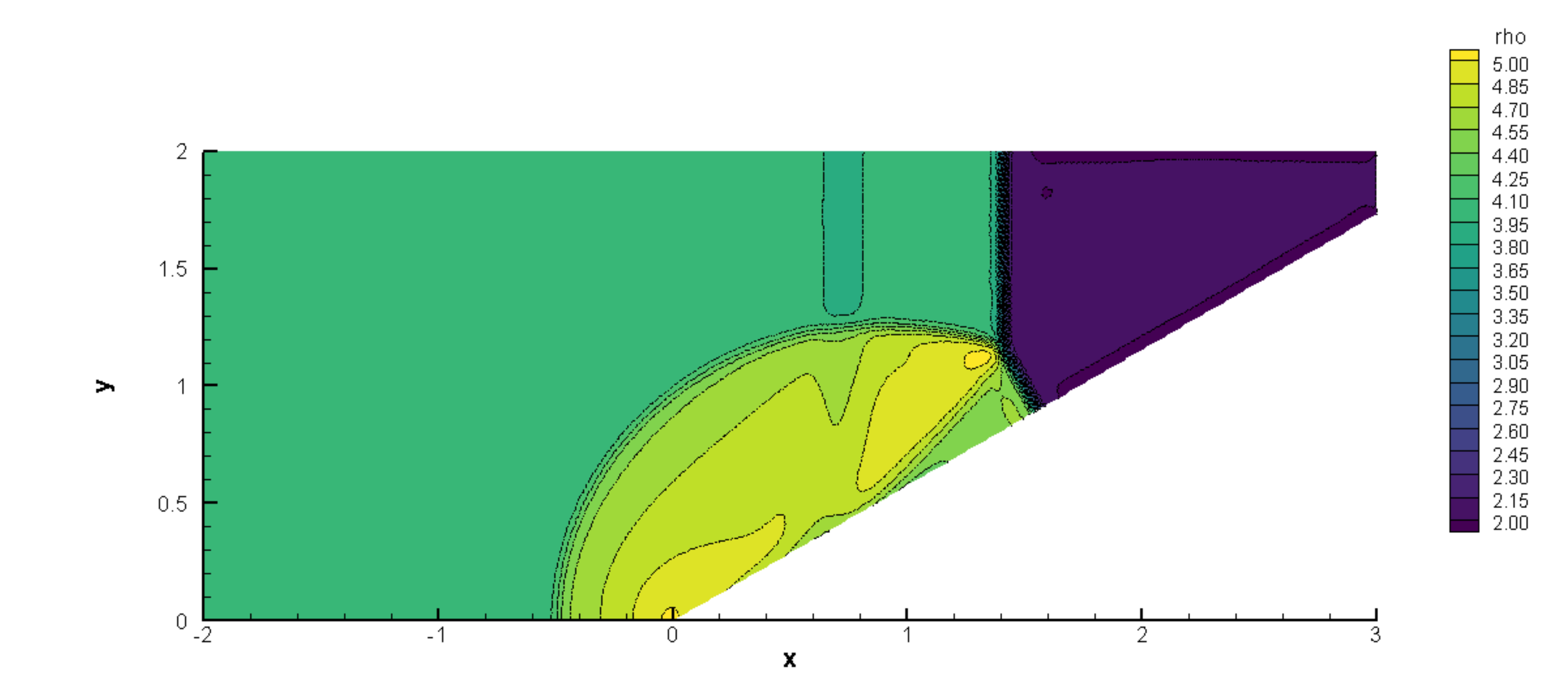}  &          
			\includegraphics[width=0.47\textwidth]{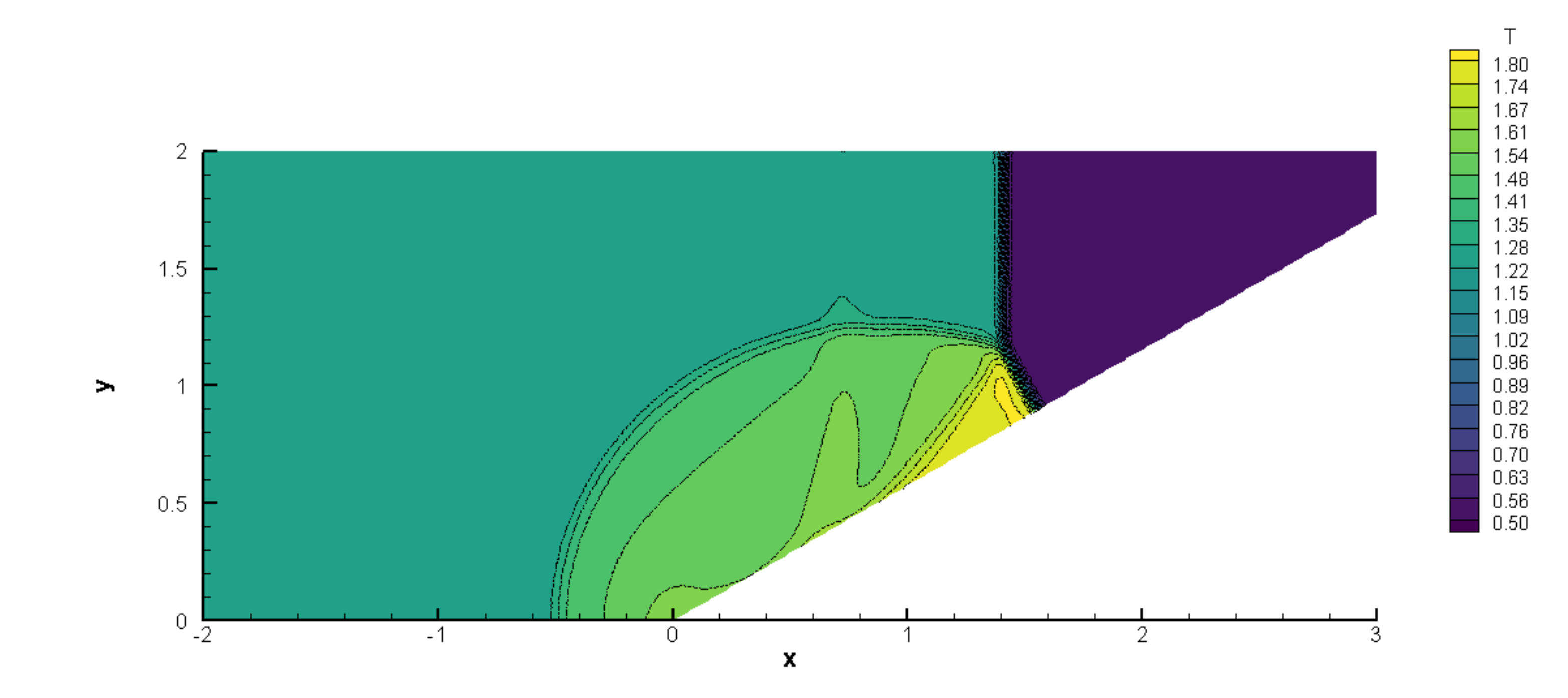}  \\           
		\end{tabular}
		\caption{Double Mach reflection problem at final time $t_f=0.7$. Third order numerical solution for density (left) and temperature (right) for the Euler equations (top row) and the Boltzmann model with $\varepsilon=5\cdot 10^{-5}$ (bottom row).}
		\label{fig.DMR-Eul-Boltz-comp}
	\end{center}
\end{figure}

\begin{figure}[!htbp]
	\begin{center}
		\begin{tabular}{cc} 
			\includegraphics[width=0.47\textwidth]{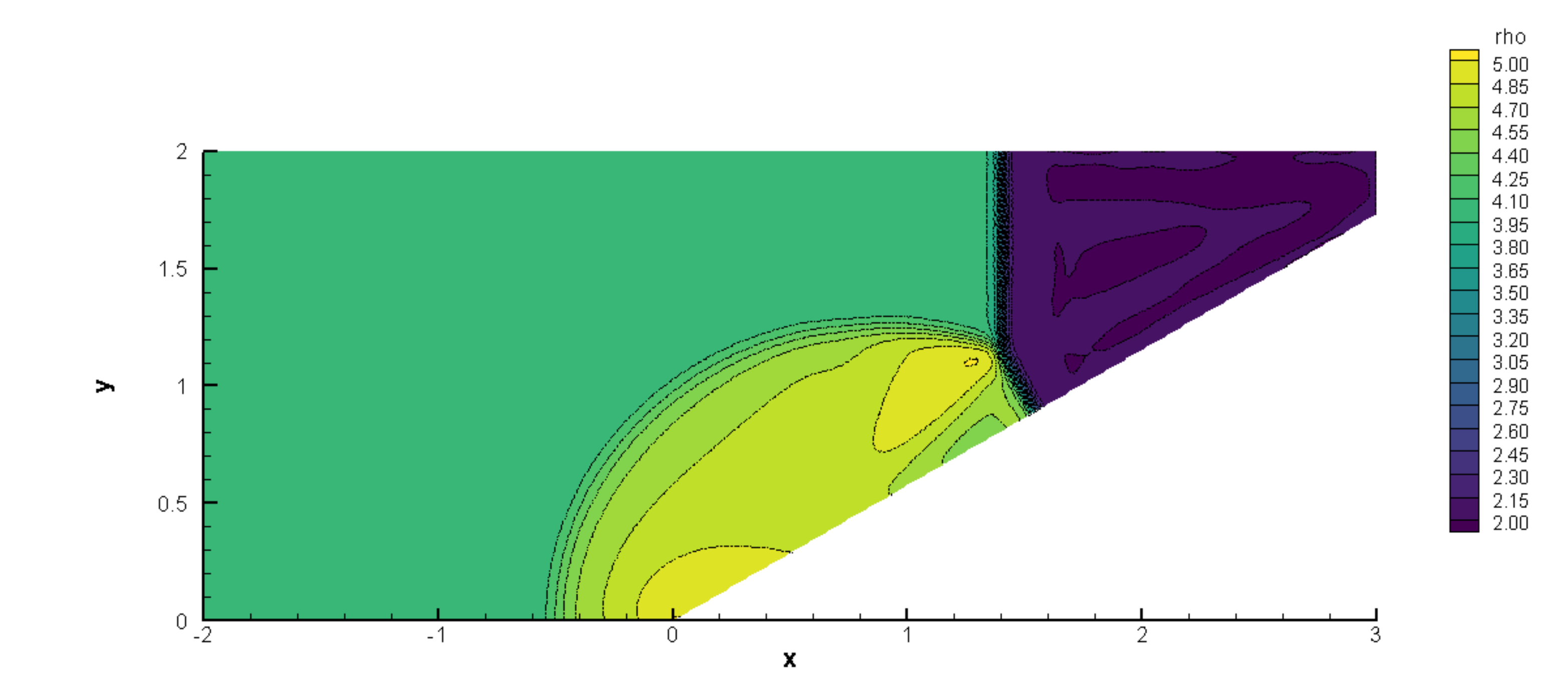}  &          
			\includegraphics[width=0.47\textwidth]{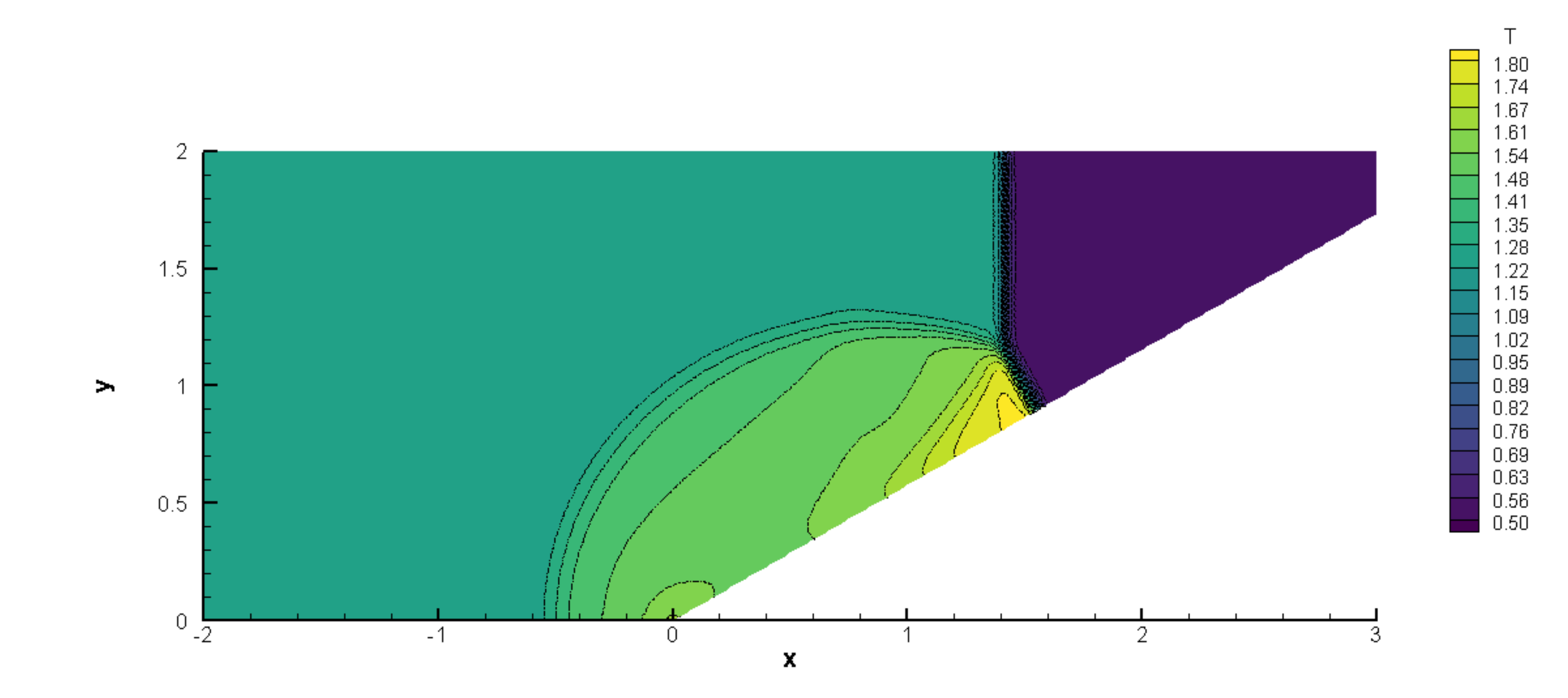}  \\
			\includegraphics[width=0.47\textwidth]{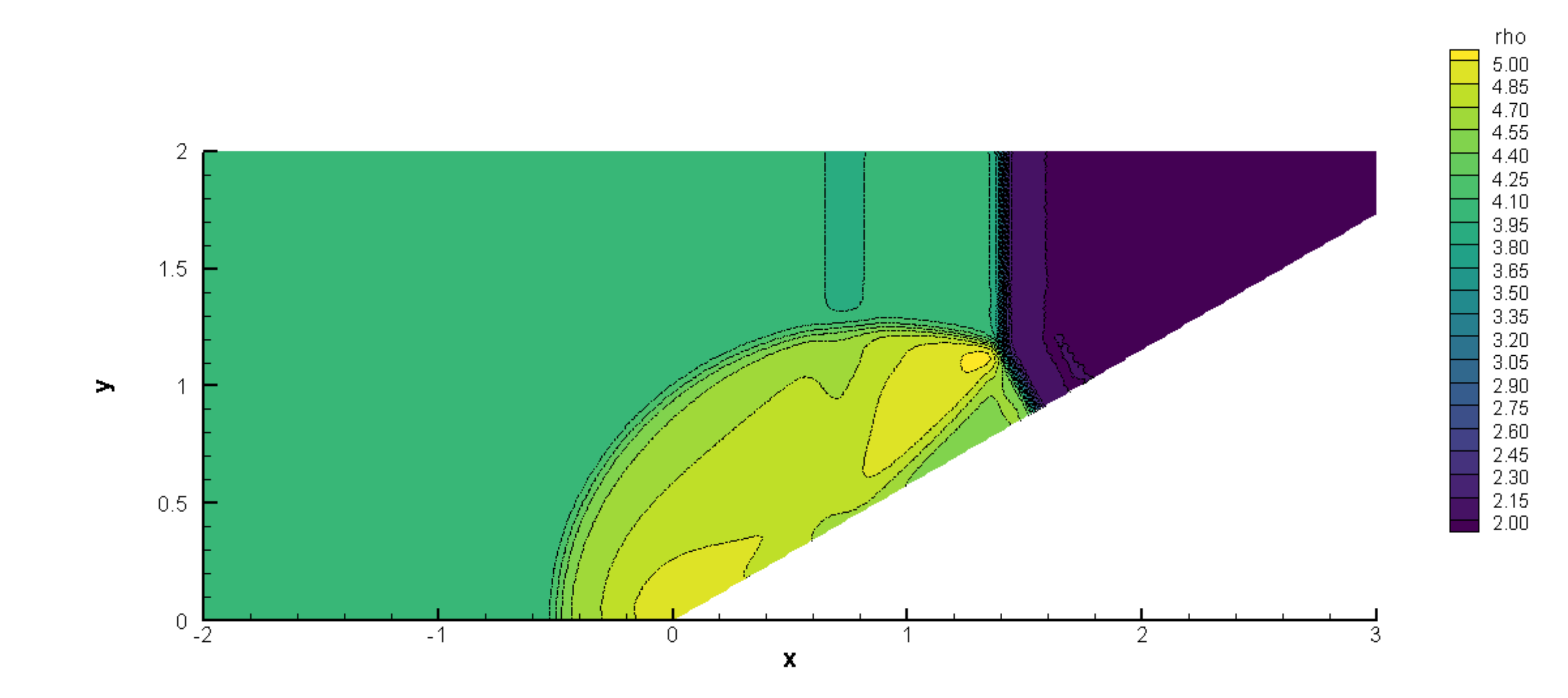}  &          
			\includegraphics[width=0.47\textwidth]{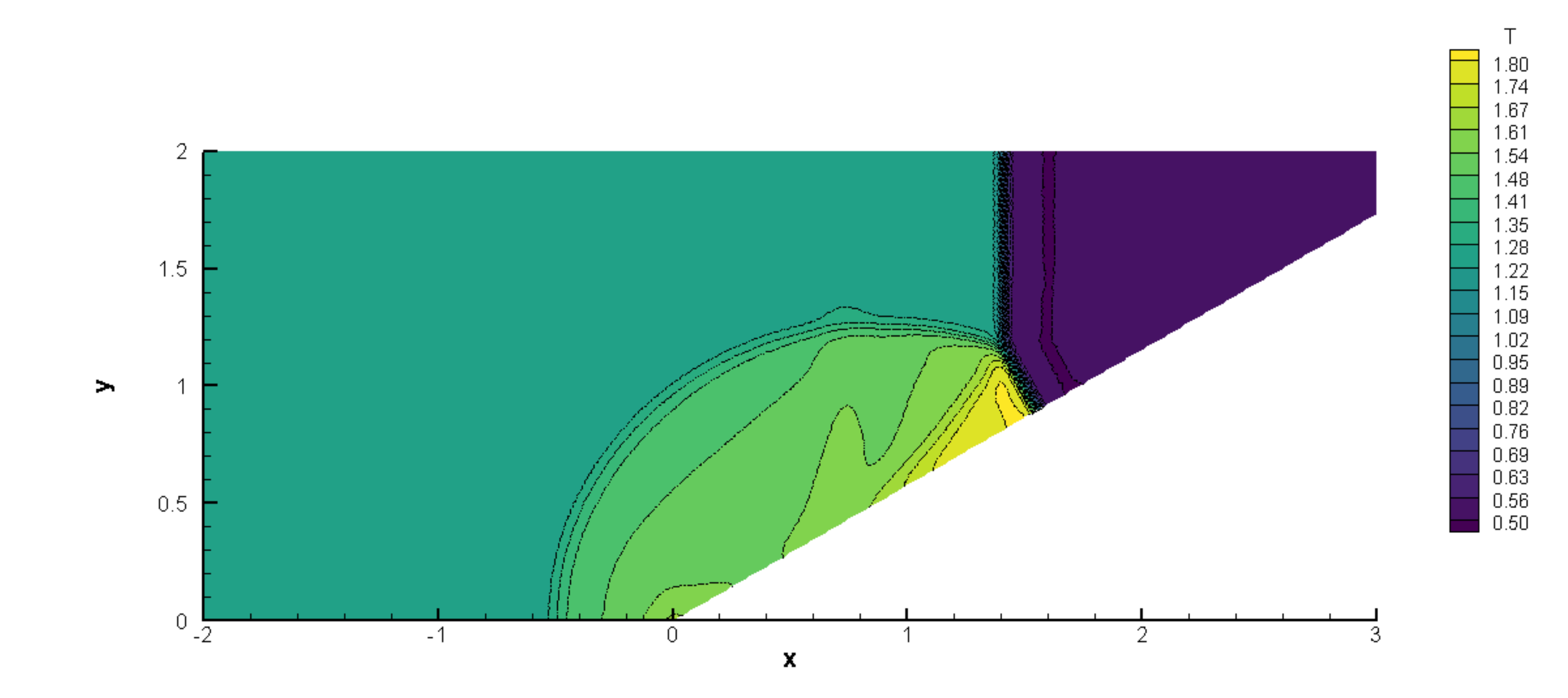}  \\           
		\end{tabular}
		\caption{Double Mach reflection problem at final time $t_f=0.7$ with $\varepsilon=5\cdot 10^{-3}$. Third order numerical solution for density (left) and temperature (right) for the Boltzmann (top row) and the BGK (bottom row) model.}
		\label{fig.DMR-BGK-Boltz-comp}
	\end{center}
\end{figure}

In this test, we propose to run a slightly modified setting of the double Mach reflection problem originally proposed in \cite{woodwardcol84}. The setup involves a shock wave with a shock Mach number $M_s=2$ which is moving along the $x-$direction of the computational domain, where a ramp with angle $\alpha=\frac{\pi}{6}$ is located. The strong shock wave that hits the ramp yields the development of other two shock waves, one travelling towards the right and the other one propagating towards the top boundary of the domain. The initial computational domain $\Omega$ is discretized with a total number of elements of $N_P=6714$.
The initial condition is given in terms of primitive variables and explicitly writes 
\begin{equation}
	\begin{cases} 
		\mathcal{U}_L = \left(4,1,0,1.25 \right) & x \leq x_d, \\
		\mathcal{U}_R = \left(2,0,0,0.50\right) & x > x_d,
	\end{cases}
\end{equation}
with the initial discontinuity located at $x_d=0$. Dirichlet boundary conditions are set on the left and rights sides, while sliding wall boundary conditions have been imposed on the remaining sides. The velocity space is bounded in the interval $[-10;10]^2$ and is discretized with $32^2=1024$ Cartesian cells. The final time is $t_{f}=0.7$. The Figure \ref{fig.DMR-Eul-Boltz-comp} shows the third order numerical solution for density and temperature with $\varepsilon=5\cdot 10^{-5}$ together with a comparison against the solution obtained by solving the compressible Euler equations with the same order of accuracy. A qualitatively agreement of the solution can be appreciated, with the physical variable contours approximated by 21 equidistant levels for density and temperature in the interval $[2;5]$ and $[0.5;1.8]$, respectively. A comparison with the BGK solution at Knudsen number $\varepsilon=5\cdot 10^{-3}$ is finally proposed in Figure \ref{fig.DMR-BGK-Boltz-comp}, highlighting some differences in the behavior of the two models. Specifically, the Boltzmann model exhibits a smoother solution for both density and temperature, with less small-scale structures around the triple point, which corresponds to the intersection of the three shock waves. Compared to the standard double Mach reflection test \cite{woodwardcol84}, here the incidence wave has quite a slower speed which causes a simpler structure of the final solution. This choice is imposed by the impossibility to correctly describe the velocity space involved in such a problem using a fixed and uniform in space discretization of the phase-space. In fact, in order to perform such a study, it would be necessary to describe with the same velocity mesh, distribution functions with low mean velocities and temperatures and distribution functions exhibiting very large mean velocities and temperatures. This problem needs an approach with velocity grids able to adapt to the time evolution of the solution both in time and space. This kind of approach is the scope of future investigations.  

%
\subsection{Fluid flow around NACA 0012 airfoil}

\begin{figure}[!htbp]
	\begin{center}
		\begin{tabular}{ccc} 
			\includegraphics[width=0.33\textwidth]{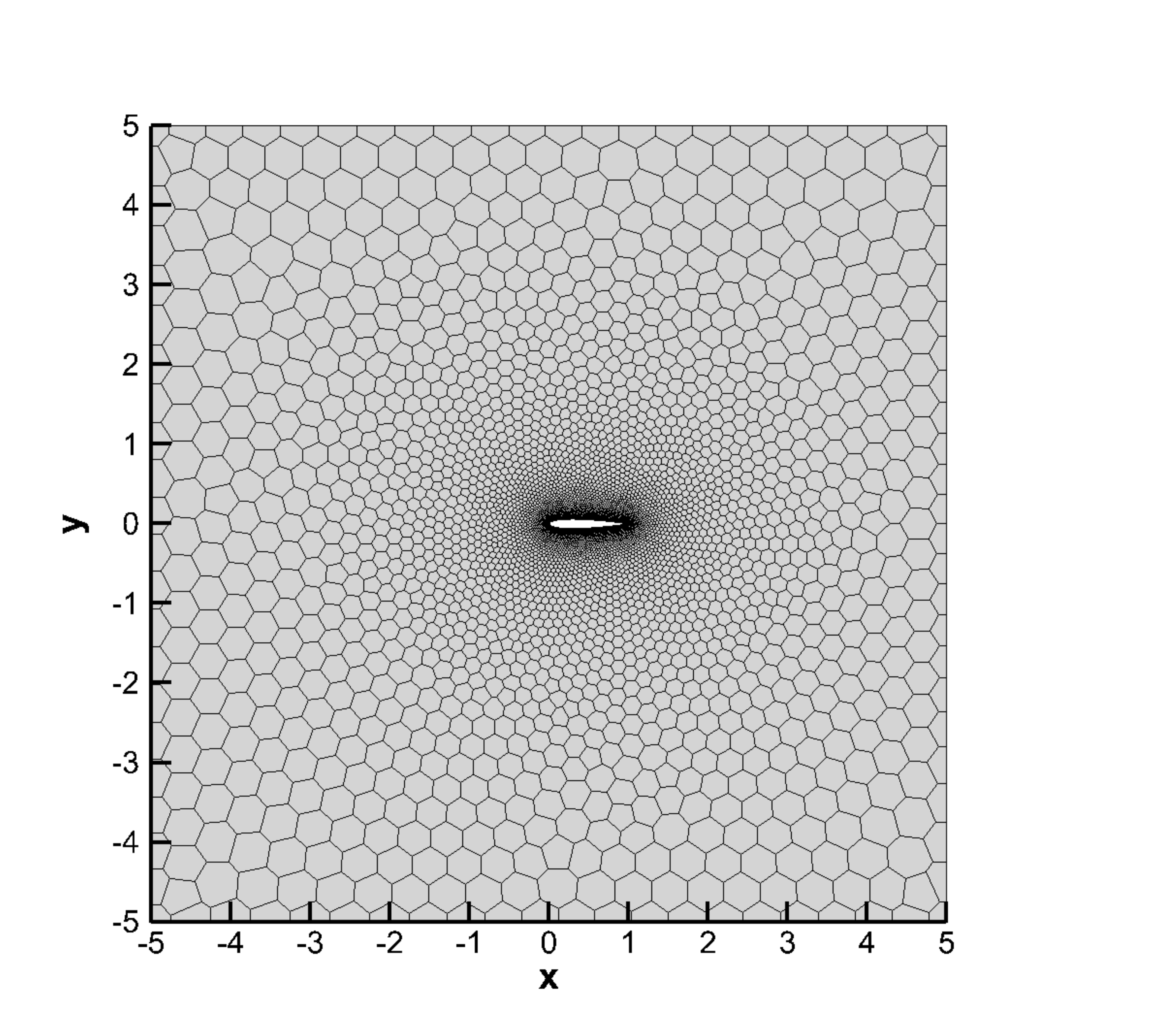}  &          
			\includegraphics[width=0.33\textwidth]{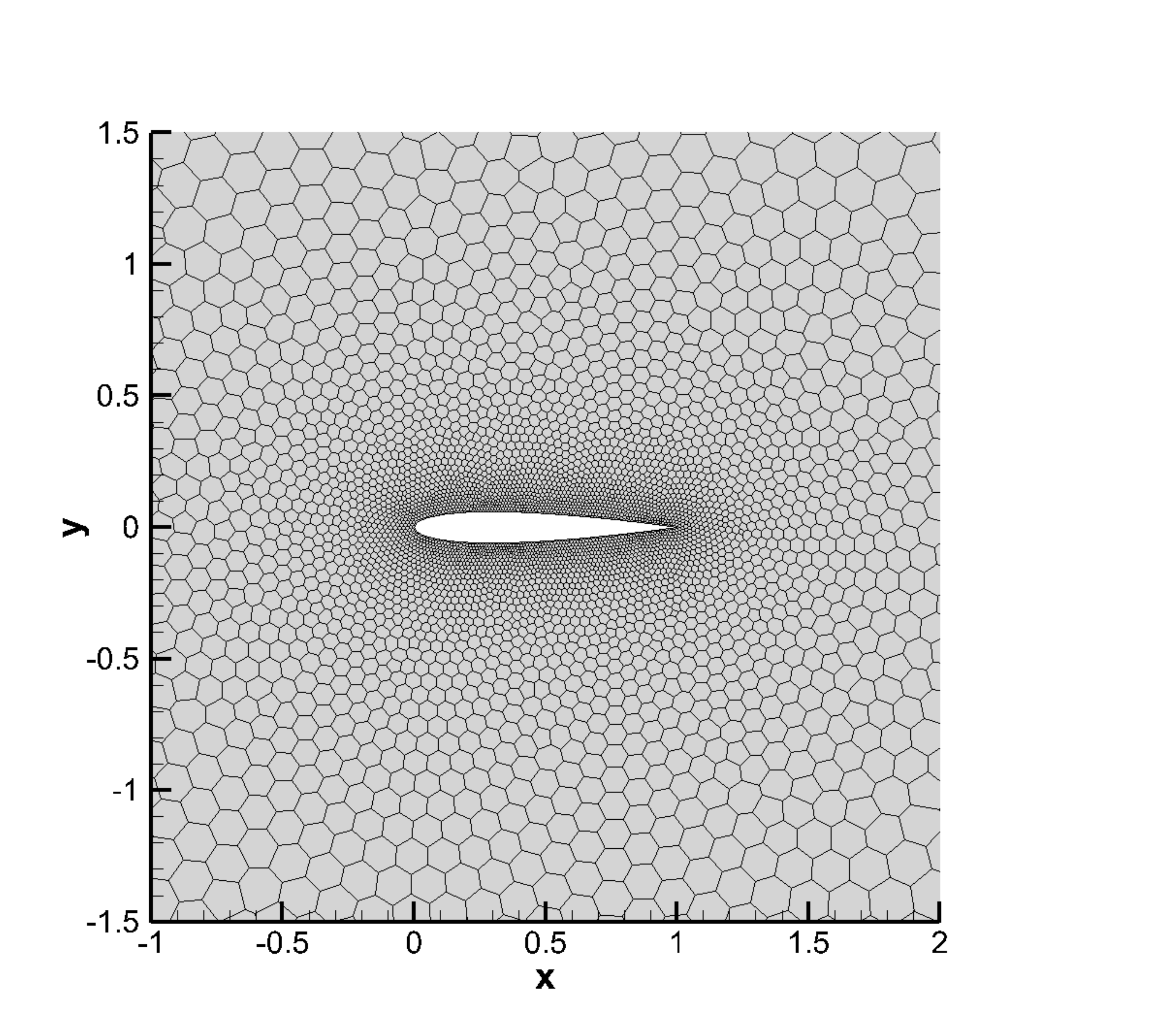}  &
			\includegraphics[width=0.33\textwidth]{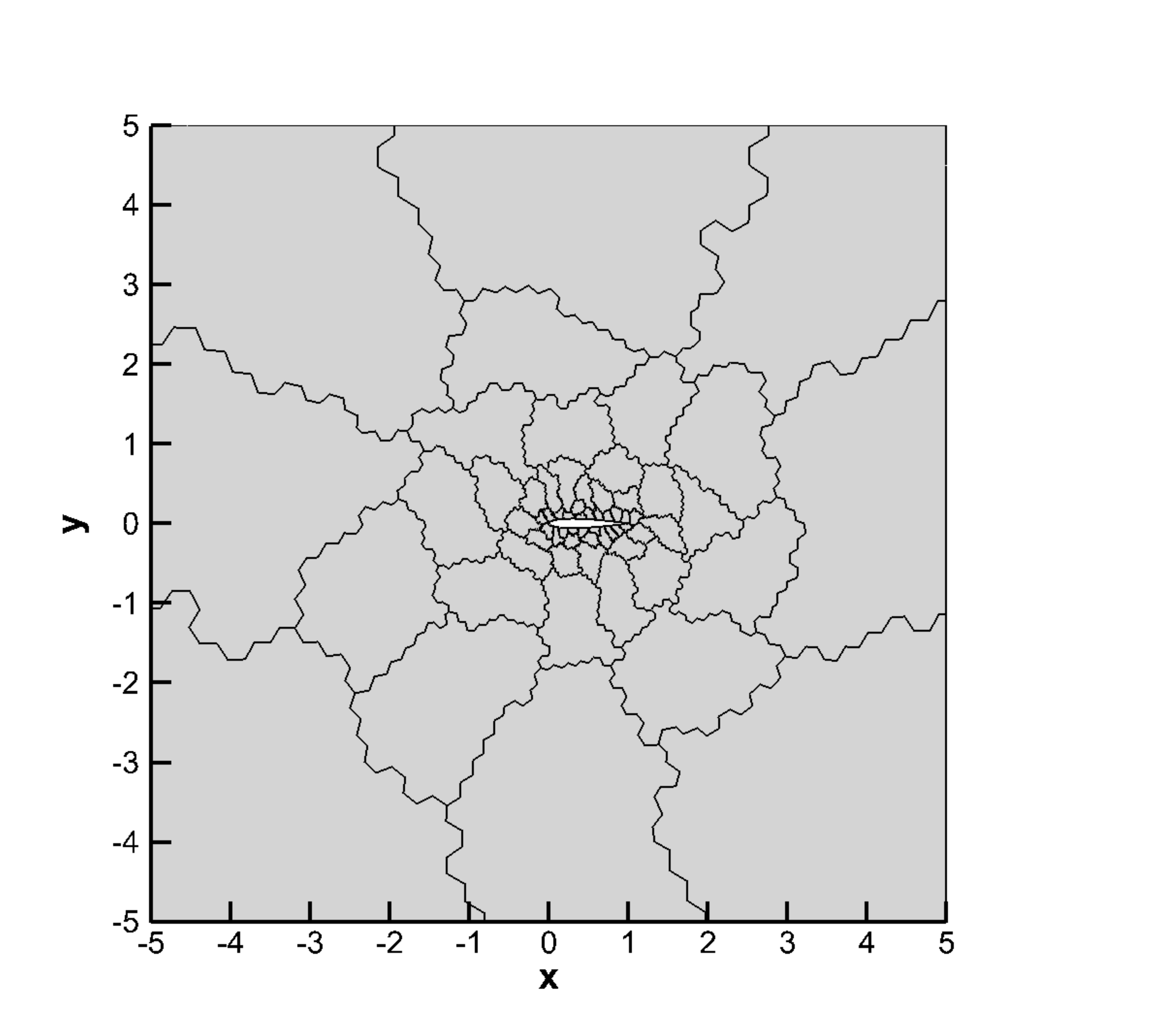}  \\           
		\end{tabular}
		\caption{Flow around NACA 0012 airfoil. Computational mesh with $N_P=10336$ polygonal cells (left), zoom around the airfoil profile discretized with 200 points (middle) and MPI domains which splits the computational domain into a total number of 64 sub-domains.}
		\label{fig.NACA-mesh}
	\end{center}
\end{figure}

The last test case concerns a more realistic application, namely the study of the flow around a NACA 0012 airfoil profile for different rarefied regimes. A slightly modified geometry from the original definition is considered, so that the airfoil closes at chord $c=1$ with a sharp trailing edge. This altered two-dimensional geometrical configuration of NACA 0012 airfoil is defined by
\begin{equation}
y = \pm 0.594689181 \, \left( 0.298222773\sqrt{x} - 0.127125232 x - 0.357907906 x^2 + 0.291984971 x^3 - 0.105174606 x^4 \right),
\end{equation}
which corresponds to a perfect scaled copy of the NACA 0012 profile with $x \in [0;1.008930411365]$. The airfoil is embedded in a square computational domain $\Omega=[-5;5]^2$ that is discretized with a total number of $N_P=10336$ polygonal control volumes. The profile of the airfoil to which is assigned a slip-wall boundary condition, is approximated with 100 equidistant points in the $x-$direction both on the upper and the lower border. The characteristic mesh size is then progressively increased linearly with the distance from the airfoil, thus ranging from $h=1/100$ to $h=1/10$. Inflow boundary conditions are imposed on the left and the bottom face of the domain, while transmissive conditions have been set elsewhere. The entire computational mesh and a zoom around the NACA profile is depicted in Figure \ref{fig.NACA-mesh}, together with the MPI partition of the domain within 64 CPUs. The velocity domain is discretized with $32^2$ Cartesian cells ranging in the interval $\mathcal{V}=[-10;10]\times [-10;10]$.

\begin{figure}[!htbp]
	\begin{center}
		\begin{tabular}{cc}
			\includegraphics[width=0.46\textwidth]{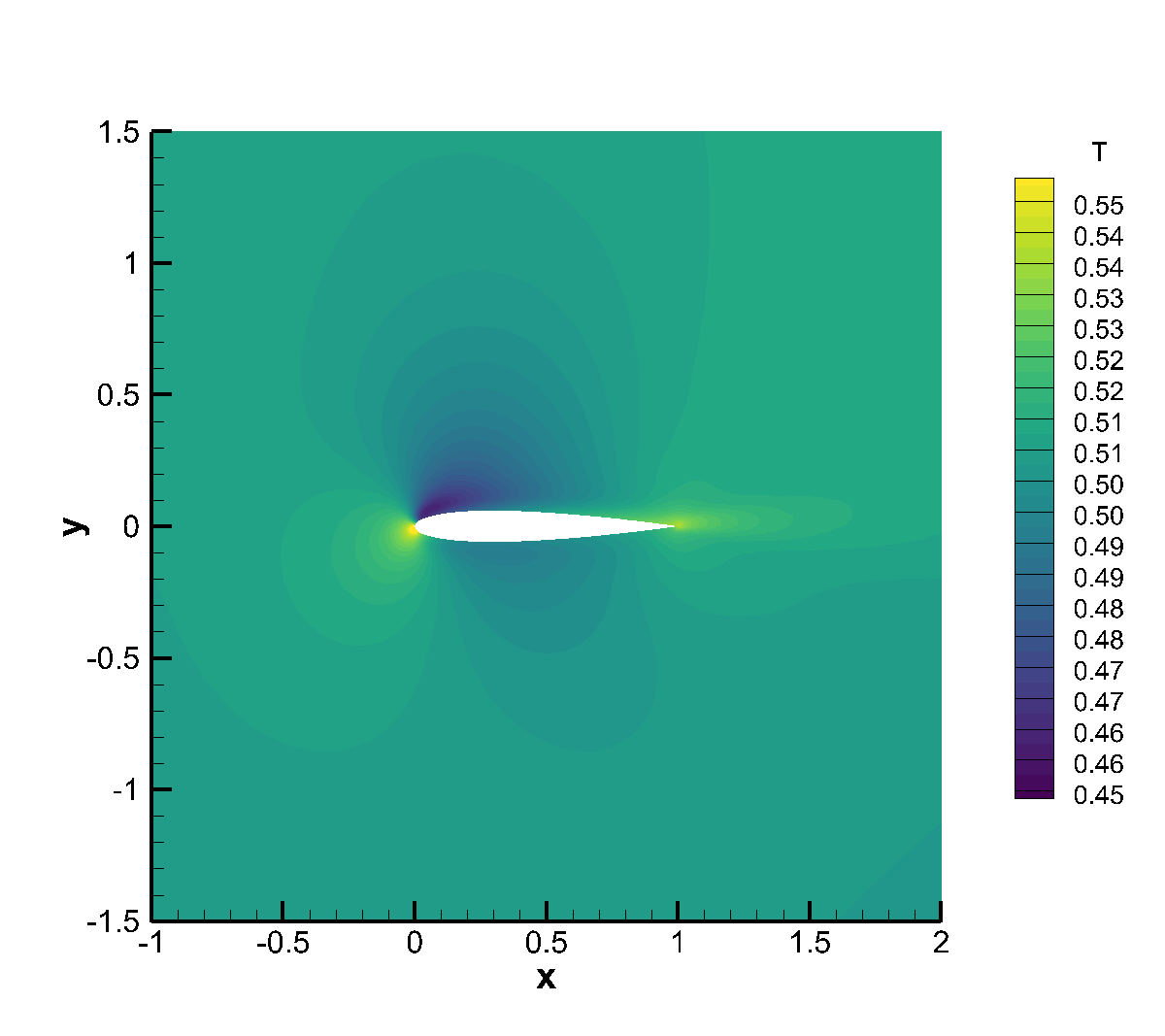}  &          
			\includegraphics[width=0.46\textwidth]{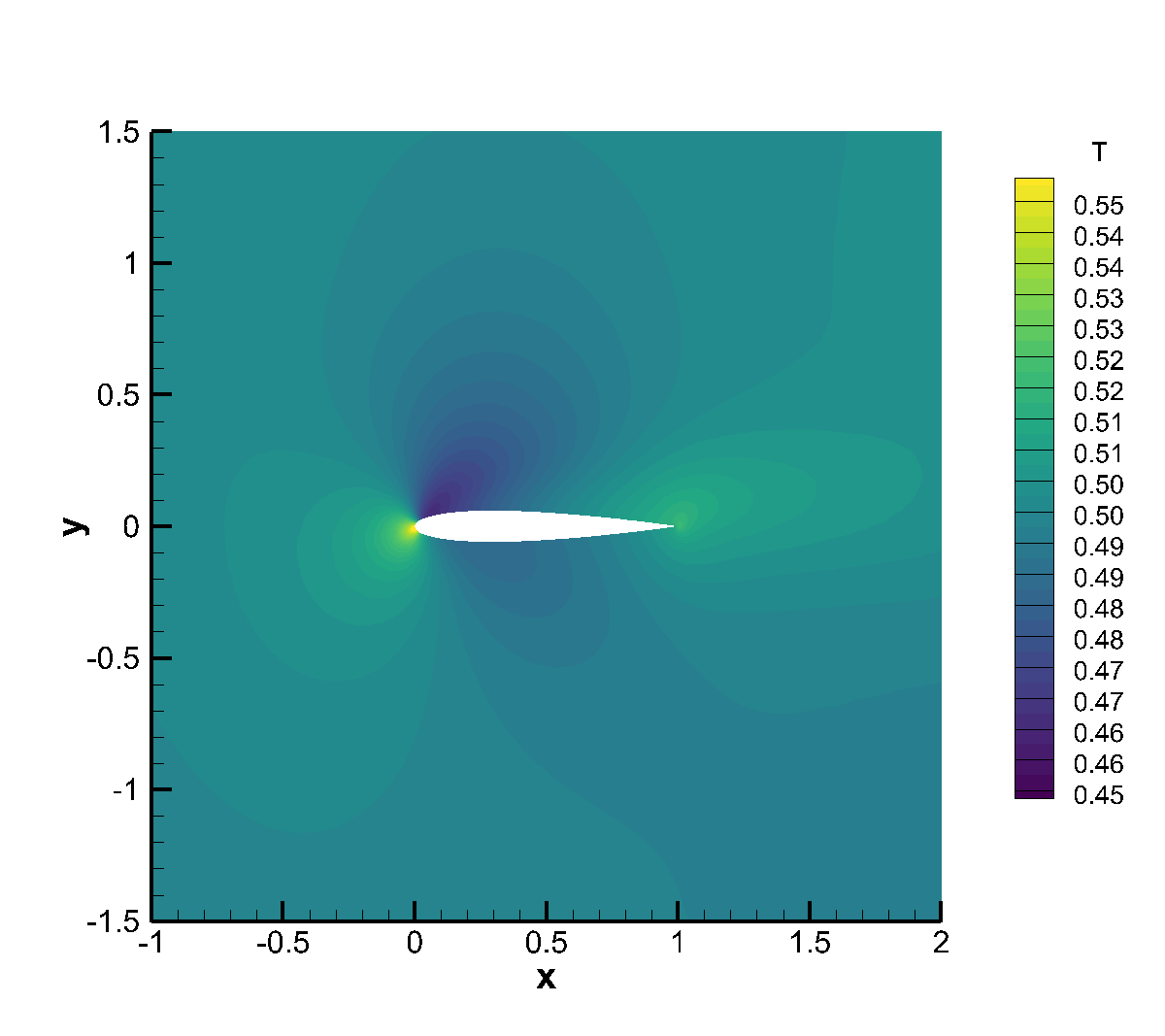} \\
			\includegraphics[width=0.46\textwidth]{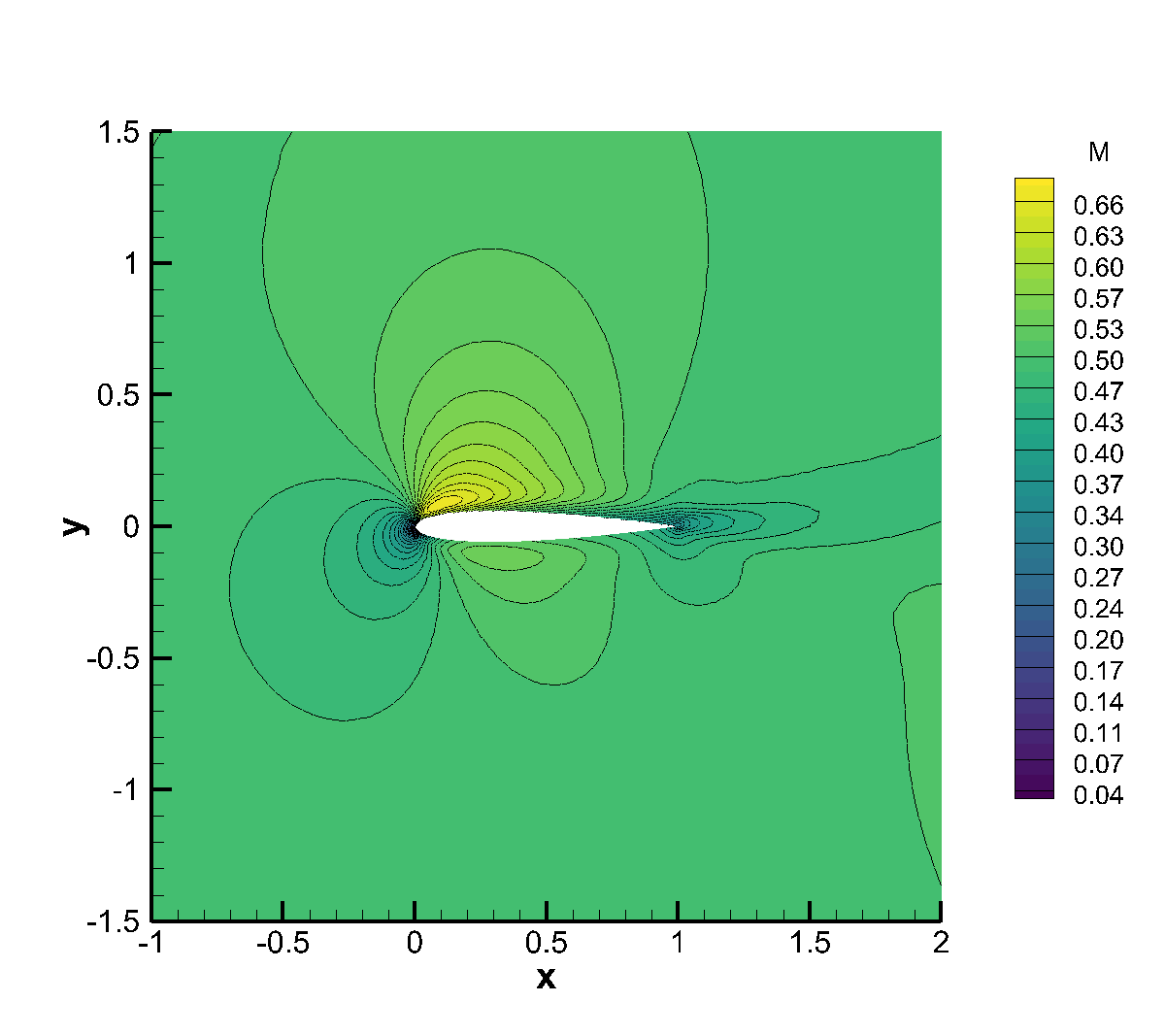}  &          
			\includegraphics[width=0.46\textwidth]{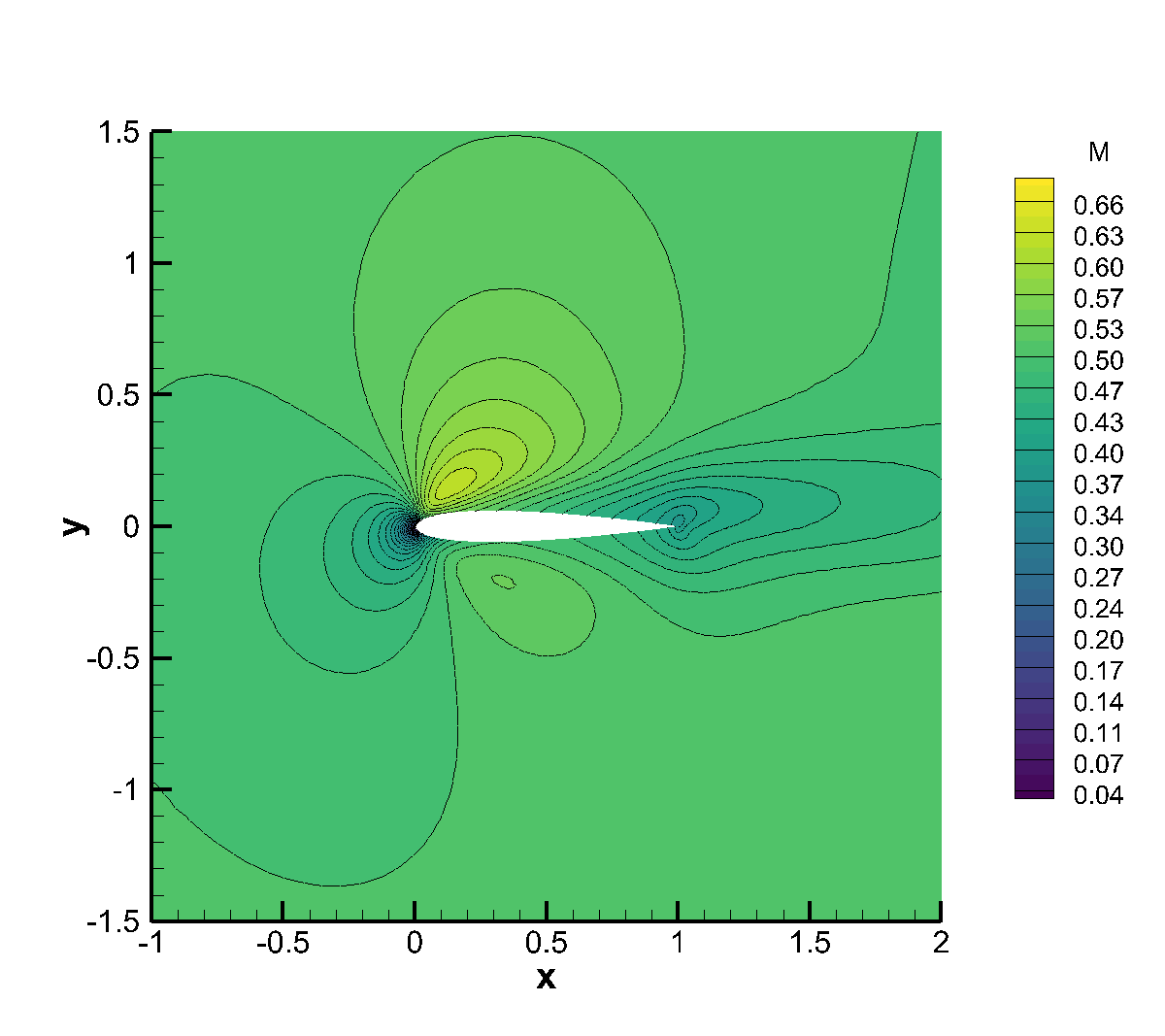} \\
			\includegraphics[width=0.46\textwidth]{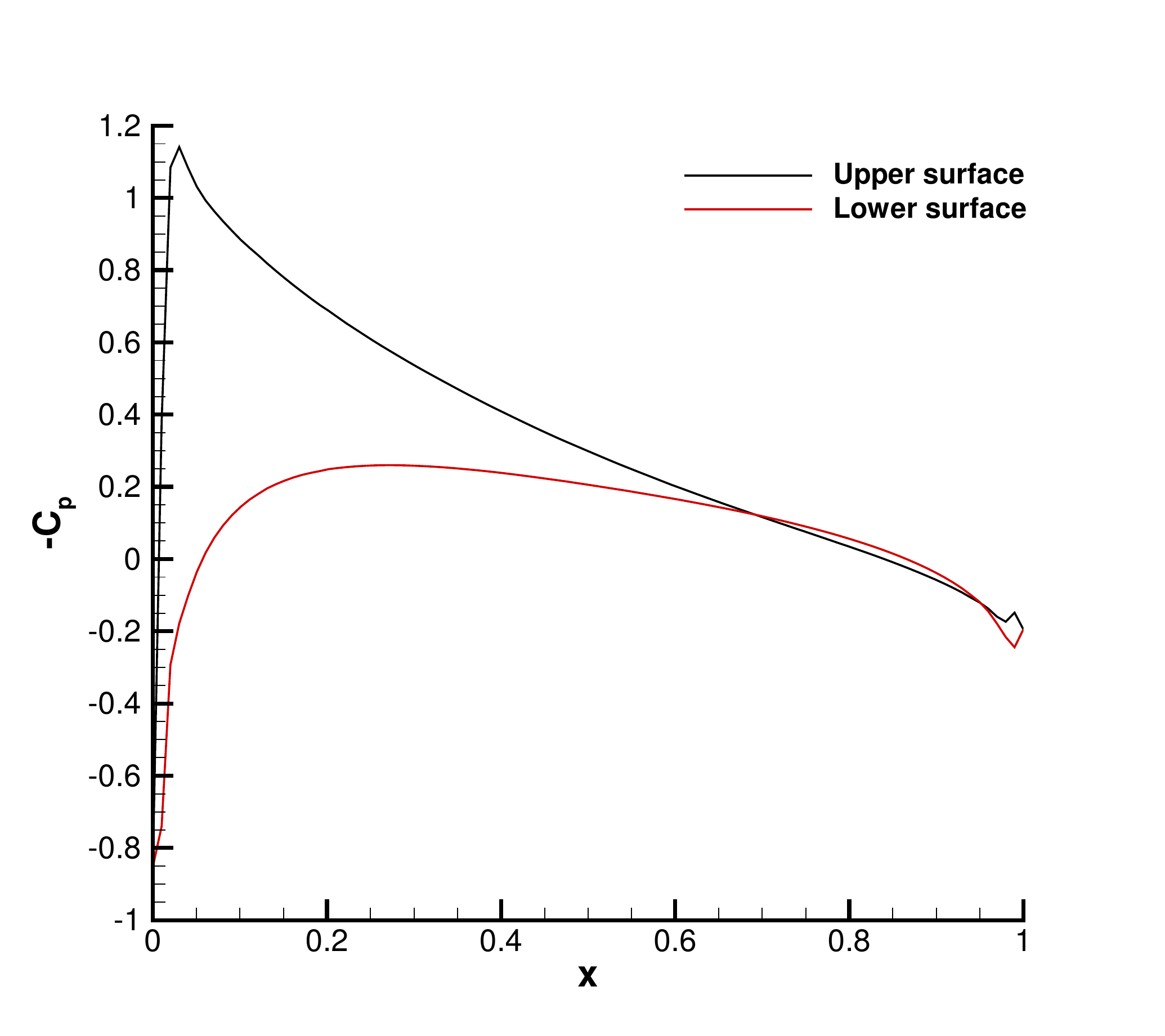}  &          
			\includegraphics[width=0.46\textwidth]{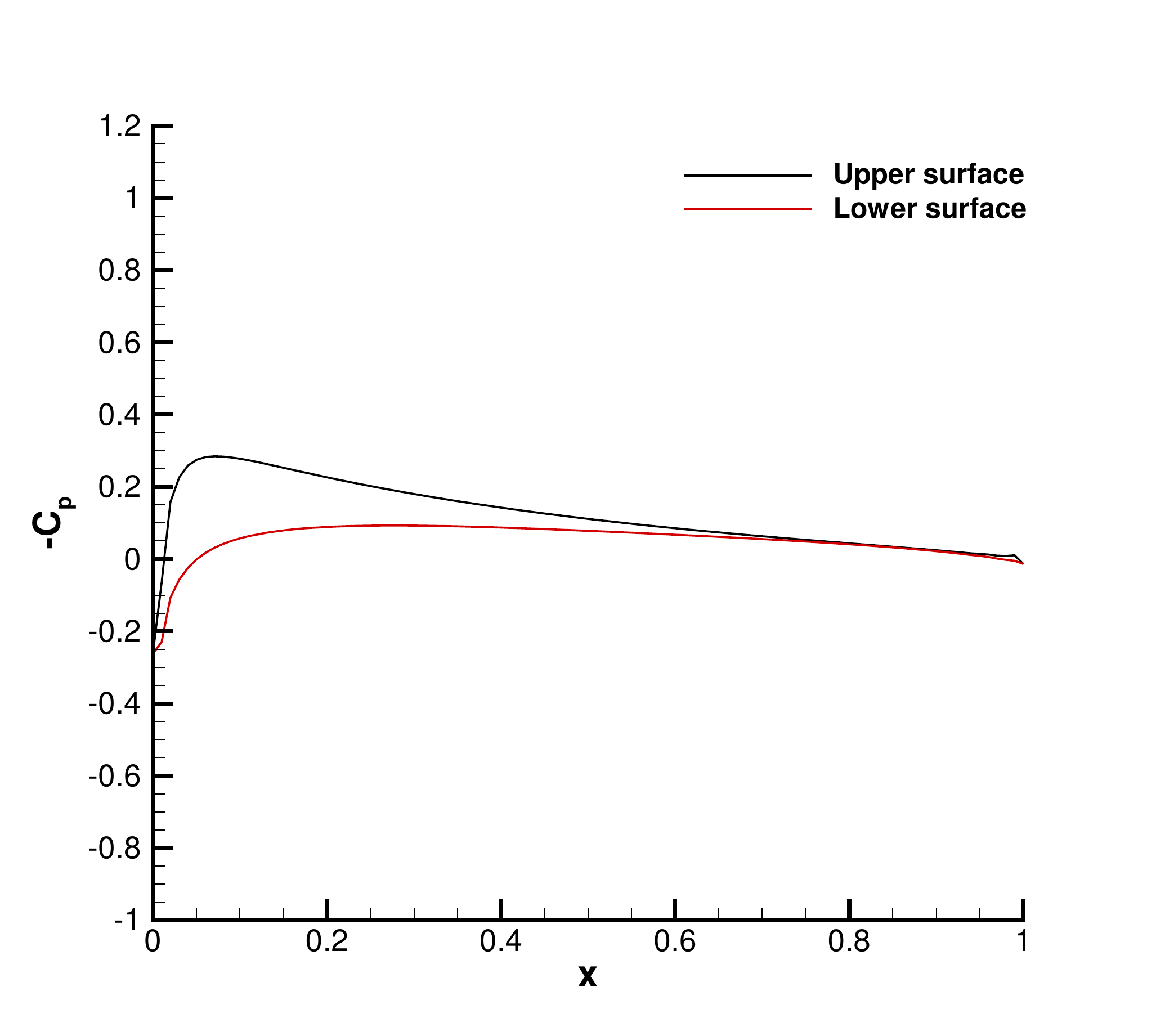} \\
		\end{tabular}
		\caption{Flow around NACA 0012 airfoil with Mach number $M=0.5$ and angle of attack $\alpha=4^{\circ}$ for $\varepsilon=5 \cdot 10^{-5}$ (left column) and $\varepsilon=5 \cdot 10^{-3}$ (right column) at $t_f=5$. Top row: 40 Mach number contours in the interval $[0.45;0.55]$. Middle row: 40 Mach number contours in the interval $[0.04;0.68]$. Bottom row: pressure coefficient distribution on the upper and lower surfaces of the airfoil. }
		\label{fig.NACA-M05-a4}
	\end{center}
\end{figure}

\begin{figure}[!htbp]
	\begin{center}
		\begin{tabular}{cc}
			\includegraphics[width=0.46\textwidth]{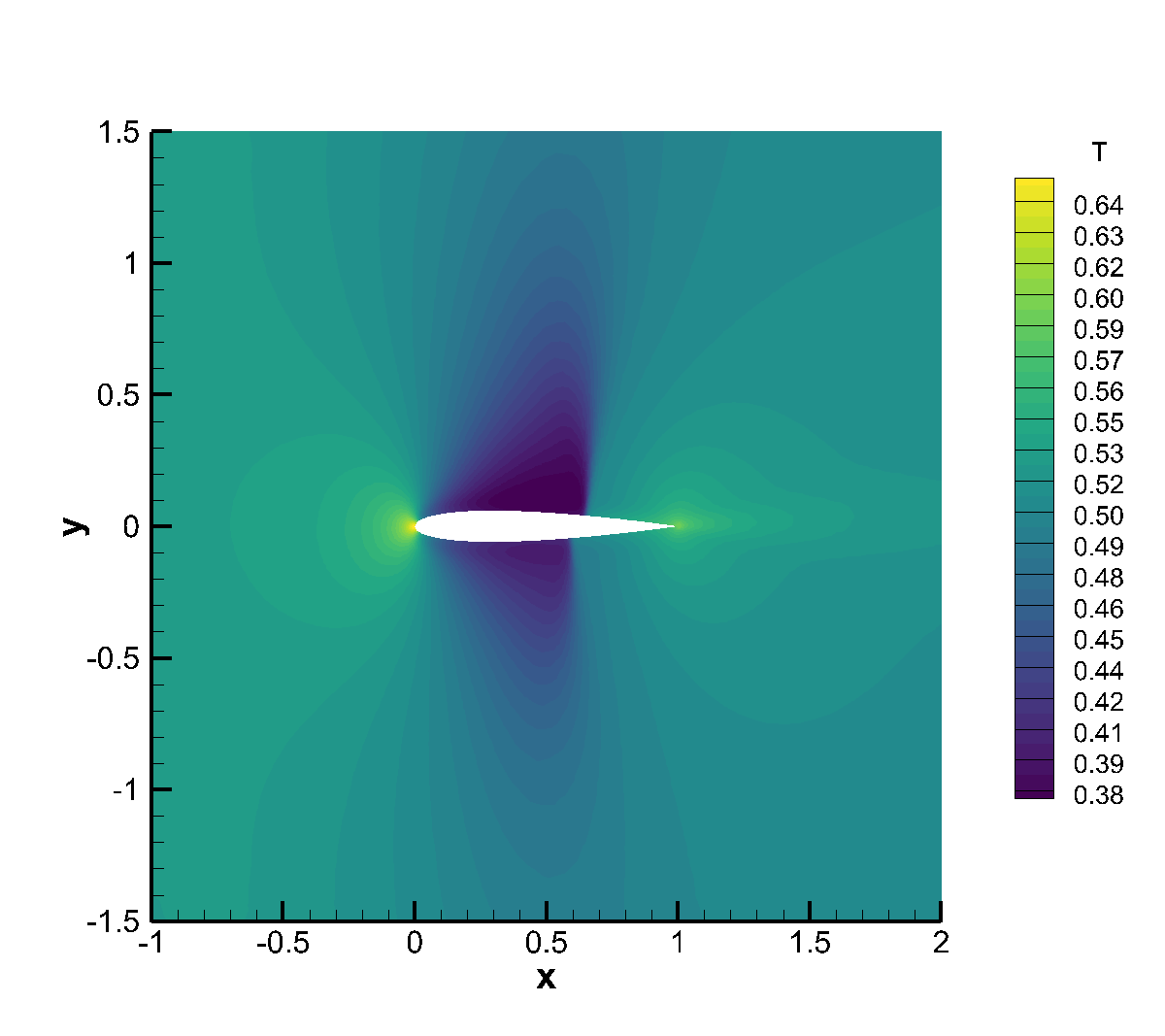}  &          
			\includegraphics[width=0.46\textwidth]{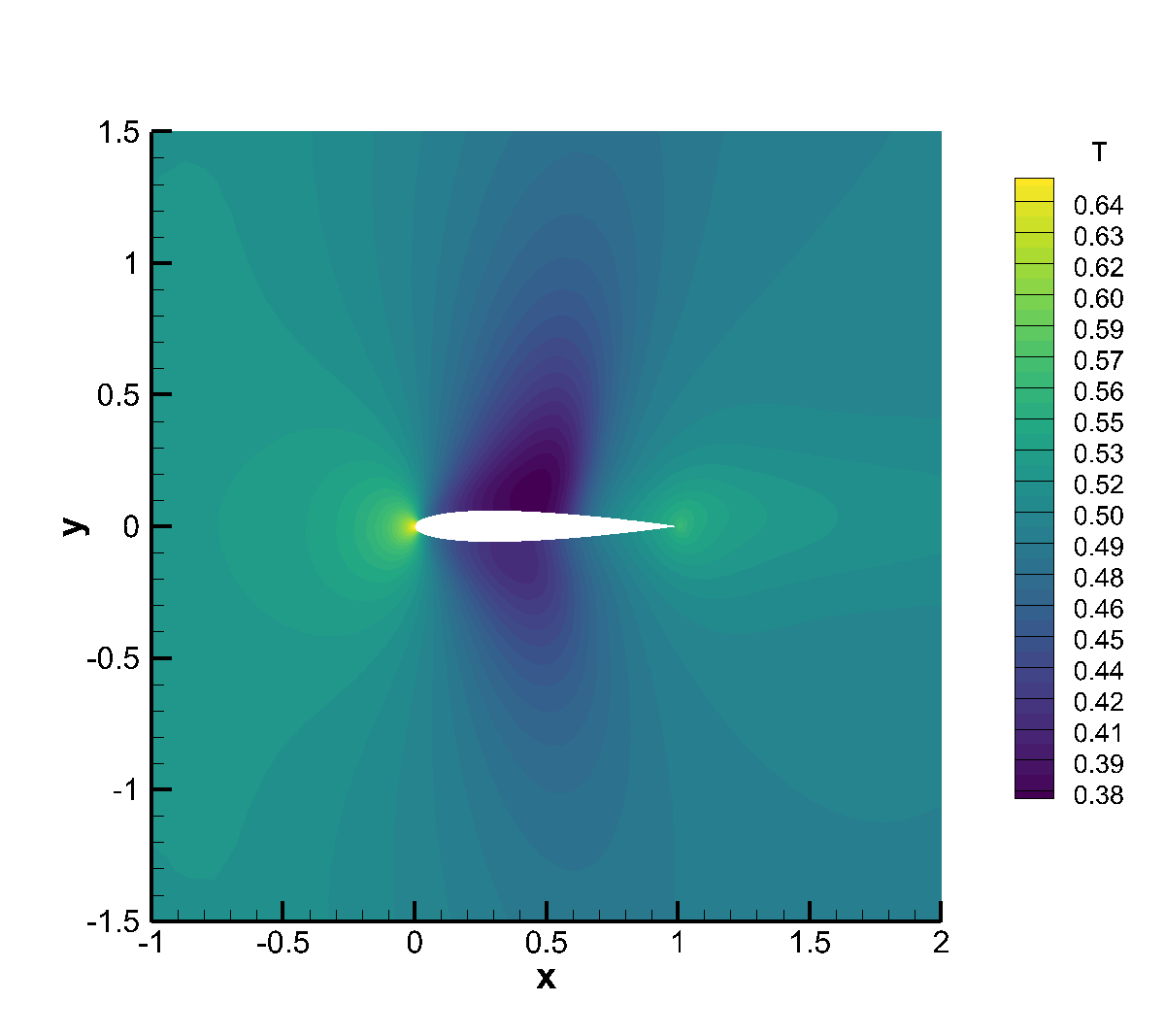} \\
			\includegraphics[width=0.46\textwidth]{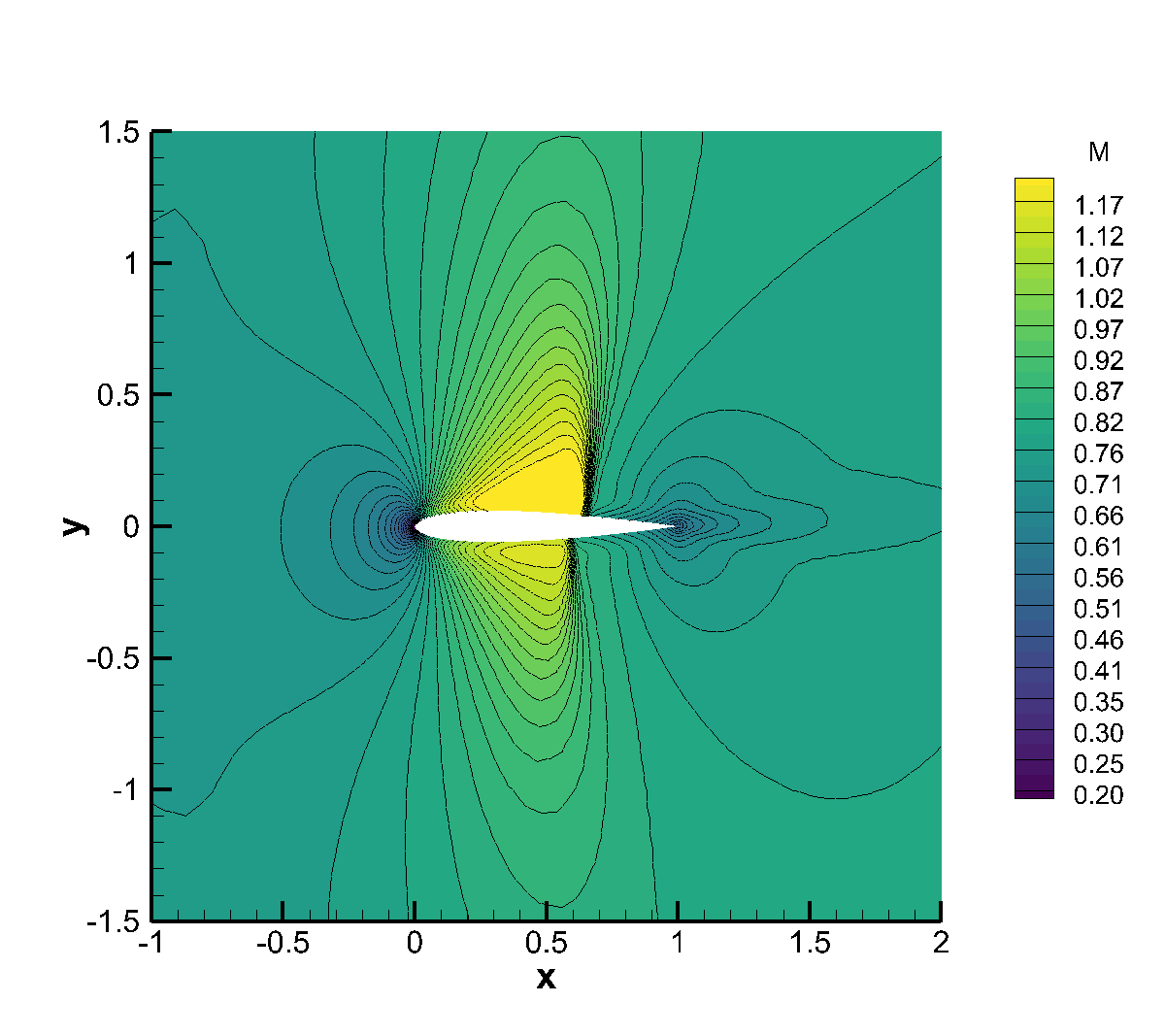}  &          
			\includegraphics[width=0.46\textwidth]{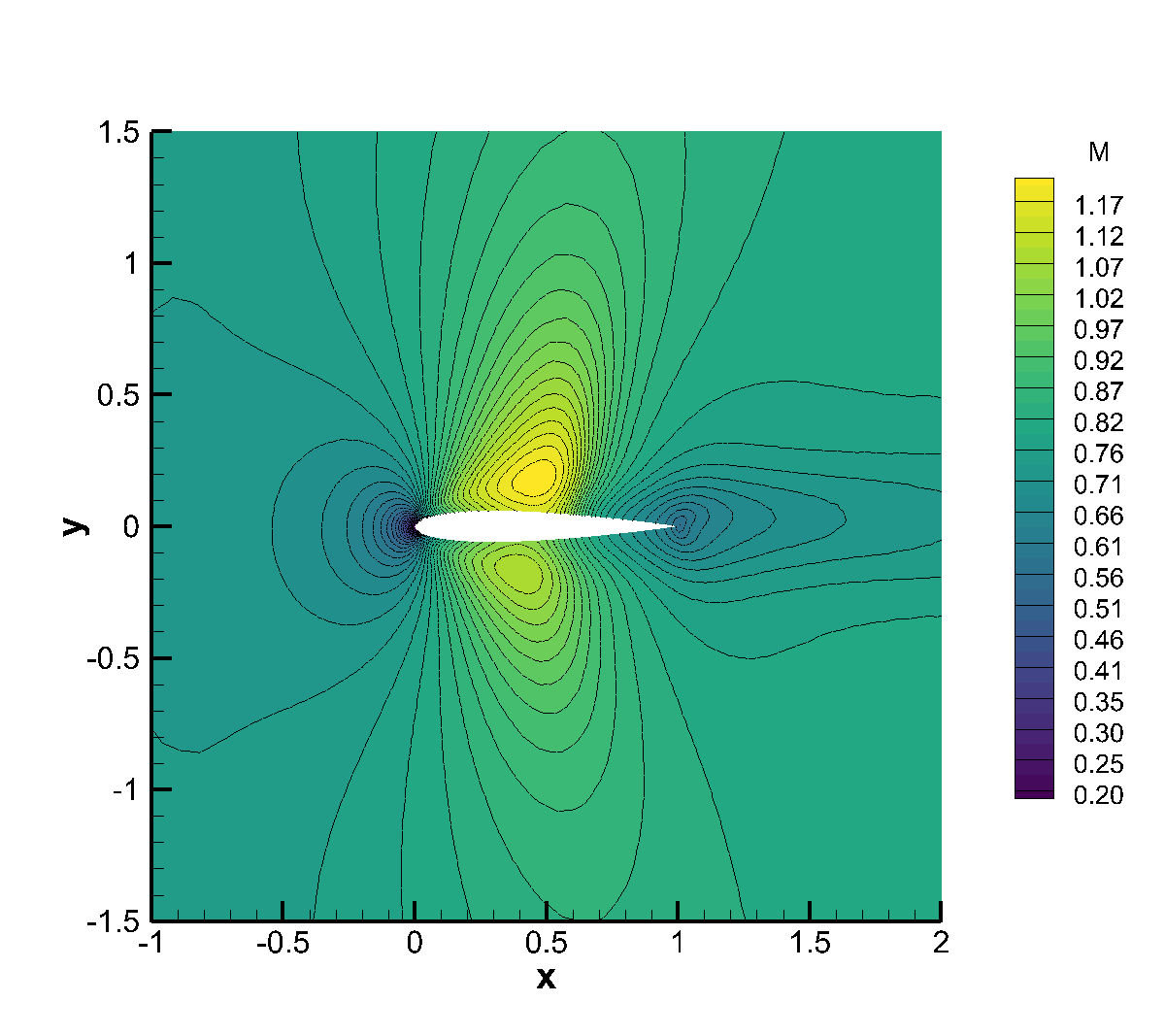} \\
			\includegraphics[width=0.46\textwidth]{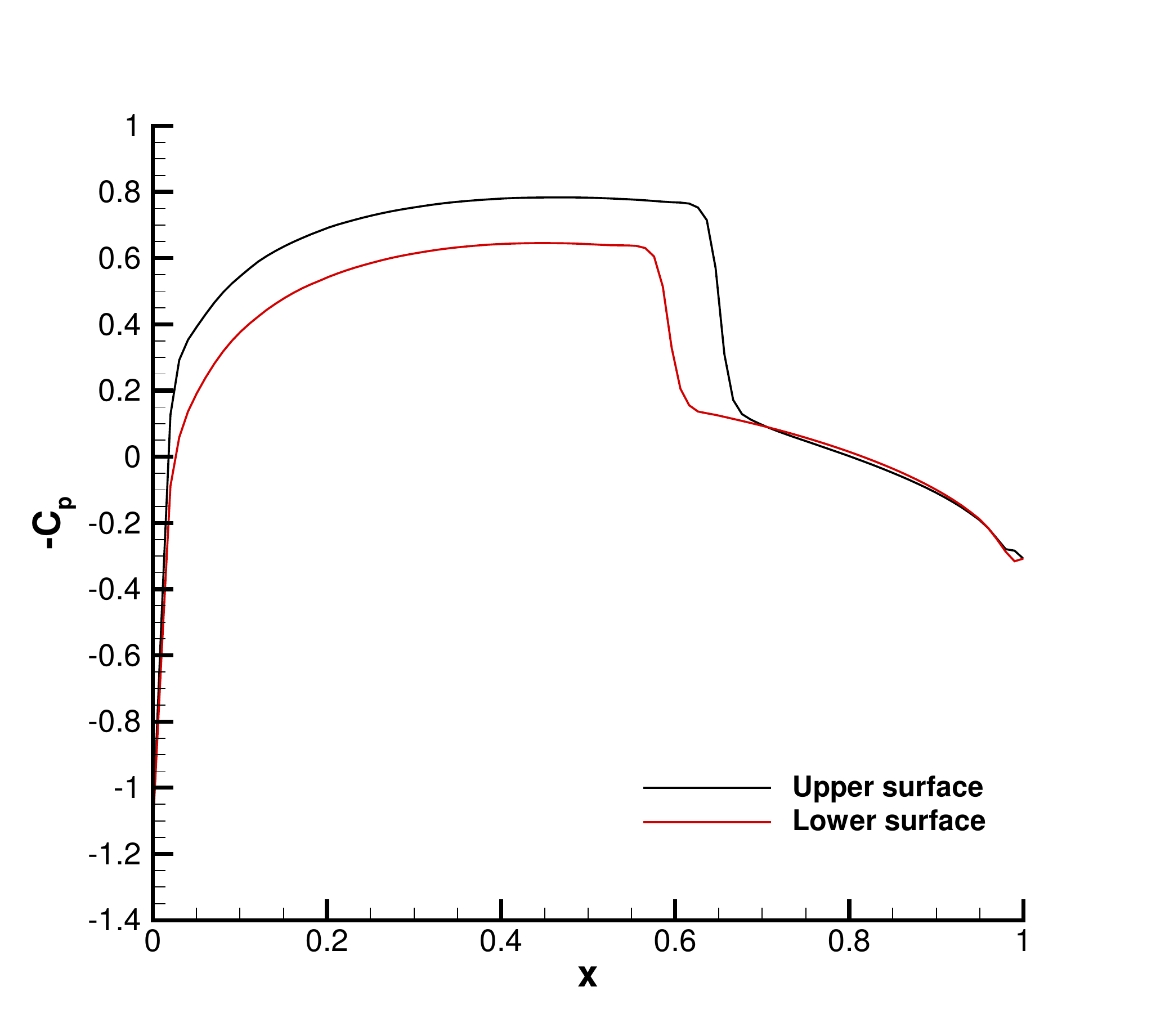}  &          
			\includegraphics[width=0.46\textwidth]{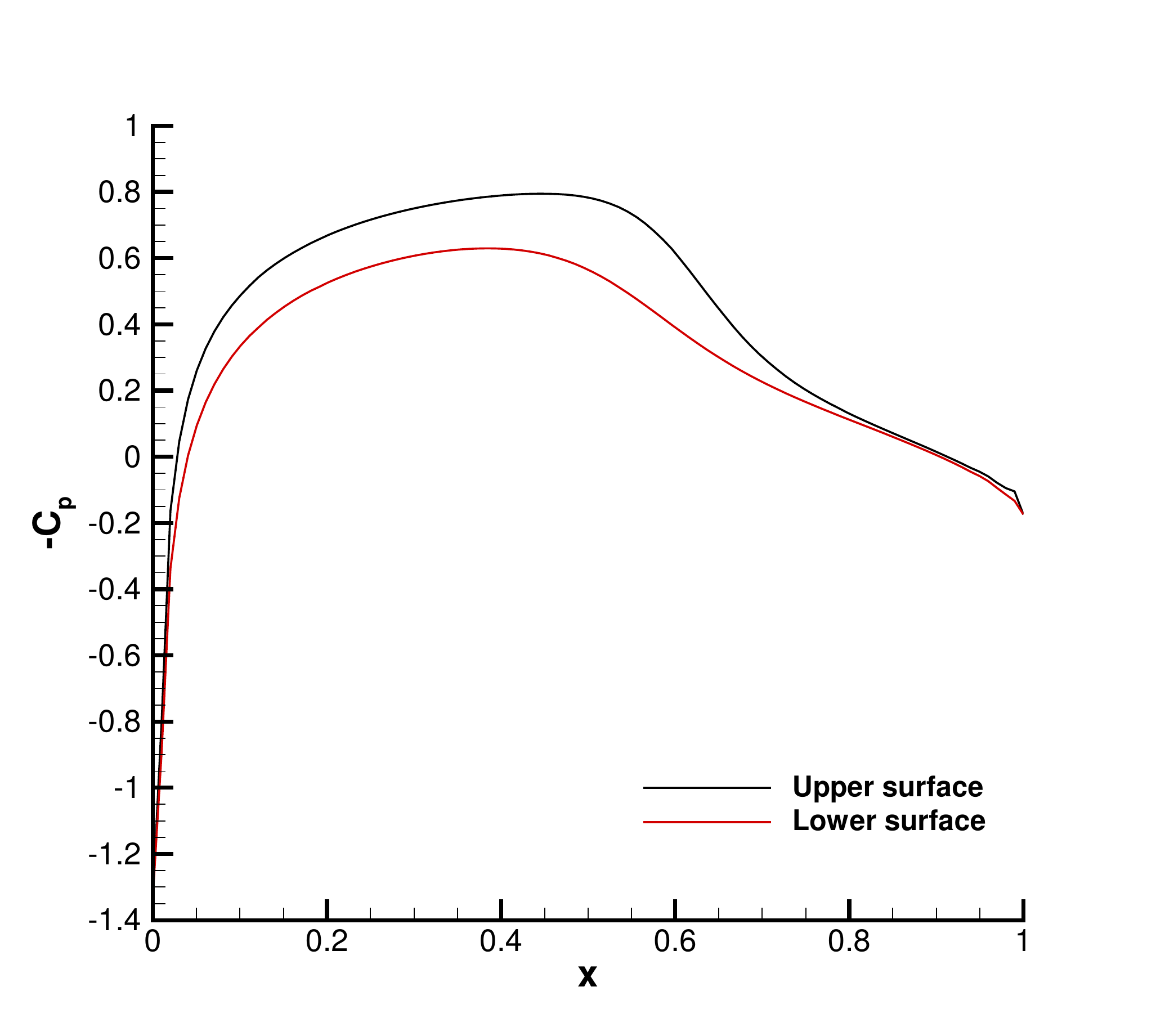} \\
		\end{tabular}
		\caption{Flow around NACA 0012 airfoil with Mach number $M=0.8$ and angle of attack $\alpha=1.25^{\circ}$ for $\varepsilon=5 \cdot 10^{-5}$ (left column) and $\varepsilon=5 \cdot 10^{-3}$ (right column) at $t_f=5$. Top row: 40 Mach number contours in the interval $[0.38;0.65]$. Middle row: 40 Mach number contours in the interval $[0.2;1.2]$. Bottom row: pressure coefficient distribution on the upper and lower surfaces of the airfoil. }
		\label{fig.NACA-M08-a125}
	\end{center}
\end{figure}

\begin{figure}[!htbp]
	\begin{center}
		\begin{tabular}{cc}
			\includegraphics[width=0.46\textwidth]{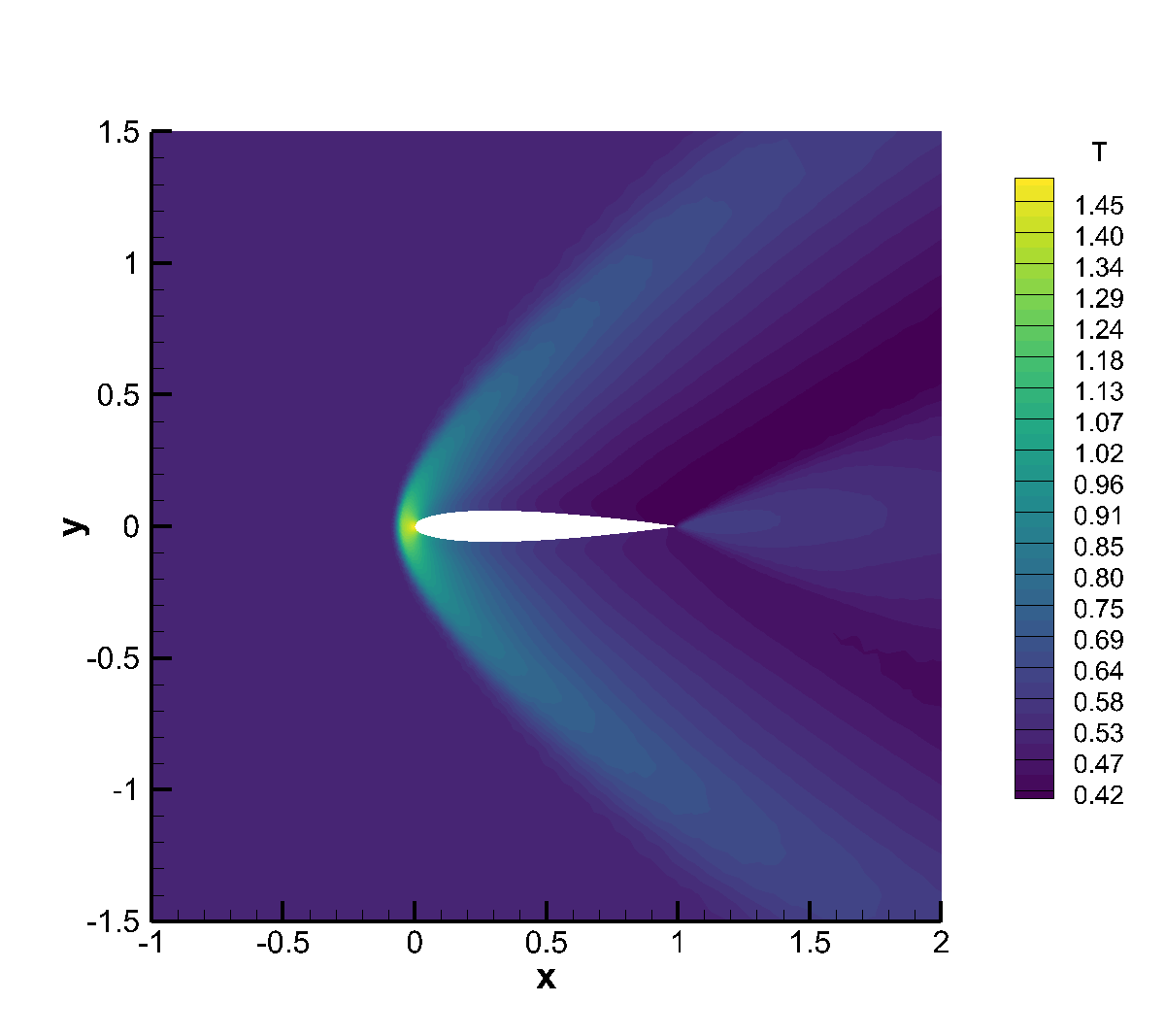}  &          
			\includegraphics[width=0.46\textwidth]{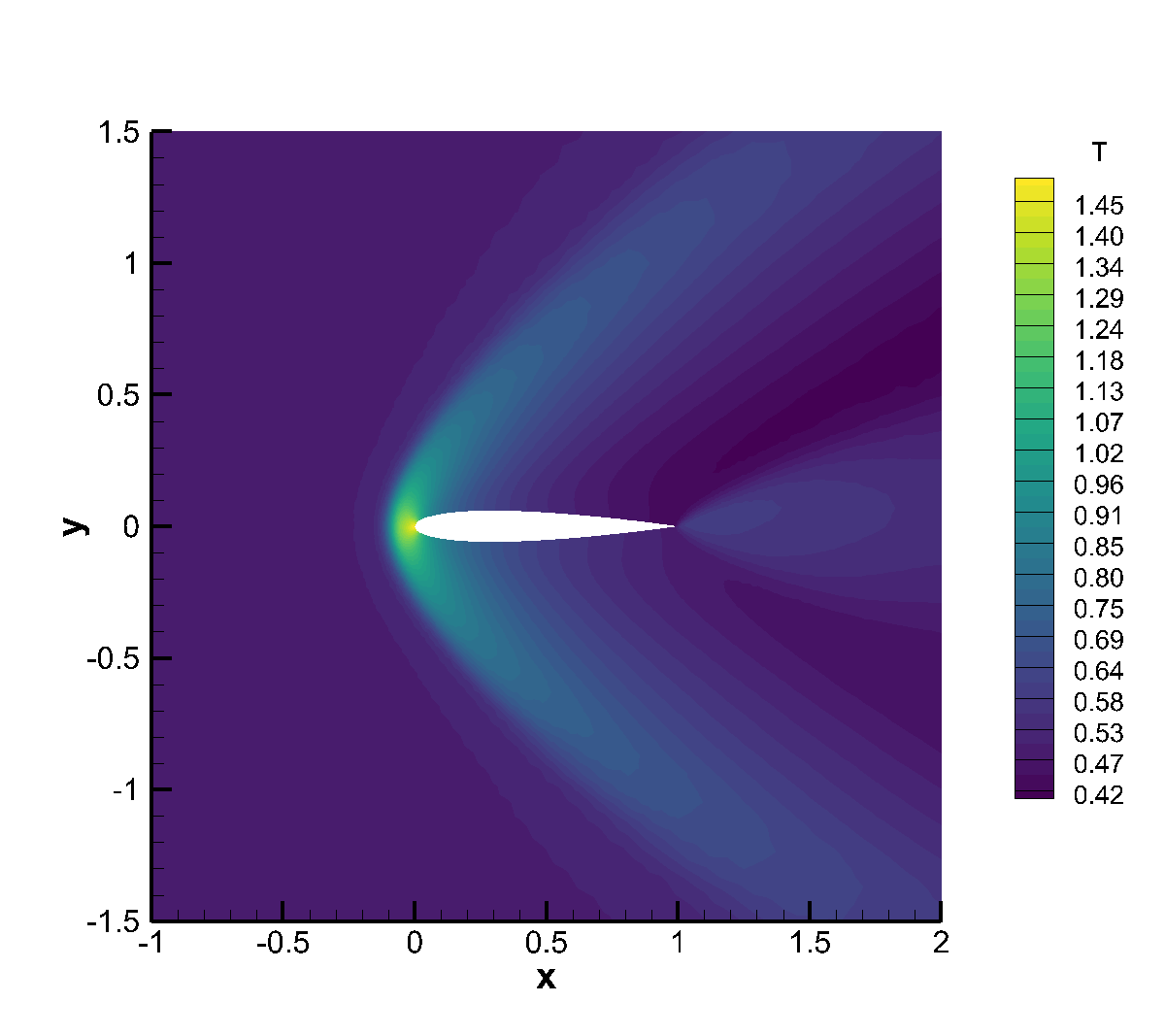} \\
			\includegraphics[width=0.46\textwidth]{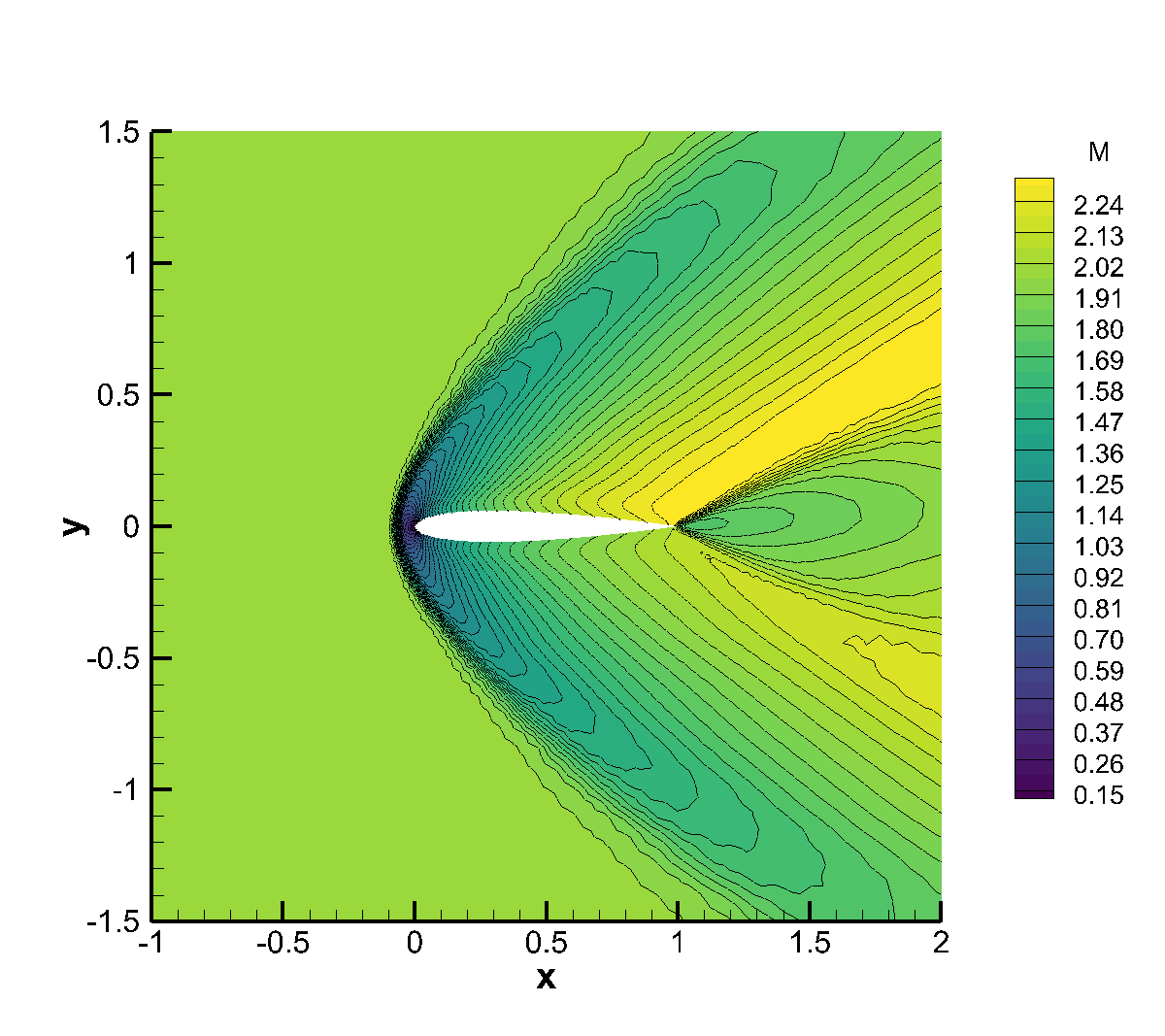}  &          
			\includegraphics[width=0.46\textwidth]{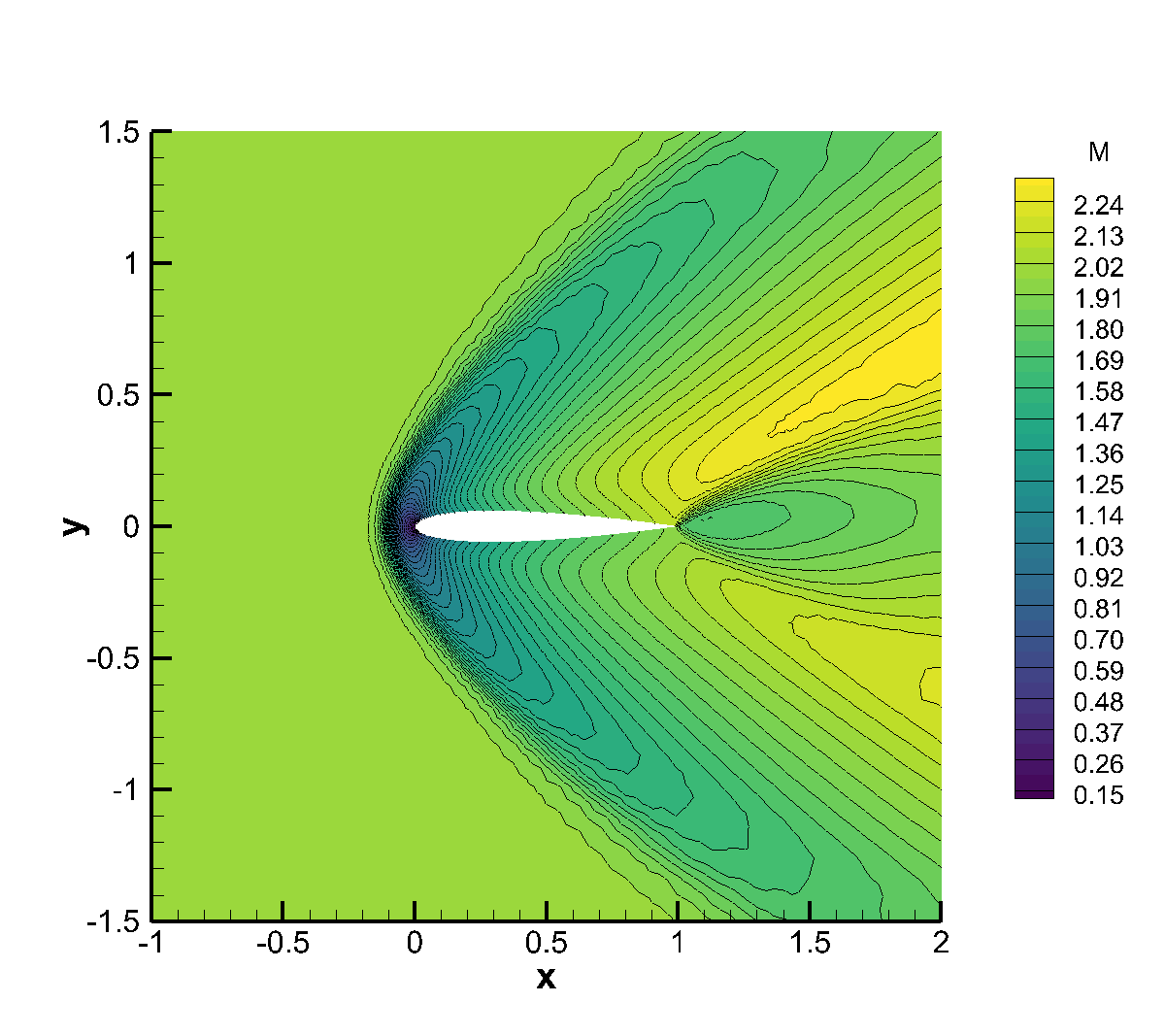} \\
			\includegraphics[width=0.46\textwidth]{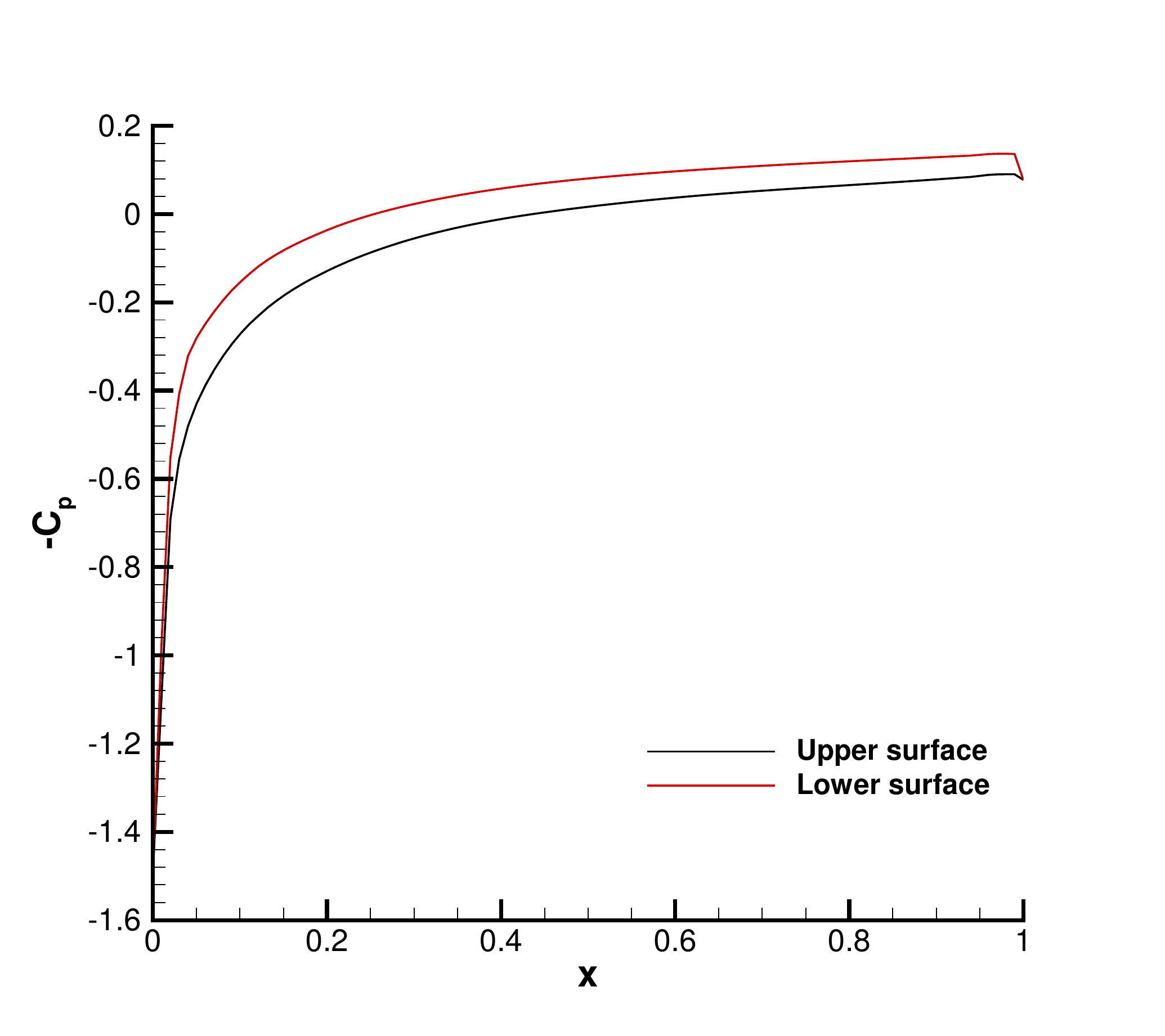}  &          
			\includegraphics[width=0.46\textwidth]{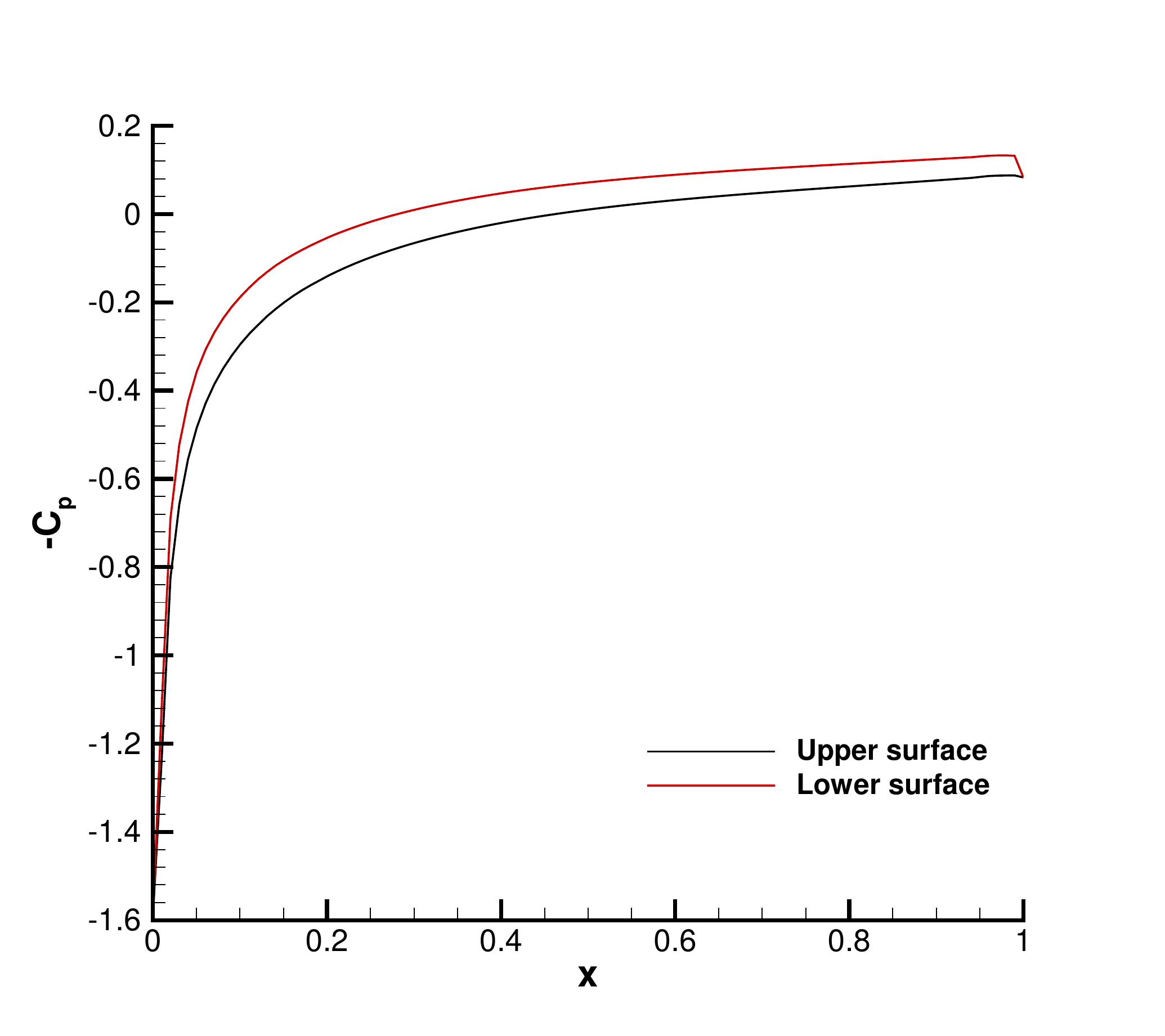} \\
		\end{tabular}
		\caption{Flow around NACA 0012 airfoil with Mach number $M=2$ and angle of attack $\alpha=2^{\circ}$ for $\varepsilon=5 \cdot 10^{-5}$ (left column) and $\varepsilon=5 \cdot 10^{-3}$ (right column) at $t_f=2$. Top row: 40 Mach number contours in the interval $[0.42;1.48]$. Middle row: 40 Mach number contours in the interval $[0.15;2.3]$. Bottom row: pressure coefficient distribution on the upper and lower surfaces of the airfoil. }
		\label{fig.NACA-M2-a2}
	\end{center}
\end{figure}

With this numerical setting, we consider then two subsonic settings and one supersonic configuration. The second order version of the novel CWENO-IMEX schemes is used to run all computations. 
The first subsonic initial condition is given by considering an inflow Mach number of $M=0.5$ and an angle of attack of $\alpha=4^{\circ}$. The free-stream conditions are as follows:
\begin{equation*}
(u,v)_\infty=(0.5,0), \qquad p_\infty=1, \qquad \rho_\infty=1,
\end{equation*}
thus the fluid pressure is set to $p=p_\infty/\gamma$. Figure \ref{fig.NACA-M05-a4} shows the numerical solution for the Mach number $M=\sqrt{u/(\gamma \, R \, T)}$ and the temperature at the final time $t_f=5$ for two different fluid regimes, namely $\varepsilon=5 \cdot 10^{-5}$ and $\varepsilon=5 \cdot 10^{-3}$. The pressure coefficient $C_p$ on the upper and the lower surface of the airfoil is also reported. It is computed as
\begin{equation}
C_p=\frac{p-p_\infty}{\frac{1}{2} \rho_\infty u_\infty^2}.
\end{equation}
In the fluid limit the results are in excellent agreement with available references in the literature \cite{Xiao17} while in the rarefied setting, the pressure coefficients become lower as expected. 
The second subsonic setting involves an inlet Mach number of $M=0.8$ and an angle of attack of $\alpha=1.25^{\circ}$, therefore a shock wave is generated both on the upper and the lower surface, as correctly shown in Figure \ref{fig.NACA-M08-a125}. This simulation is also run for two different fluid regimes and the results are in qualitatively very good agreement with \cite{Xiao17} for the case $\varepsilon=5 \cdot 10^{-5}$. The presence of the shock waves can be clearly identified also looking at the pressure coefficient distribution, while in the rarefied regime with $\varepsilon=5 \cdot 10^{-3}$ the shocks become milder and the pressure coefficients exhibit a much smoother profile.
Finally, as a last example, we consider a supersonic flow with Mach number $M=2$ up to the final time $t_f=2$. The numerical results are shown in Figure \ref{fig.NACA-M2-a2}, where it is visible the generation of shock waves both in front of the airfoil but also departing from the tail of the NACA 0012 profile. The pressure coefficients show that this simulation can be regarded as a prototype of a supersonic vehicle reentry which is approaching the soil surface. 

\section{Conclusions} \label{sec_conclu}
In this work, we have considered a high order in velocity, space and time finite volume method for solving the Boltzmann equation. The proposed approach is based on a spectral discretization of the collision operator, on a Central WENO polynomial reconstruction and on a high order implicit-explicit time discretization. The scheme is designed to work on arbitrarily unstructured control volumes. Up to our knowledge, this is one of the first examples in which the full Boltzmann model is solved using high order time-space methods on such unstructured meshes. In particular, we are not aware of any result in the case in which, in the above described framework, high order implicit methods which guarantee stability, accuracy and preservation of the asymptotic state are used in a multidimensional setting for this model.   

The numerical part consists of several measures of the theoretical convergence order for the different employed discretizations (spectral, CWENO and Runge-Kutta) and for different regimes: from dense fluids to rarefied gases. An additional measure of the computational costs involved is also reported. A second part consists in proving the capability of the designed method to deal with standard benchmark problems in gas dynamics. Here, we compare our scheme with the limit compressible Euler equations and with a Direct Simulation Monte Carlo method solving the Boltzmann equation. Comparisons between the Boltzmann and the BGK model are also proposed showing the ineffectiveness of the latter in describing far from equilibrium phenomena. Finally, a more realistic set of simulations on a NACA wing profile under different regimes and angles of attack is presented. These tests run on a MPI parallelization version of the code and provide evidence that the method can be applied to large scale simulations. A very good agreement with the existing literature is observed for these last cases.

In the future, the methods here presented will be extended to the full six dimensional case. We also aim at considering the Discontinuous Galerkin (DG) method for the space discretization and to increase the time precision and lowering the computational effort to use linear multistep approaches in time. Finally, a way to better represent the physical solutions would be to change the support of the distribution function and the nodes where it is defined as a function of the solution itself. This is a challenging research direction that we would aim to follow in the next future.

\section*{Acknowledgments}
The support of the Italian Ministry of Instruction, University and Research (MIUR) is acknowledged within PRIN Project 2017 No. 2017KKJP4X entitled "Innovative numerical methods for evolutionary partial differential equations and applications". The authors are members of the GNCS-INDAM (\textit{Istituto Nazionale di Alta Matematica}) group.

\bibliographystyle{plain}
\bibliography{biblio}

\end{document}